\documentclass[final,leqno,onefignum,onetabnum]{siamltex1213}
\usepackage{amsfonts}
\usepackage{amsmath,booktabs,ctable,threeparttable}
\usepackage{amssymb,amsfonts,boxedminipage}
\newtheorem{remark}{Remark}[section]

\title{The Low Rank Approximations and Ritz Values
in LSQR For Linear Discrete Ill-Posed Problems\thanks{This
work was supported in part by
the National Science Foundation of China (No. 11771249).}}

\author{Zhongxiao Jia\thanks{Department of Mathematical Sciences, Tsinghua
University, 100084 Beijing, China. (\email{jiazx@tsinghua.edu.cn})}}

\begin{document}
\maketitle
\slugger{sirev}{xxxx}{xx}{x}{x--x}

\begin{abstract}
LSQR and its mathematically equivalent CGLS have been popularly used over the
decades for large-scale linear discrete ill-posed problems, where
the iteration number $k$ plays the role of the regularization parameter.
It has been long known that if the Ritz values in LSQR converge
to the large singular values of $A$ in natural order until its semi-convergence
then LSQR must have the same the regularization
ability as the truncated
singular value decomposition (TSVD) method and can compute a 2-norm filtering
best possible regularized solution. However, hitherto there has been no definitive
rigorous result on the approximation behavior of the Ritz values in the context
of ill-posed problems. In this paper, for severely, moderately and mildly
ill-posed problems, we give accurate solutions of the two closely related
fundamental and highly challenging problems on the regularization of LSQR:
(i) How accurate are
the low rank approximations generated by Lanczos bidiagonalization?
(ii) Whether or not the Ritz values involved in LSQR approximate
the large singular values of $A$ in natural order? We also
show how to judge the accuracy of low rank approximations
reliably during computation without extra cost.
Numerical experiments confirm our results.
\end{abstract}


\pagestyle{myheadings}
\thispagestyle{plain}
\markboth{ZHONGXIAO JIA}{LOW RANK APPROXIMATIONS AND RITZ VALUES}

\maketitle

\section{Introduction}\label{intro}

Consider the linear discrete ill-posed problem
\begin{equation}
  \min\limits_{x\in \mathbb{R}^{n}}\|Ax-b\| \mbox{\,\ or \ $Ax=b$,}
  \ \ \ A\in \mathbb{R}^{m\times n}, \label{eq1}
  \ b\in \mathbb{R}^{m},
\end{equation}
where the norm $\|\cdot\|$ is the 2-norm of a vector or matrix, and
$A$ is extremely ill conditioned with its singular values decaying
and centered at zero without a noticeable gap,
and the right-hand side $b=b_{true}+e$ is assumed to be
contaminated by a Gaussian white noise $e$, where $b_{true}$
is noise-free and $\|e\|<\|b_{true}\|$. Without loss of generality, we
assume $m\geq n$ since the results in this paper hold for
the $m\leq n$ case. \eqref{eq1} arises from many
applications, such as image deblurring, signal processing, geophysics,
computerized tomography, heat propagation, biomedical and optical imaging,
and groundwater modeling, to name a few; see, e.g.,
\cite{aster,engl93,engl00,ito15,kaipio,kern,kirsch,natterer,vogel02}.
Due to the noise $e$ and the high
ill-conditioning of $A$, the naive solution $x_{naive}=A^{\dagger}b$
of \eqref{eq1} is generally a meaningless approximation to
the true solution $x_{true}=A^{\dagger}b_{true}$, where
$\dagger$ denotes the Moore-Penrose inverse of a matrix.
Therefore, regularization must be used to extract a
good approximation to $x_{true}$.

For a Gaussian white noise $e$,
throughout the paper, we always assume that $b_{true}$ satisfies the discrete
Picard condition $\|A^{\dagger}b_{true}\|\leq C$ with some constant $C$ for
$\|A^{\dagger}\|$ arbitrarily
large \cite{aster,gazzola15,hansen90,hansen90b,hansen98,hansen10,kern}.
Assume that $Ax_{true}=b_{true}$. Then a dominating regularization approach is
to solve the following problem
\begin{equation}\label{posed}
\min\limits_{x\in \mathbb{R}^{n}}\|Lx\| \ \ \mbox{subject to}\ \ \|Ax-b\|\leq
\tau \|e\|
\end{equation}
with $\tau>1$ slightly \cite{hansen98,hansen10},
where $L$ is a regularization matrix and its suitable choice is based on
a-prior information on $x_{true}$.
If $L\not=I$, \eqref{posed} can be mathematically transformed into a standard-form
problem \cite{hansen98,hansen10}, i.e., a 2-norm filtering regularization problem.
In this paper, we always take $L=I$.

The solutions of \eqref{eq1} and \eqref{posed} can be
analyzed by the means of the singular value decomposition (SVD) of $A$:
\begin{equation}\label{eqsvd}
  A=U\left(\begin{array}{c} \Sigma \\ \mathbf{0} \end{array}\right) V^{T},
\end{equation}
where $U = (u_1,u_2,\ldots,u_m)\in\mathbb{R}^{m\times m}$ and
$V = (v_1,v_2,\ldots,v_n)\in\mathbb{R}^{n\times n}$ are orthogonal,
$\Sigma = {\rm diag} (\sigma_1,\sigma_2,\ldots,\sigma_n)\in\mathbb{R}^{n\times n}$
with the singular values
$\sigma_1>\sigma_2 >\cdots >\sigma_n>0$ assumed to be simple
throughout the paper, and the superscript $T$
denotes the transpose of a matrix or vector.

From the SVD expansion $x_{true}=\sum\limits_{i=1}^{n}\frac{u_i^{T}b_{true}}
{\sigma_i}{v_i}$, the discrete Picard condition means that,
on average, the Fourier coefficient
$|u_i^{T}b_{true}|$ decays faster than $\sigma_i$, and it enables
regularization to compute useful approximations to $x_{true}$. The
following common model is used throughout Hansen's books
\cite{hansen98,hansen10} and the references therein as well as
\cite{jia18a} and the current paper:
\begin{equation}\label{picard}
  | u_i^T b_{true}|=\sigma_i^{1+\beta},\ \ \beta>0,\ i=1,2,\ldots,n.
\end{equation}

For the Gaussian white noise $e$, its covariance matrix
is $\eta^2 I$, the expected values $\mathcal{E}(\|e\|^2)=m \eta^2$ and
$\mathcal{E}(|u_i^Te|)=\eta,\,i=1,2,\ldots,n$, so that
$\|e\|\approx \sqrt{m}\eta$ and $|u_i^Te|\approx \eta,\
i=1,2,\ldots,n$; see, e.g., \cite[p.70-1]{hansen98} and \cite[p.41-2]{hansen10}.
Under the condition \eqref{picard}, for large singular values,
the signal term $|{u_i^{T}b_{true}}|/{\sigma_i}$ is dominant relative to
the noise term $|u_i^{T}e|/{\sigma_i}$. Once
$| u_i^T b_{true}| \leq | u_i^T e|$ from some $i$ onwards, the noise
$e$ dominates $| u_i^T b|$
for small singular values and must be suppressed.
The number $k_0$ satisfying
\begin{equation}\label{picard1}
| u_{k_0}^T b|\approx | u_{k_0}^T b_{true}|> | u_{k_0}^T e|\approx
\eta, \ | u_{k_0+1}^T b|
\approx | u_{k_0+1}^Te|
\approx \eta
\end{equation}
is called the transition point; see \cite[p.70-1]{hansen98} and
\cite[p.42, 98]{hansen10}.

The truncated SVD (TSVD) method \cite{hansen90,hansen98,hansen10} is a reliable
and commonly used method for solving a small or medium sized
\eqref{posed}. It solves a sequence of problems
\begin{equation}\label{tsvd}
\min\|x\| \ \ \mbox{subject to}\ \
 \|A_kx-b\|=\min
\end{equation}
starting with $k=1$ onwards, where $A_k=U_k\Sigma_k V_k^T$
is the 2-norm best rank $k$ approximation to $A$
with $U_k=(u_1,\ldots,u_k)$, $V_k=(v_1,\ldots,v_k)$ and $\Sigma_k=
{\rm diag}(\sigma_1,\ldots,\sigma_k)$, and
$\|A-A_k\|=\sigma_{k+1}$ \cite[p.12]{bjorck96}. The solution
$
x_{k}^{tsvd}=A_k^{\dagger}b
$
to \eqref{tsvd} is called the TSVD regularized solution,
and the index $k$ plays the role of the regularization parameter.

For the Gaussian white noise $e$, it follows
from \cite[p.71,86-8,95]{hansen10} that
$x_{k_0}^{tsvd}$ is the best TSVD solution. Moreover, it is
known from \cite{engl00,hansen98,hansen10,vogel02}
that $x_{k_0}^{tsvd}$ is a 2-norm filtering best possible regularized
solution of \eqref{eq1} when only deterministic
2-norm filtering regularization methods are taken into account.
As a result, we can take $x_{k_0}^{tsvd}$ as the standard reference
when assessing the regularization ability of a deterministic
2-norm filtering regularization method; for more general
elaborations, see \cite{jia18a}.

Over the decades Krylov solvers have been popularly used
to solve a large \eqref{eq1}. The methods project \eqref{eq1}
onto a sequence of low dimensional Krylov subspaces
and computes iterates from the subspaces to approximate $x_{true}$
\cite{aster,engl00,gilyazov,hanke95,hansen98,hansen10,kirsch}.
Of them, the CGLS method \cite{bjorck96}, which implicitly applies the CG
method to $A^TAx=A^Tb$,
and its mathematically equivalent LSQR algorithm \cite{paige82}
have been most commonly used. They are 2-norm filtering regularization
methods, have general regularizing effects
\cite{aster,eicke,gilyazov,hanke95,hanke01,hansen98,hansen10,hps16,hps09},
and exhibit typical semi-convergence \cite[p.89]{natterer}: the iterates
converge to $x_{true}$ in an initial stage; then the
noise $e$ starts to deteriorate the iterates so that they start to diverge
from $x_{true}$ and instead converge to $x_{naive}$;
see also \cite[p.314]{bjorck96}, \cite[p.733]{bjorck15},
\cite[p.135]{hansen98} and \cite[p.110]{hansen10}.

It is important to stress two special practical cases.
First, if the noise $e$ is so small that
all the $|u_i^Tb|\approx |u_i^T b_{true}|>\eta$,
then $k_0=n$ in \eqref{picard1},
meaning that $x_n^{tsvd}=x_{naive}$ is the best approximation to
$x_{true}$. Second, if $e$ is such
that all the $|u_i^T e|\approx \eta<\sigma_n$, that is, the
noise level $\|e\|$ is small relative to $\sigma_n$, then
the noise amplification 
is tolerable
even without regularization and $x_{naive}$ is a good approximation
to $x_{true}$, as has been noticed in \cite[p.7]{vogel02}. Both cases
show that for a given $e$, if $A$ is not ill conditioned enough,
regularization does not play a role in the solution process.

Indeed, we have encountered such practical
image deblurring problems resulting from two or three dimensional
continuous ill-posed problems, e.g., the image
deblurring problems {\sf fanbeamtomo} of $m=61200,\
n=14400$ with the relative noise level $\|e\|/\|b_{true}\|\leq 10^{-3}$
\cite{berisha}, {\sf blur} \cite{hansen07},
{\sf paralleltomo} \cite{hansen12} of order over ten thousands, three
{\sf GaussianBlur4XX} of $m=n=65536$ \cite{berisha}. We have found that
the singular values $\sigma_j$ of these matrices, on average,
decay more slowly than $\mathcal{O}(j^{-\alpha})$ with
$\alpha<\frac{1}{2}$ and the ratios $\sigma_1/\sigma_n$
are quite modest, i.e., $\mathcal{O}(10)\sim \mathcal{O}(10^3)$.
In view of conditioning, such problems with noise-free right-hand
sides $b_{true}$ seem to be well conditioned ordinary least
squares problems or linear systems. We have observed that
$\|e\|/\|b_{true}\|\leq 10^{-3}$ for the former three ones
or $5\times 10^{-3}$ for the latter three ones leads to
the best TSVD solutions $x_n^{tsvd}=x_{naive}$. Therefore, there is no
semi-convergence, and regularization plays no role for them. In this case,
the mentioned Krylov iterative methods
solve \eqref{eq1} in their standard manners as if they solved
an ordinary other than ill-posed problem. On the other hand, if $\varepsilon$ is
relatively bigger, say, $0.05$, then the semi-convergence of the TSVD method
and LSQR may occur.

It has long been well known
\cite{hanke93,hansen98,hansen07,hansen10} and further
addressed in \cite{jia18a}
that provided that the Ritz values involved in LSQR approximate
the large singular values of $A$ in natural order until the occurrence
of semi-convergence, the best regularized solution
obtained by LSQR is as accurate as $x_{k_0}^{tsvd}$.
Unfortunately, as stressed by Hanke and Hansen~\cite{hanke93},
Hansen~\cite{hansen07} and many others, e.g.,
Gazzola and Novati \cite{gazzola15},
a strict proof of the regularizing properties of conjugate gradients is
extremely difficult and
proving if the Ritz values converge in this order is a difficult task.
In fact, up to now there has been no either general or specific rigorous
result on the approximation behavior of the Ritz values.

As matter of fact, as we have observed in \cite{jia18a},
the Ritz values converge to
the large singular values of $A$ in natural order for severely
ill-posed problems but they may
fail to do so at some iterations $k\leq k^*$ for some moderately
and mildly ill-posed problems, where $k^*$ is
the iteration at which the semi-convergence of LSQR occurs.
In the latter case, the regularization of LSQR is much more involved,
and hitherto nothing has been theoretically known on the
accuracy of the best regularized solution by LSQR at semi-convergence,
and a common belief seems that LSQR has the partial regularization.
However, the numerical experiments in \cite{jia18a} have indicated that
the best regularized solutions by LSQR are as accurate as $x_{k_0}^{tsvd}$
even if the Ritz values fail to approximate the large singular values of $A$
in natural order for some iterations $k\leq k^*$. For
the definition of severely, moderately and mildly ill-posed problems,
we refer to \cite{hofmann86}; also see \cite{jia18a} for a supplement.

In order to assess the regularization ability of a regularization method,
the definition of {\em full} and {\em partial} regularization is
introduced in \cite{huangjia,jia18a}. For the 2-norm filtering
regularization problem \eqref{posed}, if a regularized solution
is as accurate as $x_{k_0}^{tsvd}$, then it is called a 2-norm
filtering best possible regularized solution.
If a regularization method can compute such a best possible one,
then it is said to have the {\em full} regularization. Otherwise, it
is said to have
only the {\em partial} regularization. By such a definition, a natural
and fundamental question is:
{\em  Does LSQR have the full or partial regularization for
severely, moderately and mildly ill-posed problems?} This question
was implicit for CGLS in the survey paper \cite{bjorck79}
of Bj\"{o}rck and Eld{\'{e}}n and explicitly posed
in \cite{huangjia,jia18a}.

In \cite{jia18a}, the author has established a general $\sin\Theta$ theorem for
the 2-norm distances between the underlying
$k$ dimensional Krylov subspace and the $k$ dimensional dominant right
singular subspace of $A$, and derived accurate estimates for
these distances for severely, moderately
and mildly ill-posed problems, respectively. As has been addressed
in \cite{jia18a}, these results are the first key and fundamental step towards
to answering the question of LSQR having the full or partial
regularization. This paper is a continuation of
\cite{jia18a}. On the basis of \cite{jia18a},
for these three kinds of ill-posed problems, we will give accurate solutions of
the two closely related fundamental and highly
challenging problems on the regularization of LSQR:
(i) How accurate are the low rank approximations
generated by Lanczos bidiagonalization?
(ii) Whether or not the Ritz values involved in LSQR approximate
the large singular values of $A$ in natural order?
We establish accurate estimates for the accuracy of the low rank
approximations and give definitive results on how the Ritz values converge.
In addition, notice that the accuracy of low rank approximations is computationally
infeasible for $A$ large. We show how to judge it reliably during
the Lanczos bidiagonalization process without extra cost.

The paper is organized as follows. In Section \ref{lsqr}, we describe the
Lanczos bidiagonalization process and LSQR, and state some of the results
in \cite{jia18a} that are used to analyze the results in this paper.
In Section \ref{rankapp}, for severely and
moderately problems with suitable $\rho>1$ and $\alpha>1$
we prove that the $k$-step Lanczos bidiagonalization always
generates a near best rank $k$ approximation to $A$ and the $k$ Ritz values
always approximate the large singular values of $A$ in natural order
until the occurrence of semi-convergence. For mildly ill-posed problems,
we prove that the above results generally fail to
hold. In Section \ref{alphabeta}, we establish a monotonic
property of the accuracy of rank $k$ approximations generated by
Lanczos bidiagonalization, derive bounds for the decay rates of
entries of the bidiagonal
matrices generated by Lanczos bidiagonalization, and show that they
can be used to reliably judge the
accuracy of the low rank approximations.
In Section \ref{numer}, we report numerical experiments to confirm our results
and make some observations on the regularization of LSQR.
Finally, we conclude the paper and come to the
conjecture that LSQR has the full
regularization in Section \ref{concl}.

In the paper, we use
$\mathcal{K}_{k}(C, w)= span\{w,Cw,\ldots,C^{k-1}w\}$
to denote the $k$ dimensional Krylov subspace generated
by the matrix $\mathit{C}$ and the vector $\mathit{w}$, and by $I$ and the
bold letter $\mathbf{0}$ the identity matrix
and the zero matrix with orders clear from the context, respectively.
For the matrix $B=(b_{ij})$, define $|B|=(|b_{ij}|)$,
and for $|C|=(|c_{ij}|)$, $|B|\leq |C|$ means
$|b_{ij}|\leq |c_{ij}|$ componentwise.

\section{The LSQR algorithm and the estimates
for the distances between $\mathcal{V}_k^R$ and
$span\{V_k\}$}\label{lsqr}

The LSQR algorithm is based on Lanczos bidiagonalization that
computes two orthonormal bases $\{q_i\}_{i=1}^k$ and
$\{p_i\}_{i=1}^{k+1}$  of $\mathcal{K}_{k}(A^{T}A,A^{T}b)$ and
$\mathcal{K}_{k+1}(A A^{T},b)$  for $k=1,2,\ldots,n$,
respectively. For $k=1,2,..,n$,
the $k$-step Lanczos bidiagonalization process
can be written in the matrix form
\begin{align}
  AQ_k&=P_{k+1}B_k,\label{eqmform1}\\
  A^{T}P_{k+1}&=Q_{k}B_k^T+\alpha_{k+1}q_{k+1}(e_{k+1}^{(k+1)})^{T},\label{eqmform2}
\end{align}
where $e_{k+1}^{(k+1)}$ is the $(k+1)$-th canonical basis vector of
$\mathbb{R}^{k+1}$, $P_{k+1}=(p_1,p_2,\ldots,p_{k+1})$,
$Q_k=(q_1,q_2,\ldots,q_k)$, and
\begin{equation}\label{bk}
  B_k = \left(\begin{array}{cccc} \alpha_1 & & &\\ \beta_2 & \alpha_2 & &\\ &
  \beta_3 &\ddots & \\& & \ddots & \alpha_{k} \\ & & & \beta_{k+1}
  \end{array}\right)\in \mathbb{R}^{(k+1)\times k}.
\end{equation}
It is known from \eqref{eqmform1} that
\begin{equation}\label{Bk}
B_k=P_{k+1}^TAQ_k.
\end{equation}
We remark that the singular values of $B_k$, called the Ritz values of $A$ with
respect to the left and right subspaces $span\{P_{k+1}\}$ and $span\{Q_k\}$,
are all simple. It is easily justified that
Lanczos bidiagonalization can be run to
completion without breakdown since the starting vector $b$
has nonzero components in
all the $u_i$ and the singular values of $A$ are assumed to be simple.

Write $\mathcal{V}_k^R= \mathcal{K}_k(A^TA,A^Tb)$.
At iteration $k$, LSQR solves the problem
$\|Ax_k^{lsqr}-b\|=\min_{x\in \mathcal{V}_k^R}
\|Ax-b\|$ and computes the iterate $x_k^{lsqr}=Q_ky_k^{lsqr}$ with
\begin{equation}\label{yk}
  y_k^{lsqr}=\arg\min\limits_{y\in \mathbb{R}^{k}}\|B_ky-\beta_1 e_1^{(k+1)}\|
  =\beta_1  B_k^{\dagger} e_1^{(k+1)},
\end{equation}
where $e_1^{(k+1)}$ is the first canonical basis vector of $\mathbb{R}^{k+1}$.

Note that $\beta_1 e_1^{(k+1)}=P_{k+1}^T b$. From \eqref{yk} we have
\begin{equation}\label{xk}
x_k^{lsqr}=Q_k B_k^{\dagger} P_{k+1}^Tb,
\end{equation}
which is the minimum 2-norm solution to the perturbed
problem that replaces $A$ in \eqref{eq1} by its rank $k$ approximation
$P_{k+1}B_k Q_k^T$. In \cite{jia18a}, a key point is that
LSQR has been interpreted as solving
\begin{equation}\label{lsqrreg}
\min\|x\| \ \ \mbox{ subject to }\ \ \|P_{k+1}B_kQ_k^Tx-b\|=\min
\end{equation}
for the regularized solutions $x_k^{lsqr}$ of \eqref{eq1} starting with $k=1$
onwards. Therefore, LSQR is similar to the TSVD method and replaces the best
rank $k$ approximation $A_k$ to $A$ by a rank $k$ approximation
$P_{k+1}B_kQ_k^T$. Based on this connection,
the author \cite{jia18a} has
shown that the near best rank $k$ approximation of $P_{k+1}B_k Q_k^T$
to $A$ and the approximations of the $k$ singular values of $B_k$ to the
large ones of $A$ in natural order for $k=1,2,\ldots,k_0$ are
{\em sufficient} conditions for LSQR to have the full regularization.

Regarding the semi-convergence points $k^*$ and $k_0$
of LSQR and the TSVD method,
the author \cite{jia18a} has proved the following basic
property, which is useful to analyze some results
and numerical experiments in this paper.

\begin{theorem}\label{semicon}
The semi-convergence of LSQR must occur at some
iteration
$$
k^*\leq k_0.
$$
If the Ritz values $\theta_i^{(k)}$
do not converge to the large singular values of
$A$ in natural order for some $k\leq k^*$, then
$k^*<k_0$, and vice versa.
\end{theorem}

In terms of the canonical angles $\Theta(\mathcal{X},\mathcal{Y})$ between
two subspaces $\mathcal{X}$ and $\mathcal{Y}$ of equal
dimension (cf. \cite[p.74-5]{stewart01} and \cite[p.43]{stewartsun}),
the author \cite{jia18a} has established
accurate estimates on the accuracy of $\mathcal{V}_k^R$
approximating the $k$ dimensional dominant right singular
subspace $\mathcal{V}_k=span\{V_k\}$
for severely, moderately and mildly ill-posed problems. Since these estimates
play a central role in {\em analyzing}
the results in the next three sections, we state them
as Lemma~\ref{lemma1} and Theorems~\ref{thm2}--\ref{moderate}. To this end,
we introduce some notation that appeared in \cite{jia18a}
and will be used in this paper. Define
\begin{equation}\label{defdelta}
\Delta_k=D_2T_{k2}T_{k1}^{-1}D_1^{-1}\in \mathbb{R}^{(n-k)\times k},
\end{equation}
in which the matrices involved are
\begin{equation*}\label{}
  D={\rm diag}(\sigma_j u_j^Tb)\in\mathbb{R}^{n\times n},\ \
  T_k=\left(\begin{array}{cccc} 1 &
  \sigma_1^2&\ldots & \sigma_1^{2k-2}\\
1 &\sigma_2^2 &\ldots &\sigma_2^{2k-2} \\
\vdots & \vdots&&\vdots\\
1 &\sigma_n^2 &\ldots &\sigma_n^{2k-2}
\end{array}\right)
\end{equation*}
with the partitions
\begin{equation*}\label{}
  D=\left(\begin{array}{cc} D_1 & 0 \\ 0 & D_2 \end{array}\right),\ \ \
  T_k=\left(\begin{array}{c} T_{k1} \\ T_{k2} \end{array}\right)
\end{equation*}
and $D_1, T_{k1}\in\mathbb{R}^{k\times k}$. Then we have
the following precise $\sin\Theta$ theorem on the 2-norm distance between
$\mathcal{V}_k^R$ and $\mathcal{V}_k$.

\begin{lemma}[\cite{jia18a}]\label{lemma1}
For $k=1,2,\ldots,n-1$ we have
\begin{align}
\|\sin\Theta(\mathcal{V}_k,\mathcal{V}_k^R)\|&=
\frac{\|\Delta_k\|}{\sqrt{1+\|\Delta_k\|^2}}
\label{deltabound}
\end{align}
with $\Delta_k \in \mathbb{R}^{(n-k)\times k}$ defined by \eqref{defdelta}.
\end{lemma}

From the lemma it is direct to get
\begin{equation}\label{tangent}
\|\tan\Theta(\mathcal{V}_k,\mathcal{V}_k^R)\|=\|\Delta_k\|.
\end{equation}
The following two theorems give estimates for $\|\Delta_k\|$ for
the three kinds of ill-posed problems, which have been shown and numerically
confirmed to be accurate in \cite{jia18a}.

\begin{theorem}[\cite{jia18a}]\label{thm2}
Let the SVD of $A$ be as \eqref{eqsvd}. Assume that \eqref{eq1} is severely
ill-posed with $\sigma_j=\mathcal{O}(\rho^{-j})$ and $\rho>1$, $j=1,2,\ldots,n$.
Then
\begin{align}
\|\Delta_k\|&\leq
\frac{\sigma_{k+1}}{\sigma_k}\frac{\max_{k+1\leq i\leq n}|u_i^Tb|}{\min_{1\leq i\leq k}
  |u_i^Tb|}
  \left(1+\mathcal{O}(\rho^{-2})\right).\label{eqres1}
\end{align}
\end{theorem}

\begin{theorem}[\cite{jia18a}]\label{moderate}
For a moderately or mildly ill-posed  \eqref{eq1} with $\sigma_j=
\zeta j^{-\alpha},\ j=1,2,\ldots,n$,
where $\alpha>\frac{1}{2}$ and $\zeta>0$ is some constant, we have
\begin{align}
\|\Delta_1\|&\leq \frac{\max_{2\leq i\leq n}|u_i^Tb|}{| u_1^Tb|}
\sqrt{\frac{1}{2\alpha-1}},\label{mod1}\\
\|\Delta_k\|&\leq \frac{\max_{k+1\leq i\leq n}|u_i^T b|}{
\min_{1\leq i\leq k}| u_i^T b|}
\sqrt{\frac{k^2}{4\alpha^2-1}+\frac{k}{2\alpha-1}}|L_{k_1}^{(k)}(0)|,
\ k=2,3,\ldots,n-1. \label{modera2}
\end{align}
where
\begin{equation}\label{lk}
|L_{k_1}^{(k)}(0)|=\max_{j=1,2,\ldots,k}|L_j^{(k)}(0)|,
\ |L_j^{(k)}(0)|=\prod\limits_{i=1,i\ne j}^k\frac{\sigma_i^2}{|\sigma_j^2-
\sigma_i^2|},\,j=1,2,\ldots,k.
\end{equation}
\end{theorem}

The author in \cite{jia18a} has derived estimates for $|L_j^{(k)}(0)|,\
j=1,2,\ldots,k$ and $|L_{k_1}^{(k)}(0)|$ for $k=2,3,\ldots,n-1$
for the moderately and mildly ill-posed problems.
It is shown that $|L_{k_1}^{(k)}(0)|>1$ becomes large soon for
mildly ill-posed problems as $k$ increases and it is $1+\mathcal{O}(k)$
for moderately ill-posed problems with suitable $\alpha>1$.

\section{The rank $k$ approximation $P_{k+1}B_kQ_k^T$ to $A$}\label{rankapp}

This section is devoted to the study of the quality of rank $k$
approximations generated by Lanczos bidiagonalization, which
replaces the highly ill-conditioned $A$ in \eqref{eq1} and
has fundamental effects on the regularization of LSQR.
We first establish the following
intermediate results, which play a key role in deriving
our main results.

\begin{theorem}\label{thm3}
Let $\Delta_k\in \mathbb{R}^{(n-k)\times k}$ and $L_{k_1}^{(k)}(0)$
be defined as \eqref{defdelta} and \eqref{lk}, and
$\Sigma_k={\rm diag}(\sigma_1,\sigma_2,\ldots,\sigma_k)$. Then for severely
ill-posed problems we have
\begin{equation} \label{prodnorm}
\|\Sigma_k\Delta_k^T\|\leq \left\{\begin{array}{ll}
\sigma_{k+1}\frac{\max_{k+1\leq i\leq n}|u_i^T b|}{
\min_{1\leq i\leq k}| u_i^T b|}
\left(1+\mathcal{O}(\rho^{-2})\right)
& \mbox{ for } 1\leq k\leq k_0,\\
\sigma_{k+1}\frac{\max_{k+1\leq i\leq n}|u_i^T b|}{
\min_{1\leq i\leq k}| u_i^T b|}\sqrt{k-k_0+1}\left(1+\mathcal{O}(\rho^{-2})\right)
& \mbox{ for } k_0<k\leq n-1,
\end{array}
\right.
\end{equation}
and for moderately or mild ill-posed problems with the singular values
$\sigma_j=\zeta j^{-\alpha}$, $\alpha>\frac{1}{2}$
and $\zeta$ a positive constant we have
\begin{equation}\label{prodnorm2}
\|\Sigma_k\Delta_k^T\|\leq \left\{\begin{array}{ll}
\sigma_1\frac{\max_{2\leq i\leq n}| u_i^Tb|}{| u_1^T b|}
\sqrt{\frac{1}{2\alpha-1}} &
 \mbox{ for } k=1,\\
\sigma_k\frac{\max_{k+1\leq i\leq n}|u_i^T b|}{
\min_{1\leq i\leq k}| u_i^T b|}\sqrt{\frac{k^2}{ (4\alpha^2-1)}+
\frac{k}{2\alpha-1}}
|L_{k_1}^{(k)}(0)|& \mbox{ for } 1<k\leq k_0, \\
\sigma_k\frac{\max_{k+1\leq i\leq n}|u_i^T b|}{
\min_{1\leq i\leq k}| u_i^T b|}\sqrt{\frac{k k_0}{(4\alpha^2-1)}+
\frac{k(k-k_0+1)}{2\alpha-1}}|L_{k_1}^{(k)}(0)|
& \mbox{ for } k_0<k\leq n-1.
\end{array}
\right.
\end{equation}
\end{theorem}

{\em Proof}.
It has been proved in \cite{jia18a} (cf. the proof of Theorem 4.2 there) that
\begin{align}
|\Delta_k|&=|D_2T_{k2}T_{k1}^{-1}D_1^{-1}|
\leq |L_{k_1}^{(k)}(0)||\tilde\Delta_k| \label{amplify}
\end{align}
with the definition $L_{1}^{(1)}(0)=1$, where
\begin{equation}
|\tilde\Delta_k|=\left|(\sigma_{k+1} u_{k+1}^T b,\sigma_{k+2}u_{k+2}^Tb,
\ldots,\sigma_n u_n^T b)^T
\left(\frac{1}{\sigma_1 u_1^Tb},\frac{1}{\sigma_2 u_2^Tb},\ldots,
\frac{1}{\sigma_k u_k^Tb}\right)\right|  \label{delta1}
\end{equation}
is a rank one matrix. Therefore, we have
$$
|\Delta_k\Sigma_k|\leq |L_{k_1}^{(k)}(0)|\left|(\sigma_{k+1} u_{k+1}^T b,
\sigma_{k+2}u_{k+2}^Tb,\ldots,\sigma_n u_n^T b)^T
\left(\frac{1}{u_1^Tb},\frac{1}{u_2^Tb},\ldots,
\frac{1}{u_k^Tb}\right)\right|,
$$
from which we obtain
\begin{align}
\|\Sigma_k\Delta_k^T\|&=\|\Delta_k\Sigma_k\|\leq \left\||\Delta_k\Sigma_k|\right\|
\notag\\
&\leq |L_{k_1}^{(k)}(0)|\left(\sum_{j=k+1}^n\sigma_j^2| u_j^Tb|^2\right)^{1/2}
\left(\sum_{j=1}^k \frac{1}{| u_j^Tb|^2}\right)^{1/2}. \label{sigdel}
\end{align}
For the severely ill-posed problems and moderately or mildly ill-posed problems,
it has been proved in \cite{jia18a} (cf. the proofs of Theorems 4.2
and 4.4, respectively) that
\begin{align}
\left(\sum_{j=k+1}^n\sigma_j^2| u_j^Tb|^2\right)^{1/2}
&=\sigma_{k+1}\max_{k+1\leq i\leq n}|u_i^Tb| \left(1+\mathcal{O}(\rho^{-2})\right)
\label{severe1}
\end{align}
with $1+\mathcal{O}(\rho^{-2})=1$ for $k=n-1$
and
\begin{align}
\left(\sum_{j=k+1}^n\sigma_j^2| u_j^Tb|^2\right)^{1/2}
&\leq \sigma_k \max_{k+1\leq i\leq n}|u_i^Tb|\sqrt{\frac{k}
{2\alpha-1}},
\label{modeest}
\end{align}
respectively, from
which and \eqref{sigdel} we obtain \eqref{prodnorm} and \eqref{prodnorm2}
for $k=1$.

In order to bound $\|\Sigma_k\Delta_k^T\|$ for $k>1$, we need to estimate
$\left(\sum_{j=1}^k\frac{1}{| u_j^Tb|^2}\right)^{1/2}$.
Next we do this for severely and moderately or mildly ill-posed
problems, respectively, for each kind of which we need to consider the cases
$k\leq k_0$ and $k>k_0$ separately.

By the discrete Picard condition \eqref{picard}, \eqref{picard1} and the
properties on $e$, it is known from \cite[p.70-1]{hansen98}
and \cite[p.41-2]{hansen10} that
$
| u_i^T b|\approx| u_i^T b_{true}|=\sigma_i^{1+\beta}>\eta
$
monotonically decreases with $i=1,2,\ldots,k_0$,
and becomes stabilized for $i>k_0$ since
$| u_i^T b|\approx | u_i^T e |$ with the expected values
$
\mathcal{E}(| u_i^T e |)=\eta.
$
Therefore,
to present our derivation and results compactly and elegantly,
in the later proof we will use the {\em ideal} equality
\begin{equation}\label{idealeq}
|u_i^Tb|= |u_i^Tb_{true}|=\sigma_i^{1+\beta},\ i=1,2,\ldots,k_0
\end{equation}
by dropping the negligible $|u_i^Te|$.

The case $k\leq k_0$ for severely ill-posed problems:
From \eqref{idealeq}, we have
$
\min_{1\leq i\leq k}| u_i^Tb|=| u_k^Tb|\leq |u_{j+1}^T b|
$
for $j=1,2,\ldots,k-1$. Making use of \eqref{picard} and
\eqref{idealeq}, we obtain
\begin{align*}
\sum_{j=1}^k
\frac{1}{| u_j^Tb|^2}&= \frac{1}{\min_{1\leq i\leq k} |u_i^Tb|^2} \sum_{j=1}^k
\frac{\min_{1\leq i\leq k} |u_i^Tb|^2}{| u_j^Tb|^2}\\
&\leq \frac{1}{\min_{1\leq i\leq k} |u_i^Tb|^2}\left(
\frac{| u_k^Tb|^2}{| u_k^Tb|^2}+ \sum_{j=1}^{k-1}
\frac{| u_{j+1}^Tb|^2}{| u_j^Tb|^2}\right)\\
&\leq \frac{1}{\min_{1\leq i\leq k} |u_i^Tb|^2}\left(1+\sum_{j=1}^{k-1}
\frac{\sigma_{j+1}^2}{\sigma_j^2}\right)\\
&\leq \frac{1}{\min_{1\leq i\leq k} |u_i^Tb|^2}
\left(1+\mathcal{O}\left(\sum_{j=1}^{k-1}\rho^{2(j-k)}\right)\right)\\
&=\frac{1}{\min_{1\leq i\leq k} |u_i^Tb|^2}
\left(1+\mathcal{O}(\rho^{-2}) \right).
\end{align*}

The case  $k> k_0$ for severely ill-posed problems: Exploiting the above result
for $k\leq k_0$ and
$\min_{1\leq i\leq k}| u_i^Tb| \leq | u_j^Tb|$ for $j=k_0+1,\ldots,k$ for $k>k_0$,
we obtain
\begin{align*}
\sum_{j=1}^k
\frac{1}{| u_j^Tb|^2}&= \frac{1}{\min_{1\leq i\leq k}| u_i^Tb|^2}
\left( \sum_{j=1}^{k_0}
\frac{\min_{1\leq i\leq k}| u_i^Tb|^2}{| u_j^Tb|^2}+\sum_{j=k_0+1}^{k}
\frac{\min_{1\leq i\leq k}| u_i^Tb|^2}{| u_j^Tb|^2}\right)\\
&\leq \frac{1}{\min_{1\leq i\leq k}| u_i^Tb|^2}
\left( \sum_{j=1}^{k_0}
\frac{\min_{1\leq i\leq k}| u_i^Tb|^2}{| u_j^Tb|^2}+\sum_{j=k_0+1}^{k}
\frac{\min_{1\leq i\leq k}| u_i^Tb|^2}{| u_j^Tb|^2}\right)\\
&\leq \frac{1}{\min_{1\leq i\leq k}| u_i^Tb|^2}
\left(1+\mathcal{O}\left(\sum_{j=1}^{k_0-1}
\rho^{2(j-k_0)}\right)+k-k_0\right)\\
&=\frac{1}{\min_{1\leq i\leq k} |u_i^Tb|^2}
\left(1+\mathcal{O}(\rho^{-2})+k-k_0\right).
\end{align*}
Substitute the above two relations for the two cases into \eqref{sigdel} and
combine them with \eqref{severe1} and $|L_{k_1}^{(k)}(0)|=1+
\mathcal{O}(\rho^{-2})$ proved in \cite{jia18a}. We then obtain \eqref{prodnorm}.

The case $k\leq k_0$ for moderately or mildly ill-posed problems:
In a similar way to the above, we have
\begin{align*}
\sum_{j=1}^k
\frac{1}{| u_j^Tb|^2}&= \frac{1}{\min_{1\leq i\leq k} |u_i^Tb|^2} \sum_{j=1}^k
\frac{\min_{1\leq i\leq k} |u_i^Tb|^2}{| u_j^Tb|^2}\\
&\leq  \frac{1}{\min_{1\leq i\leq k} |u_i^Tb|^2} \left(1+\sum_{j=1}^{k-1}
\frac{| u_{j+1}^Tb|^2}{| u_j^Tb|^2}\right)\\
&\leq \frac{1}{\min_{1\leq i\leq k} |u_i^Tb|^2}\left(1+\sum_{j=1}^{k-1}
\frac{\sigma_{j+1}^2}{\sigma_j^2}\right)\\
&\leq \frac{1}{\min_{1\leq i\leq k} |u_i^Tb|^2}\left(1+\sum_{j=1}^{k-1}
\left(\frac{j}{k}\right)^{2\alpha}\right)\\
&=\frac{1}{\min_{1\leq i\leq k}| u_i^T b|^2} k\left(1+\sum_{j=1}^{k-1}
\frac{1}{k}\left(\frac{j}{k}
\right)^{2\alpha}\right) \\
&\leq \frac{1}{\min{1\leq i\leq k}| u_i^T b|^2} \left(1+k \int_0^1
x^{2\alpha }dx \right) \\
&=\frac{1}{\min_{1\leq i\leq k}| u_i^Tb|^2}
\left(1+\frac{k}{2\alpha+1}\right).
\end{align*}

The case $k> k_0$ for moderately or mildly ill-posed problems:
Exploiting the above manner, we have
\begin{align*}
\sum_{j=1}^k\frac{1}{| u_j^Tb|^2}
&= \frac{1}{\min_{1\leq i\leq k} |u_i^Tb|^2} \left(\sum_{j=1}^{k_0}
\frac{\min_{1\leq i\leq k} |u_i^Tb|^2}{| u_j^Tb|^2}
+\sum_{j=k_0+1}^{k}
\frac{\min_{1\leq i\leq k} |u_i^Tb|^2}{| u_j^Tb|^2}\right)\\
&\leq \frac{1}{\min_{1\leq i\leq k} |u_i^Tb|^2} \left(\sum_{j=1}^{k_0}
\left(\frac{j}{k_0}\right)^{2\alpha }+k-k_0\right)\\
&\leq \frac{1}{\min_{1\leq i\leq k} |u_i^Tb|^2}
\left(1+\frac{k_0}{2\alpha +1}+k-k_0\right).
\end{align*}
Substitute the above two bounds for the two cases into \eqref{sigdel} and
combine them with \eqref{modeest}. We then obtain \eqref{prodnorm2}.
\qquad\endproof

Regarding the factor
$\frac{\max_{k+1\leq i\leq n}|u_i^Tb|}{\min_{1\leq i\leq k} |u_i^Tb|}$,
based on \eqref{picard} and the properties of $e$ described
in the introduction,
we can easily justify (cf. \cite{jia18a}) that
\begin{align}
\frac{\max_{k+1\leq i\leq n}|u_i^Tb|}{\min_{1\leq i\leq k} |u_i^Tb|}
&\approx \frac{|u_{k+1}^Tb|}{| u_k^Tb|}\approx \frac{\sigma_{k+1}^{1+\beta}}
{\sigma_k^{1+\beta}}<1, \ k=1,2,\ldots,k_0,\label{replace}\\
\frac{\max_{k+1\leq i\leq n}|u_i^Tb|}{\min_{1\leq i\leq k} |u_i^Tb|}
&\approx \frac{|u_{k+1}^Tb|}{| u_k^Tb|}\approx \frac{\eta}{\eta}=1, \
k=k_0+1,\ldots,n-1.\label{replace1}
\end{align}

\eqref{prodnorm} and \eqref{prodnorm2} indicate that $\|\Sigma_k\Delta_k^T\|$
decays swiftly as $k$ increases up to $k_0$ for severely and moderately
ill-posed problems. Trivially, we have
$$
\sigma_k\|\Delta_k\|\leq \|\Sigma_k\Delta_k^T\|\leq \sigma_1 \|\Delta_k\|.
$$
By carefully
comparing the accurate estimates for
$\|\Sigma_k\Delta_k^T\|$ with those for $\|\Delta_k\|$
in Theorems~\ref{thm2}--\ref{moderate}, for $k\leq k_0$, Theorem~\ref{thm3}
indicates that
$$
\sigma_k\|\Delta_k\|\leq \|\Sigma_k\Delta_k^T\|
\approx \sigma_k\|\Delta_k\|,
$$
that is, the true $\|\Sigma_k\Delta_k^T\|$ approximately attains its
sharpest lower bound, and it is impossible to improve the estimate and
get a smaller one.
In contrast, for $k=2,3,\ldots, k_0$,
if we had estimated it by its simple upper bound
\begin{equation}\label{rough}
\|\Sigma_k\Delta_k^T\|\leq \|\Sigma_k\|\|\Delta_k^T\|=\sigma_1\|\Delta_k\|,
\end{equation}
we would have obtained a bound which not only does not decay but also
increases for moderately and mildly ill-posed problems as $k$ increases.
Such bound is useless to derive those accurate bounds to be presented.
The subtlety to bound $\|\Sigma_k\Delta_k^T\|$ consists in deriving
\eqref{amplify} and \eqref{sigdel} and bounding
$\|\Sigma_k\Delta_k^T\|$ as a whole.

Making use of Theorems~\ref{thm2}--\ref{moderate} and Theorem~\ref{thm3},
in what follows we prove
that, at iteration $k$, Lanczos bidiagonalization generates a near best
rank $k$ approximation to $A$ and the $k$ Ritz values $\theta_i^{(k)}$
approximate the large singular values $\sigma_i$ of $A$ in natural order
for severely or moderately ill-posed problems with suitable $\rho>1$ or
$\alpha>1$ for $k\leq k_0$, but these two results
fail to hold for mildly ill-posed problems for some $k\leq k^*$.
By Theorem~\ref{semicon},
this means that $k^*<k_0$ for mildly ill-posed problems and for severely
or moderately ill-problems with $\rho>1$ or $\alpha>1$ not enough.

Define
\begin{equation}\label{gammak}
\gamma_k = \|A-P_{k+1}B_kQ_k^T\|,
\end{equation}
which measures the accuracy of the rank $k$ approximation $P_{k+1}B_kQ_k^T$ to
$A$. We will introduce the precise definition of near best
rank $k$ approximation to $A$ later on.

\begin{theorem}\label{main1}
Let $|L_{k_1}^{(k)}(0)|$ be defined by \eqref{lk}.
Then
we have
\begin{equation}\label{final}
  \sigma_{k+1}\leq \gamma_k\leq \sqrt{1+\eta_k^2}\sigma_{k+1}
\end{equation}
with
\begin{equation} \label{const1}
\eta_k\leq \left\{\begin{array}{ll}
\xi_k\frac{\max_{k+1\leq i\leq n}|u_i^Tb|}{\min_{1\leq i\leq k}| u_i^T b|}
\left(1+\mathcal{O}(\rho^{-2})\right)
& \mbox{ for } 1\leq k\leq k_0,\\
\xi_k\frac{\max_{k+1\leq i\leq n}|u_i^Tb|}{\min_{1\leq i\leq k}| u_i^T b|}
\sqrt{k-k_0+1}\left(1+\mathcal{O}(\rho^{-2})\right)
& \mbox{ for } k_0<k \leq n-1
\end{array}
\right.
\end{equation}
for severely ill-posed problems with $\sigma_i=\mathcal{O}(\rho^{-i})$,
$i=1,2,\ldots,n$ and
\begin{equation}\label{const2}
\eta_k\leq \left\{\begin{array}{ll}
\xi_1\frac{\sigma_1}{\sigma_2}\frac{\max_{2\leq i\leq n}| u_i^Tb|}{| u_1^Tb|}
\sqrt{\frac{1}{2\alpha-1}} & \mbox{ for } k=1, \\
\xi_k\frac{\sigma_k}{\sigma_{k+1}}\frac{\max_{k+1\leq i\leq n}|u_i^Tb|}
{\min_{1\leq i\leq k}| u_i^T b|}
\sqrt{\frac{k^2}{4\alpha^2-1}+\frac{k}{2\alpha-1}}
|L_{k_1}^{(k)}(0)|& \mbox{ for } 1< k\leq k_0, \\
\xi_k\frac{\sigma_k}{\sigma_{k+1}}\frac{\max_{k+1\leq i\leq n}|u_i^Tb|}
{\min_{1\leq i\leq k}| u_i^T b|}\sqrt{\frac{k k_0}{4\alpha^2-1
}+\frac{k(k-k_0+1)}{2\alpha-1}}|L_{k_1}^{(k)}(0)|
& \mbox{ for } k_0<k\leq n-1
\end{array}
\right.
\end{equation}
for moderately or mildly ill-posed problems with $\sigma_i=\zeta i^{-\alpha},\
i=1,2,\ldots,n$, where
$\xi_k=\sqrt{\left(\frac{\|\Delta_k\|}{1+\|\Delta_k\|^2}\right)^2+1}$ for
$\|\Delta_k\|<1$ and $\xi_k\leq\frac{\sqrt{5}}{2}$ for $\|\Delta_k\|\geq 1$.
\end{theorem}

{\em Proof}.
Since $A_k$ is the best rank $k$ approximation to $A$ and $\|A-A_k\|=\sigma_{k+1}$,
the lower bound in \eqref{final} holds trivially. Next we prove the upper bound.

From \eqref{eqmform1}, we obtain
\begin{align}
\gamma_k
&= \|A-P_{k+1}B_kQ_k^T\|= \|A-AQ_kQ_k^T\|= \|A(I-Q_kQ_k^T)\|. \label{gamma2}
\end{align}
From the proof of Lemma~\ref{lemma1} (cf. Lemma 4.1 \cite{jia18a}), we obtain
$$
\mathcal{V}_k^R
=\mathcal{K}_{k}(A^{T}A,A^{T}b)=span\{Q_k\}=span\{\hat{Z}_k\},
$$
where the orthonormal $Q_k$ is generated by Lanczos bidiagonalization and
\begin{equation}\label{decomp}
\hat{Z}_k=(V_k+V_k^{\perp}\Delta_k)(I+\Delta_k^T\Delta_k)^{-\frac{1}{2}}.
\end{equation}
Therefore, the orthogonal projector onto $\mathcal{V}_k^R$ is
$Q_kQ_k^T=\hat{Z}_k\hat{Z}_k^T$.
Keep in mind that $A_k=U_k\Sigma_k V_k^T$. It is direct to justify that
$(U_k\Sigma_k V_k^T)^T(A-U_k\Sigma_k V_k^T)=\mathbf{0}$ for $k=1,2,\ldots,n-1$.
Exploiting this and noting that $\|I-\hat{Z}_k\hat{Z}_k^T\|=1$ and
$V_k^TV_k^{\perp}=\mathbf{0}$  for $k=1,2,\ldots,n-1$,
we get from \eqref{gamma2} and \eqref{decomp} that
\begin{align}
\gamma_k^2 &= \|(A-U_k\Sigma_kV_k^T+U_k\Sigma_kV_k^T)(I-\hat{Z}_k\hat{Z}_k^T)\|^2
  \notag\\
  &=\max_{\|y\|=1}\|(A-U_k\Sigma_kV_k^T+U_k\Sigma_kV_k^T)
  (I-\hat{Z}_k\hat{Z}_k^T)y\|^2 \notag\\
  &=\max_{\|y\|=1}\|(A-U_k\Sigma_kV_k^T)(I-\hat{Z}_k\hat{Z}_k^T)y+
 U_k\Sigma_kV_k^T(I-\hat{Z}_k\hat{Z}_k^T)y\|^2\notag\\
 &=\max_{\|y\|=1}\left(\|(A-U_k\Sigma_kV_k^T)(I-\hat{Z}_k\hat{Z}_k^T)y\|^2+
 \| U_k\Sigma_kV_k^T(I-\hat{Z}_k\hat{Z}_k^T)y\|^2\right)\notag\\
 &\leq \|(A-U_k\Sigma_kV_k^T)(I-\hat{Z}_k\hat{Z}_k^T)\|^2+
 \| U_k\Sigma_kV_k^T(I-\hat{Z}_k\hat{Z}_k^T)\|^2 \notag \\
 &\leq \sigma_{k+1}^2+\| \Sigma_kV_k^T(I-\hat{Z}_k\hat{Z}_k^T)\|^2 \notag\\
 &= \sigma_{k+1}^2+\|\Sigma_kV_k^T\left(I-(V_k+V_k^{\perp}\Delta_k)(I+
 \Delta_k^T\Delta_k)^{-1}(V_k+V_k^{\perp}\Delta_k)^T\right)\|^2\notag\\
 &= \sigma_{k+1}^2 + \left\|\Sigma_k\left(V_k^T-(I+
 \Delta_k^T\Delta_k)^{-1}(V_k+V_k^{\perp}\Delta_k)^T\right)\right\|^2 \notag\\
  &= \sigma_{k+1}^2 + \left\|\Sigma_k(I+
 \Delta_k^T\Delta_k)^{-1}\left((I+\Delta_k^T\Delta_k)V_k^T-
 \left(V_k+V_k^{\perp}\Delta_k\right)^T\right)\right\|^2 \notag\\
  &= \sigma_{k+1}^2+ \|\Sigma_k(I+
 \Delta_k^T\Delta_k)^{-1}\left(\Delta_k^T\Delta_kV_k^T-\Delta_k^T
 (V_k^{\perp})^T\right)\|^2\notag\\
  &= \sigma_{k+1}^2 + \|\Sigma_k(I+
 \Delta_k^T\Delta_k)^{-1}\Delta_k^T\Delta_kV_k^T-\Sigma_k(I+
 \Delta_k^T\Delta_k)^{-1}\Delta_k^T(V_k^{\perp})^T\|^2 \notag\\
    &\leq \sigma_{k+1}^2 +\|\Sigma_k(I+
 \Delta_k^T\Delta_k)^{-1}\Delta_k^T\Delta_k\|^2+\|\Sigma_k(I+
 \Delta_k^T\Delta_k)^{-1}\Delta_k^T\|^2\notag\\
 &=\sigma_{k+1}^2+\epsilon_k^2.
 \label{estimate1}
\end{align}

We estimate $\epsilon_k$ accurately below. To this end, we need to use two
key identities and some results related. By the SVD of $\Delta_k$, it is direct
to justify that
\begin{equation}\label{inden1}
(I+
 \Delta_k^T\Delta_k)^{-1}\Delta_k^T\Delta_k=\Delta_k^T\Delta_k(I+
 \Delta_k^T\Delta_k)^{-1}
\end{equation}
and
\begin{equation}\label{inden2}
(I+
 \Delta_k^T\Delta_k)^{-1}\Delta_k^T=\Delta_k^T(I+
 \Delta_k\Delta_k^T)^{-1}.
\end{equation}
Define the function $f(\lambda)=\frac{\lambda}{1+\lambda^2}$ with
$\lambda\in [0,\infty)$. Since the derivative
$f^{\prime}(\lambda)=\frac{1-\lambda^2}{(1+\lambda^2)^2}$,
$f(\lambda)$ is monotonically increasing for $\lambda\in [0,1]$
and decreasing for $\lambda\in [1,\infty)$, and the maximum of
$f(\lambda)$ over $\lambda\in [0,\infty)$ is $\frac{1}{2}$, which attains at
$\lambda=1$. Based on these properties and exploiting the SVD of $\Delta_k$,
we obtain
\begin{equation}\label{compact}
\|\Delta_k(I+\Delta_k^T\Delta_k)^{-1}\|=\frac{\|\Delta_k\|}{1+\|\Delta_k\|^2}
\end{equation}
for $\|\Delta_k\|<1$ and
\begin{equation}\label{noncomp}
\|\Delta_k(I+\Delta_k^T\Delta_k)^{-1}\|\leq\frac{1}{2}
\end{equation}
for $\|\Delta_k\|\geq 1$ (Note: in this case, since $\Delta_k$ may have at least
one singular value smaller than one, we do not
have the expression \eqref{compact}). It then follows
from \eqref{estimate1}, \eqref{compact}, \eqref{noncomp}
and $\|(1+\Delta_k\Delta_k^T)^{-1}\|\leq 1$ that
\begin{align}
\epsilon_k^2&=\|\Sigma_k \Delta_k^T\Delta_k(I+
 \Delta_k^T\Delta_k)^{-1}\|^2+\|\Sigma_k \Delta_k^T (I+
 \Delta_k\Delta_k^T)^{-1}\|^2 \label{separa}\\
 &\leq \|\Sigma_k\Delta_k^T\|^2\|\Delta_k(I+\Delta_k^T\Delta_k)^{-1}\|^2+
 \|\Sigma_k\Delta_k^T\|^2 \|(1+\Delta_k\Delta_k^T)^{-1}\|^2 \notag\\
 &\leq \|\Sigma_k\Delta_k^T\|^2\left(\|\Delta_k
 (I+\Delta_k^T\Delta_k)^{-1}\|^2+1\right)\notag\\
  &=\|\Sigma_k\Delta_k^T\|^2\left(\left(\frac{\|\Delta_k\|}
  {1+\|\Delta_k\|^2}\right)^2+1\right)
  =\xi_k^2\|\Sigma_k\Delta_k^T\|^2 \notag
\end{align}
for $\|\Delta_k\|<1$ and
$$
\epsilon_k\leq \|\Sigma_k\Delta_k^T\|\sqrt{\|\Delta_k
(I+\Delta_k^T\Delta_k)^{-1}\|^2+1}=\xi_k\|\Sigma_k\Delta_k^T\|
\leq \frac{\sqrt{5}}{2}\|\Sigma_k\Delta_k^T\|
$$
for $\|\Delta_k\|\geq 1$. Replace $\|\Sigma_k\Delta_k^T\|$ by
its bounds \eqref{prodnorm} and \eqref{prodnorm2} in the
above, insert the resulting bounds for $\epsilon_k$
into \eqref{estimate1}, and let $\epsilon_k=\eta_k\sigma_{k+1}$.
Then we obtain the upper bound in \eqref{final} with $\eta_k$
satisfying \eqref{const1} and \eqref{const2} for severely and moderately
or mildly ill-posed problems, respectively.
\qquad\endproof


\begin{remark}\label{decayrate}
For severely ill-posed problems, from Theorem~\ref{thm2} and
\eqref{replace}--\eqref{replace1} we approximately have
\begin{align}
\|\Delta_k\|&\leq\frac{\sigma_{k+1}^{2+\beta}}{\sigma_k^{2+\beta}}
\left(1+\mathcal{O}(\rho^{-2})\right)\sim
\rho^{-2-\beta},\ k=1,2,\ldots, k_0,
\label{case3}\\
\|\Delta_k\|&\leq\frac{\sigma_{k+1}}{\sigma_k}\left(1+\mathcal{O}(\rho^{-2})\right)
\sim \rho^{-1},
\ k=k_0+1,\ldots,n-1, \label{case4}
\end{align}
from which and the definition of $\xi_k$ it follows that
$$
\xi_k(1+\mathcal{O}(\rho^{-2}))
=1+\mathcal{O}(\rho^{-2})
$$
for both $k\leq k_0$ and $k>k_0$. Therefore, from
\eqref{const1}, for $k\leq k_0$ we obtain
\begin{equation}\label{etak0}
\eta_k\leq \xi_k
\frac{|u_{k+1}^Tb |}{|u_k^T b |}\left(1+\mathcal{O}(\rho^{-2})\right)
=\frac{|u_{k+1}^Tb |}{|u_k^T b |}\left(1+\mathcal{O}(\rho^{-2})\right)
=\rho^{-1-\beta}<1
\end{equation}
by dropping the smaller $\mathcal{O}(\rho^{-3-\beta})$. From \eqref{replace1},
for $k>k_0$ we obtain
\begin{equation}\label{incres}
\eta_k\leq \xi_k\frac{\max_{k+1\leq i\leq n}|u_i^Tb|}
{\min_{1\leq i\leq k}| u_i^T b|}
\sqrt{k-k_0+1}\left(1+\mathcal{O}(\rho^{-2})\right)
\sim \sqrt{k-k_0+1}.
\end{equation}
\end{remark}

\begin{remark}\label{decayrate2}
From \eqref{final}, \eqref{const1} and \eqref{etak0},
for severely ill-posed problems and $k\leq k_0$ we have
$$
1<\sqrt{1+\eta_k^2}<1+\frac{1}{2}{\eta_k^2}\leq
1+\frac{1}{2}\frac{\sigma_{k+1}^{2(1+\beta)}}{\sigma_k^{2(1+\beta)}}
\sim 1+\frac{1}{2}\rho^{-2(1+\beta)},
$$
which and \eqref{final} indicate that $\gamma_k$ is an accurate
approximation to $\sigma_{k+1}$.
Thus, the rank $k$ approximation $P_{k+1}B_kQ_k^T$ is as accurate as
the best rank $k$ approximation $A_k$ within the
factor $\sqrt{1+\eta_k^2}\approx 1$ for suitable $\rho>1$ and $k\leq k_0$.
In contrast, \eqref{final} and \eqref{incres} shows that $\gamma_k$
is a marginally less accurate approximation to $\sigma_{k+1}$ for $k>k_0$.
\end{remark}

\begin{remark}
For the moderately or mildly ill-posed problems with
$\sigma_i=\zeta i^{-\alpha}$,
from the estimate \eqref{const2} for $\eta_k$,
for $k\leq k_0$ we approximately have
\begin{equation}\label{etadelta2}
\frac{\sigma_k}{\sigma_{k+1}}\|\Delta_k\|\leq
\eta_k\leq \frac{\sqrt{5}}{2}\frac{\sigma_k}{\sigma_{k+1}}\|\Delta_k\|.
\end{equation}
Therefore, based on Theorems~\ref{moderate}-\ref{moderate},
it is known that $P_{k+1}B_kQ_k^T$ is almost as accurate as
$A_k$ for suitable $\alpha>1$ and $k\leq k^*$ with $k^*$ not large.
\end{remark}

\begin{remark}
For both severely and moderately ill-posed problems, we notice that the
situation is not so satisfying for increasing $k>k_0$. But at that time,
a possibly big $\eta_k$ does not do harm to our regularization purpose
since Theorem~\ref{semicon} shows that the semi-convergence of LSQR must
occur at some iteration $k^*\leq k_0$, and LSQR is stopped once
its semi-convergence is practically identified.
\end{remark}

\begin{remark}\label{mildre}
For mildly ill-posed problems, the situation is fundamentally different.
We have
$\sqrt{\frac{k^2}{4\alpha^2-1}+\frac{k}{2\alpha-1}}>1$ and
$|L_{k_1}^{(k)}(0)|>1$ considerably as $k$ increases up to $k_0$
because of $\frac{1}{2}<\alpha\leq 1$,
leading to $\eta_k>1$ substantially. This means that for some
iterations $k\leq k^*$,
$\gamma_k$ may be substantially bigger than $\sigma_{k+1}$ and can
well lie between $\sigma_k$ and $\sigma_1$. In this case,
the rank $k$ approximation $P_kB_kQ_k^T$ is considerably less accurate
than the best rank $k$ approximation $A_k$.
\end{remark}

\begin{remark}
This theorem is different from Theorem 2.7 in
\cite{huangjia}. There are several subtle treatments
in the proof of Theorem~\ref{main1}, and ignoring or missing any one
of them would make it impossible to obtain accurate
estimates for $\gamma_k$. The first is to treat
$\|U_k\Sigma_kV_k^T(I-\hat{Z}_k\hat{Z}_k^T)\|$ as a whole.
If we amplified it by
$$
\|U_k\Sigma_kV_k^T(I-\hat{Z}_k\hat{Z}_k^T)\|
\leq \|\Sigma_k\|\|V_k^T(I-\hat{Z}_k\hat{Z}_k^T)\|=
\sigma_1\|\sin\Theta(\mathcal{V}_k,\mathcal{V}_k^R)\|,
$$
as in the proof of Theorem 2.7 in \cite{huangjia},
we would obtain a too large overestimate, which is almost a
constant for severely
ill-posed problems for $k=1,2,\ldots,k^*$ and increases with
$k$ for moderately and mildly ill-posed problems.
The second is the use of
\eqref{inden1} and \eqref{inden2}. The third is the extraction of
$\|\Sigma_k\Delta_k^T\|$ from \eqref{separa} other than
amplify it to $\|\Sigma_k\|\|\Delta_k\|=\sigma_1\|\Delta_k\|$.
The fourth is accurate estimates for $\|\Sigma_k\Delta_k^T\|$;
see \eqref{prodnorm} and \eqref{prodnorm2} in Theorem~\ref{thm3}.
For example, without using \eqref{inden1} and \eqref{inden2}, we would
obtain
\begin{align*}
\epsilon_k^2 &\leq \|\Sigma_k\|^2\|(I+
 \Delta_k^T\Delta_k)^{-1}\Delta_k^T\Delta_k\|^2+\|\Sigma_k\|^2\|(I+
 \Delta_k^T\Delta_k)^{-1}\Delta_k^T\|^2\\
 &=\sigma_1^2\left(\frac{\|\Delta_k\|^2}{1+\|\Delta_k\|^2}\right)^2+\sigma_1^2
 \|(I+\Delta_k^T\Delta_k)^{-1}\Delta_k^T\|^2.
\end{align*}
From \eqref{noncomp} and the previous
estimates for $\|\Delta_k\|$, such bound is too pessimistic, and
it even does not decrease and become small as $k$ increases, while
our estimates for $\epsilon_k=\eta_k\sigma_{k+1}$
in Theorem~\ref{main1} are optimal and can decay swiftly with $k$.
\end{remark}

As it will turn out in the next section, there are intimate
relationships between the quality of the rank $k$ approximation
$P_{k+1}B_kQ_k^T$ and the approximation behavior of the Ritz values
$\theta_i^{(k)},\ i=1,2,\ldots,k$. To be precise, we introduce
the following definition of a near best rank $k$ approximation
to $A$: By definition \eqref{gammak}, the rank $k$ matrix $P_{k+1}B_kQ_k^T$
is called a near best rank $k$ approximation to $A$ if it satisfies
\begin{equation}\label{near}
\sigma_{k+1}\leq \gamma_k<\frac{\sigma_k+\sigma_{k+1}}{2},
\end{equation}
that is, $\gamma_k$ lies between $\sigma_{k+1}$ and $\sigma_k$ and is closer to
$\sigma_{k+1}$.

Based on Theorem~\ref{main1}, for the severely and moderately or mildly ill-posed
problems with the model singular values $\sigma_i=\zeta\rho^{-i}$ and
$\sigma_i=\zeta i^{-\alpha}$, $i=1,2,\ldots,n$ where $\rho>1$ and
$\alpha>\frac{1}{2}$, we next derive the sufficient conditions on
$\rho$ and $\alpha$ that guarantee that $P_{k+1}B_kQ_k^T$ is a near best rank $k$
approximation to $A$ for $k=1,2,\ldots,k^*$.

\begin{theorem}\label{nearapprox}
For a given \eqref{eq1}, $P_{k+1}B_kQ_k^T$ is a near best rank $k$ approximation
to $A$ if
\begin{equation}\label{condition}
\sqrt{1+\eta_k^2}<\frac{1}{2}\frac{\sigma_k}{\sigma_{k+1}}+\frac{1}{2}.
\end{equation}
Furthermore, $P_{k+1}B_kQ_k^T$ is a near best rank $k$ approximation to $A$
if $\rho>2$ for the severely ill-posed problems with $\sigma_i=\zeta\rho^{-i}, \
i=1,2,\ldots,n$ or $\alpha$ satisfies
\begin{equation}\label{condition1}
2\sqrt{1+\eta_k^2}-1<\left(\frac{k+1}{k}\right)^{\alpha},\
k=1,2,\ldots,k^*
\end{equation}
for the moderately and mildly ill-posed problems with $\sigma_i=\zeta i^{-\alpha}$
and $\alpha>\frac{1}{2}$, $i=1,2,\ldots,n$.
\end{theorem}

{\em Proof}.
By \eqref{final}, we see that $\gamma_k\leq \sqrt{1+\eta_k^2}\sigma_{k+1}$.
Therefore, $P_{k+1}B_kQ_k^T$ is a near best rank $k$ approximation to $A$ in
the sense of \eqref{near} provided that
$$
\sqrt{1+\eta_k^2}\sigma_{k+1}<\frac{\sigma_k+\sigma_{k+1}}{2},
$$
from which \eqref{condition} follows.

From \eqref{etak0}, for the severely ill-posed problems with
$\sigma_i=\zeta\rho^{-i}, \ i=1,2,\ldots,n$ we have
\begin{equation}\label{simp}
\sqrt{1+\eta_k^2}<1+\frac{1}{2}\eta_k^2\leq 1+\frac{1}{2}\rho^{-2(1+\beta)}
<1+\rho^{-1}, \ k=1,2,\ldots,k^*.
\end{equation}
Since $\sigma_k/\sigma_{k+1}=\rho$, \eqref{condition} holds provided that
$$
1+\rho^{-1}<\frac{1}{2}\rho+\frac{1}{2},
$$
i.e., $\rho^2-\rho-2>0$, solving which for $\rho$ we get $\rho>2$. For the
moderately or mildly ill-posed problems with $\sigma_i=\zeta i^{-\alpha}$
and $\alpha>\frac{1}{2}$, $i=1,2,\ldots,n$, we have $\frac{\sigma_k}{\sigma_{k+1}}=\left(\frac{k+1}{k}\right)^{\alpha}$,
from which and \eqref{condition} it is direct to obtain \eqref{condition1}.
\qquad\endproof

\begin{remark}
For severely ill-posed problems with the model singular values
$\sigma_i=\zeta\rho^{-i}$ and $k=1,2,\ldots,k^*$, one always obtains
a near best rank $k$ approximation $P_{k+1}B_kQ_k^T$ to $A$ provided $\rho>2$.
\end{remark}

\begin{remark}\label{apprbe}
For the moderately ill-posed problems with $\sigma_i=\zeta i^{-\alpha}$,
on the one hand, for each fixed $k\leq k^*$,
there must be $\alpha>1$ such that \eqref{condition1} holds since,
by \eqref{const2}, its left-hand side decreases
to zero and the right-hand side increases to infinity with respect to $\alpha$;
the bigger $k$ is, the bigger $\alpha>1$ is required. Therefore,
there exists a suitable $\alpha>1$ that guarantees that $P_{k+1}B_kQ_k^T$
is a near best rank $k$ approximation to $A$ for all $k\leq k^*$.
On the other hand, for a {\em given} problem \eqref{eq1}
with $\alpha>1$ fixed, the smaller $k$ is, the more easily \eqref{condition1}
is met as the left-hand side decreases and the right-hand side
increases with respect to $k$. As a consequence, it is more possible
that $P_{k+1}B_kQ_k^T$ is a near best
rank $k$ approximations to $A$ for $k$ small and it
may not be a near best rank $k$ approximation at some
iterations $k\leq k^*$ if $\alpha$ is not big enough.
\end{remark}

\begin{remark}
For the mildly ill-posed problems with $\sigma_i=\zeta i^{-\alpha}$,
Theorem~\ref{moderate} has shown that
$\|\Delta_k\|$ is generally not small and can be large
as $k$ increases up to $k^*$. From \eqref{etadelta2}, we see
that the size of $\eta_k$ is comparable to $\|\Delta_k\|$. Note that
$\left(\frac{k+1}{k}\right)^{\alpha}\leq 2$ for
$\frac{1}{2}<\alpha\leq 1$ and any $k\geq 1$.
Consequently, \eqref{condition1} may be satisfied only for $k$ very small
and $\alpha$ not close to $\frac{1}{2}$, and it
cannot be met generally as $k$ increases.
Hence $P_{k+1}B_kQ_k^T$ may be a near best rank $k$ approximation
to $A$ no longer soon as $k$ increases.
\end{remark}

\section{The approximation properties of the Ritz values $\theta_i^{(k)}$}
\label{ritzapprox}

In this section, we make an in-depth analysis on the approximation
behavior of the Ritz values $\theta_i^{(k)}$. This problem
has been highly concerned since the use of LSQR in the context of
ill-posed problems, but has remained open.

Starting with Theorem~\ref{main1}, we prove that, under
certain sufficient conditions on $\rho$ and $\alpha$ for the severely
and moderately ill-posed problems with $\sigma_i=\zeta\rho^{-i}$ and
$\sigma_i=\zeta i^{-\alpha}$, respectively,
the $k$ Ritz values $\theta_i^{(k)}$ approximate the
large singular values $\sigma_i$ of $A$ in natural order
for $k=1,2,\ldots,k^*$,
which means that the semi-convergence of LSQR occurs
at iteration $k^*=k_0$ and no Ritz value smaller than
$\sigma_{k_0+1}$ appears before $k\leq k_0$.
Combining this result with Theorem~\ref{nearapprox},
we come to the definitive conclusion that LSQR has the full
regularization for these two kinds of problems for suitable
$\rho>1$ and $\alpha>1$. On the other hand,
we will show that for some $k\leq k^*$
the Ritz values generally do not approximate
the large singular values of $A$ in natural order
for severely or moderately ill-posed problems with $\rho>1$ or
$\alpha>1$ not enough
and mildly ill-posed problems, which, by Theorem~\ref{semicon}, means
that $k^*<k_0$.

\begin{theorem}\label{ritzvalue}
Assume that \eqref{eq1} is severely ill-posed with
$\sigma_i=\zeta\rho^{-i}$ and $\rho>1$ or moderately and
mildly ill-posed with $\sigma_i=\zeta i^{-\alpha}$, $i=1,2,\ldots,n$, and
let the Ritz values $\theta_i^{(k)}$ be labeled
as $\theta_1^{(k)}>\theta_2^{(k)}>\cdots>\theta_{k}^{(k)}$.
Then
\begin{align}
0<\sigma_i-\theta_i^{(k)} &\leq \sqrt{1+\eta_k^2}\sigma_{k+1},\
i=1,2,\ldots,k.\label{error}
\end{align}
For $k=1,2,\ldots, k^*$, if $\rho\geq 1+\sqrt{2}$ or $\alpha$ satisfies
\begin{equation}\label{condm}
1+\sqrt{1+\eta_{k}^2}<\left(\frac{k+1}{k}\right)^{\alpha},
\end{equation}
then the $k$ Ritz values $\theta_i^{(k)}$ strictly interlace
the first large $k+1$ singular values of $A$:
\begin{align}
\sigma_{i+1}&<\theta_i^{(k)}<\sigma_i,\,i=1,2,\ldots,k,
\label{error2}
\end{align}
which means that the $\theta_i^{(k)}$ approximate the first $k$ large $\sigma_i$ in
natural order.
\end{theorem}

{\em Proof}.
Note that, for $k=1,2,\ldots,k^*$, the $k$ Ritz values $\theta_i^{(k)}$ are
just the nonzero singular values of $P_{k+1}B_kQ_k^T$, whose other $n-k$
singular values are zeros. We write
$$
A=P_{k+1}B_k Q_k^T+(A-P_{k+1}B_k Q_k^T).
$$
Since $\|A-P_{k+1}B_k Q_k^T\|=\gamma_k$, by the Mirsky's theorem of
singular values \cite[p.204, Theorem 4.11]{stewartsun} we have
\begin{equation}\label{errbound}
| \sigma_i-\theta_i^{(k)}|\leq \gamma_k\leq
\sqrt{1+\eta_k^2}\sigma_{k+1},\ i=1,2,\ldots,k.
\end{equation}
Since the singular values of $A$ are simple and $b$ has components in all the
left singular vectors $u_1,u_2,\ldots, u_n$ of $A$, Lanczos bidiagonalization
can be run to completion without breakdown,
producing $P_{n+1},\ Q_n$ and
the lower bidiagonal $B_n\in \mathbb{R}^{(n+1)\times n}$ such that
\begin{equation}\label{fulllb}
P^TAQ_n=\left(\begin{array}{c}
B_n\\
\mathbf{0}
\end{array}
\right)
\end{equation}
with $P=(P_{n+1},\hat{P}) \in \mathbb{R}^{m\times m}$ and
$Q_n\in \mathbb{R}^{n\times n}$ being orthogonal and
the diagonals $\alpha_i$ and subdiagonals
$\beta_{i+1}$, $i=1,2,\ldots,n$, of $B_n$ being positive.\footnote{
If $m=n$, it is easily justified that $\beta_{n+1}=0$, Lanczos
bidiagonalization
generates the orthogonal matrices $P_n$, $Q_n$ and the $n\times n$ lower
bidiagonal $B_n$ with the positive diagonals $\alpha_i$ and subdiagonals
$\beta_i$. This does not affect
the derivation and results followed, and we only need to replace
$P_{n+1}$ by $P_n$.} Notice that the singular values of $B_k,\ k=1,2,\ldots,n,$
are all simple and $B_k$ consists of the first $k$ columns of $B_n$
with the last $n-k$ {\em zero} rows deleted. Applying the Cauchy's {\em strict}
interlacing theorem \cite[p.198, Corollary 4.4]{stewartsun} to the singular
values of $B_k$ and $B_n$, we have
\begin{align}
\sigma_{n-k+i}< \theta_i^{(k)}&< \sigma_i,\ i=1,2,\ldots,k.
\label{interlace}
\end{align}
Therefore, \eqref{errbound} becomes
\begin{equation}\label{ritzapp}
0< \sigma_i-\theta_i^{(k)}\leq\gamma_{k}\leq
\sqrt{1+\eta_k^2}\sigma_{k+1},\ i=1,2,\ldots,k,
\end{equation}
which proves \eqref{error}.

For $i=1,2,\ldots,k$, notice that $\rho^{-k+i}\leq 1$. Then
from \eqref{ritzapp}, \eqref{simp} and $\sigma_i=\zeta\rho^{-i}$ we obtain
\begin{align*}
\theta_i^{(k)}&\geq \sigma_i-\gamma_k>\sigma_i-
(1+\rho^{-1})\sigma_{k+1}\\
&=\zeta\rho^{-i}-\zeta (1+\rho^{-1})\rho^{-(k+1)}\\
&=\zeta\rho^{-(i+1)}(\rho-(1+\rho^{-1})\rho^{-k+i})\\
&\geq \zeta\rho^{-(i+1)}(\rho-(1+\rho^{-1}))\\
&\geq\zeta\rho^{-(i+1)}=\sigma_{i+1},
\end{align*}
provided that $\rho-(1+\rho^{-1})=\rho-\rho^{-1}-1\geq 1$. Solving the inequality
gives $\rho\geq 1+\sqrt{2}$. Together with the upper bound
of \eqref{interlace}, we have proved \eqref{error2}.

For the moderately and mildly ill-posed problems with
$\sigma_i=\zeta i^{-\alpha},\ i=1,2,\ldots,n$ and $k=1,2,\ldots,k^*$,
we get
\begin{align*}
\theta_i^{(k)}&\geq \sigma_i-\gamma_k\geq\sigma_i-\sqrt{1+\eta_k^2}
\sigma_{k+1}\\
&=\zeta i^{-\alpha}-\zeta \sqrt{1+\eta_k^2}(k+1)^{-\alpha}\\
&=\zeta (i+1)^{-\alpha}\left(\left(\frac{i+1}{i}\right)^{\alpha}
-\sqrt{1+\eta_k^2}\left(\frac{i+1}{k+1}\right)^{\alpha}\right)\\
&>\zeta (i+1)^{-\alpha}=\sigma_{i+1},
\end{align*}
i.e., \eqref{error2} holds, provided that $\eta_k$ and $\alpha>1$
satisfy
$$
\left(\frac{i+1}{i}\right)^{\alpha}
-\sqrt{1+\eta_k^2}\left(\frac{i+1}{k+1}\right)^{\alpha}>1,
$$
which means that
$$
\sqrt{1+\eta_k^2}<\left(\left(\frac{i+1}{i}\right)^{\alpha}-1\right)
\left(\frac{k+1}{i+1}\right)^{\alpha}=
\left(\frac{k+1}{i}\right)^{\alpha}-\left(\frac{k+1}{i+1}\right)^{\alpha},\
i=1,2,\ldots,k.
$$
It is easily justified that the above right-hand side monotonically
decreases with respect to $i=1,2,\ldots,k$, whose minimum
attains at $i=k$ and equals $\left(\frac{k+1}{k}\right)^{\alpha}-1$,
from which we obtain the condition \eqref{condm}.
\qquad\endproof

\begin{remark}
For each $k\leq k^*$, since the left-hand side of \eqref{condm} tends to
two and its right-side tends to infinity with respect to $\alpha>1$,
there must be $\alpha>1$ such that \eqref{condm} holds.
Comparing Theorem~\ref{nearapprox} with Theorem~\ref{ritzvalue}, we
find out that, as far as the severely and moderately ill-posed problems are
concerned, Lanczos bidiagonalization generates the near best rank $k$
approximations $P_{k+1}B_kQ_k^T$ to $A$, and the singular values
$\theta_i^{(k)}$ of $B_k$ approximate the
large singular values $\sigma_i$ of $A$ in natural order for suitable
$\rho>1$ and $\alpha>1$ for $k=1,2,\ldots,k^*$.
On the other hand, for a {\em given} problem with $\alpha>1$,
the smaller $k$ is, the more easily the condition \eqref{condm} is met,
and thus the more possible is for the $\theta_i^{(k)}$ to
approximate the large singular values $\sigma_i$ in natural order.
The $\theta_i^{(k)}$ may not approximate the large singular values $\sigma_i$
in natural order at some iterations $k\leq k^*$ if $\alpha>1$ is not enough.
\end{remark}

\begin{remark}
For the mildly ill-posed problems $\sigma_i=\zeta i^{-\alpha}$,
the sufficient condition \eqref{condm} for \eqref{error2}
is never met because its left-hand side is always bigger than two but the
right-hand side
$
\left(\frac{k+1}{k}\right)^{\alpha}\leq 2
$
for any $k\geq 1$ and $\frac{1}{2}< \alpha\leq 1$. This indicates that
it is hard that the $\theta_i^{(k)}$ approximate
the $k$ large singular values $\sigma_i$ in natural order.
\end{remark}

\section{Monotonicity of $\gamma_k$, decay rates of $\alpha_k$ and $\beta_{k+1}$ and
their practical importance}
\label{alphabeta}

In this section, we present a number of results on $\gamma_k$ and
the decay rates of
$\alpha_k$ and $\beta_{k+1}$, and highlight their implications and practical
importance.

\begin{theorem}\label{main2}
With the notation defined previously, the following results hold:
\begin{eqnarray}
  \alpha_{k+1}&<&\gamma_k\leq \sqrt{1+\eta_k^2}\sigma_{k+1},
  \ k=1,2,\ldots,n-1,\label{alpha}\\
 \beta_{k+2}&<& \gamma_k\leq\sqrt{1+\eta_k^2}\sigma_{k+1}, \ k=1,2,\ldots,n-1,
 \label{beta}\\
\alpha_{k+1}\beta_{k+2}&\leq &
\frac{\gamma_k^2}{2}\leq
   \frac{(1+\eta_k^2)\sigma_{k+1}^2}{2}, \ k=1,2,\ldots,n-1,
  \label{prod2}\\
\gamma_{k+1}&<&\gamma_k,\  \ k=1,2,\ldots,n-2. \label{gammamono}
\end{eqnarray}
\end{theorem}

{\em Proof}.
From \eqref{fulllb}, since $P$ and $Q_n$ are orthogonal matrices,
we have
\begin{align}
\gamma_k &=\|A-P_{k+1}B_kQ_k^T\|=\|P^T(A-P_{k+1}B_kQ_k^T)Q_n\| \label{invar}\\
&=
\left\| \left(\begin{array}{c}
B_n \\
\mathbf{0}
\end{array}
\right)-(I,\mathbf{0} )^TB_k (I,\mathbf{0} )\right\|=\|G_k\| \label{gk}
\end{align}
with
\begin{align}\label{gk1}
G_k&=\left(\begin{array}{cccc}
\alpha_{k+1} & & & \\
\beta_{k+2}& \alpha_{k+2} & &\\
&
  \beta_{k+3} &\ddots & \\& & \ddots & \alpha_{n} \\
  & & & \beta_{n+1}
  \end{array}\right)\in \mathbb{R}^{(n-k+1)\times (n-k)}
\end{align}
resulting from deleting the $(k+1)\times k$ leading principal matrix $B_k$
of $B_n$ and the first $k$ zero rows and columns of the resulting matrix.
For $k=1,2,\ldots,n-1$ we have
\begin{align}
\alpha_{k+1}^2+\beta_{k+2}^2&=\|G_ke_1\|^2\leq \|G_k\|^2=\gamma_k^2,
\label{alphabetasum1}
\end{align}
which proves the lower bounds in \eqref{alpha} and \eqref{beta}
since $\alpha_{k+1}>0$ and $\beta_{k+2}>0$.
Furthermore, from \eqref{final}, we obtain the upper bounds in \eqref{alpha}
and \eqref{beta}. Noting that
\begin{align*}
2\alpha_{k+1}\beta_{k+2}&\leq \alpha_{k+1}^2+\beta_{k+2}^2= \gamma_k^2,
\end{align*}
we prove \eqref{prod2}.

By $\gamma_k=\|G_k\|$ and \eqref{gk1}, observe that $\gamma_{k+1}=\|G_{k+1}\|$
equals the 2-norm of the submatrix deleting the first column of $G_k$.
Applying the Cauchy's strict interlacing theorem to the singular values
of this submatrix and $G_k$, we obtain \eqref{gammamono}.
\qquad\endproof

\begin{remark}
The strict decreasing property \eqref{gammamono} of $\gamma_k$ and
the lower bounds for $\gamma_k$ in \eqref{alpha}--\eqref{prod2} hold
unconditionally for any general $A$,
independent of the degree of ill-posedness of \eqref{eq1}.
\end{remark}

\begin{remark}
It is generally impractical to compute $\gamma_k$ for $A$ large.
However, the proofs of \eqref{alpha} and \eqref{beta} indicate that
$\alpha_{k+1}+\beta_{k+2}$ decays as fast as $\gamma_k$. Hence,
strikingly, we can reliably judge the decay rates of $\gamma_k$
by those of $\alpha_{k+1}+\beta_{k+2}$ during computation without extra cost.
\end{remark}

\begin{remark}\label{severmod}
For severely and moderately ill-posed problems with suitable $\rho>1$ and
$\alpha>1$, based on the previous results,
\eqref{alpha} and \eqref{beta} show that $\alpha_{k+1}+\beta_{k+2}$
decays as fast as $\sigma_{k+1}$. For mildly ill-posed problems,
since generally $\eta_k>1$ considerably,
$\alpha_{k+1}+\beta_{k+2}$ decays more slowly than
$\sigma_{k+1}$.
\end{remark}


\section{Numerical experiments}\label{numer}

We choose several 1D and 2D ill-posed problems
and a random ill-posed problem with the prescribed singular values in the
toolboxes \cite{gazzola18,hansen07,hansen12}. Table~\ref{tab1} lists
the test problems, which take default parameter(s) and
include severely, moderately and mildly
ill-posed problems  as well as the ones with the singular
values decaying more slowly than those of the mildly ill-posed problems
for {\em given} $m$ and $n$. The random mildly ill-posed problem
{\sf regutm.m} \cite{hansen95,hansen07} is constructed with
the prescribed singular values $\sigma_i=i^{-0.6}$
and the left and right singular vectors $u_i,\ v_i,
\, i=1,2,\ldots,n$ having exactly $i-1$ sign changes, and we
set $x_{true}=ones(n,1)$ and generate the noise-free
right hand side $b_{true}=Ax_{true}$.

We mention that
it is hard to find available 2D real-life suitable severely and moderately
ill-posed problems for justifying our results.
Gazzola, Hansen and Nagy \cite{gazzola18} have
presented a number of 2D test problems, where
the image deblurring problem {\sf PRblurgauss}, the inverse diffusion problem
{\sf PRdiffusion} and the nuclear magnetic resonance (NMR) relaxometry
problem {\sf PRnmr} are severely ill-posed. But the latter two matrices are
only available as a function handle, for which we cannot compute their SVDs.
Setting the parameter {\sf options.BlurLevel='severe'}, we
have computed the SVD of {\sf PRblurgauss} with $m=n=10000$ and
found $\sigma_1/\sigma_{1500}\approx 1.99\times 10^{14}=
\mathcal{O}(\frac{1}{\epsilon_{\rm mach}})$.
Unfortunately, we have found out that about half of the first 1500
large singular values are genuinely or numerically {\em multiple}.
For example, among the first 40 singular values,
$\sigma_3,\sigma_6, \sigma_8, \sigma_{11}, \sigma_{13}, \sigma_{15}, \sigma_{17},
\sigma_{19}, \sigma_{22}, \sigma_{24}, \sigma_{25}, \sigma_{26}, \sigma_{28},
\sigma_{30}, \sigma_{33}, \sigma_{35},
\sigma_{37}$ and $\sigma_{39}$ are multiple. Therefore, these
problems are either unsuitable or inaccessible for our propose.
Meanwhile, there is no
2D moderately ill-posed problem, and all the other problems are mildly
ill-posed in \cite{gazzola18}.
We will test the 2D mildly ill-posed problems {\sf PRblurrotation}
and {\sf PRblurspeckle} of $m=n=14400$, which simulate a spatially
variant rotational motion blur around the center of the image and
spatially invariant blurring caused by atmospheric
turbulence, respectively.

For our propose of justifying the sharpness of our estimates,
it is enough to test any severely, moderately and mildly ill-posed
problems. In the meantime, for each test problem we
compare the accuracy of the best LSQR regularized solution
$x_{k^*}^{lsqr}$ with that of the best TSVD solution $x_{k_0}^{tsvd}$.
In addition, we attempt
to show that for the 2D problems {\sf blur}, {\sf fanbeamtomo} and
{\sf paralleltomo} whose singular values decay more
slowly than those of a mildly ill-posed problem for given $m$ and $n$, LSQR
and the TSVD method behave as if these problems were ordinary ones
if the noise level
$$
\varepsilon=\frac{\|e\|}{\|b_{true}\|}
$$
is fairly small, e.g., $10^{-3}$,
that is, the noise $e$ does not play any effect
in regularization, and the best regularized solution by the
TSVD method is simply $x_n^{tsvd}=x_{naive}=A^{\dagger}b$,
and LSQR works in its regular way and ultimately finds a converged approximation
to $x_{naive}$. In this case, no semi-convergence
occurs, and the
discrete ill-posed problems are actually ordinary ones; see
the elaboration in the introduction. However, we will report that if
$\varepsilon$ is relatively larger, e.g., 0.05, then
the semi-convergence of LSQR and the TSVD method occurs.

Keep in mind that for an underlying linear compact operator equation
$Kf=g$, the singular values of the operator $K$ must decay at least as fast as
$\mathcal{O}(k^{-\alpha}),\ k=1,2,\ldots,\infty$ with
$\alpha>\frac{1}{2}$; see, e.g., \cite{hanke93,hansen98}. Therefore,
provided that such a continuous problem is discretized finely enough
by a suitable scheme, the resulting linear discrete ill-posed
problem \eqref{eq1} will inherit such property, that is,
the singular values $\sigma_k,\ k=1,2,\ldots,n$ of $A$ ultimately
decay faster than $\mathcal{O}(k^{-1/2})$
when $n$ is sufficiently large; otherwise,
the solution norm $\|x_{true}\|$ of \eqref{eq1} is
unbounded when $n\rightarrow \infty$.
For the 2D image deblurring  problems {\sf blur} \cite{hansen07},
{\sf fanbeamtomo} and {\sf paralleltomo} \cite{hansen12},
we have found that their singular values $\sigma_k$ decay slowly
and considerably more slowly than $\mathcal{O}(k^{-\frac{1}{2}})$
even if $m$ and $n$ are {\em a few ten thousands}. Actually, the
$\frac{\sigma_1}{\sigma_n}$
lie between $\mathcal{O}(10)\sim \mathcal{O}(10^3)$, very modest!
We will report precise details later. This indicates that the discretizations
are not sufficiently fine. In Table~\ref{tab1}, we take
$n=150^2, \ 120^2,\ 100^2$ for {\sf blur, fanbeamtomo} and {\sf paralleltomo}.
For {\sf blur}, it means that the 2D rectangular domain is discretized into
$150^2$ cells; for {\sf fanbeamtomo} and {\sf paralleltomo}, they mean that
the 2D domains are divided into 120 and 100 equally spaced
intervals in both dimensions, which create $120^2$ and $100^2$ cells,
respectively. Obviously, the discretizations are not fine enough, though
the discrete problems are seemingly already large.

{\small
\begin{table}[htp]
    \centering
    \caption{The description of test problems}
    \begin{tabular}{llll}
     \hline
     Problem  & Description                                & Ill-posedness &
     $m\times n$\\
     \hline
     {\sf shaw}     & 1D image restoration model & severe & $10240\times 10240$\\
     {\sf gravity}  & 1D gravity surveying problem  & severe  & $10240\times 10240$\\
     {\sf heat}     & Inverse heat equation                      & moderate & $10240\times 10240$\\
     {\sf deriv2}   & Computation of second derivative           & moderate & $10000\times 10000$\\
     {\sf regutm}   & Random ill-posed problem    & mild &$10000\times 10000$\\
     {\sf PRblurrotation} & 2D image deblurring problem &mild & $14400\times 14400$\\
    {\sf blur}      & 2D image restoration                 &unknown  & $22500\times 22500$\\
   {\sf fanbeamtomo} & 2D fan-beam tomography problem  &unknown & $61200\times 14400$\\
    {\sf paralleltomo}& 2D tomography problem &unknown& $25380\times 10000$\\
\hline
   \end{tabular}
   \label{tab1}
\end{table}
}

For each example, we use the code of \cite{hansen07,hansen12}
to generate $A$, $x_{true}$ and $b_{true}$.
In order to simulate the noisy data, we generate Gaussian white noise
vectors $e$ such that the relative noise levels
$\varepsilon$ equal some prescribed values.
To simulate exact arithmetic,
we use full reorthogonalization in Lanczos bidiagonalization.
All the computations are carried out in Matlab R2017b on the
Intel Core i7-4790k with CPU 4.00 GHz processor and 16 GB RAM
with the machine precision
$\epsilon_{\rm mach}= 2.22\times10^{-16}$ under the Miscrosoft
Windows 8 64-bit system.

Our numerical experiments consist of three subsections. The first two
subsections are devoted to the justification of our results,
and the third subsection pay special attention to the regularization
behavior and ability of LSQR.

\subsection{The accuracy $\gamma_k$ of rank $k$ approximations and the
approximation properties of the Ritz values $\theta_i^{(k)}$}

\begin{figure}
\begin{minipage}{0.48\linewidth}
  \centerline{\includegraphics[width=6.0cm,height=4.5cm]{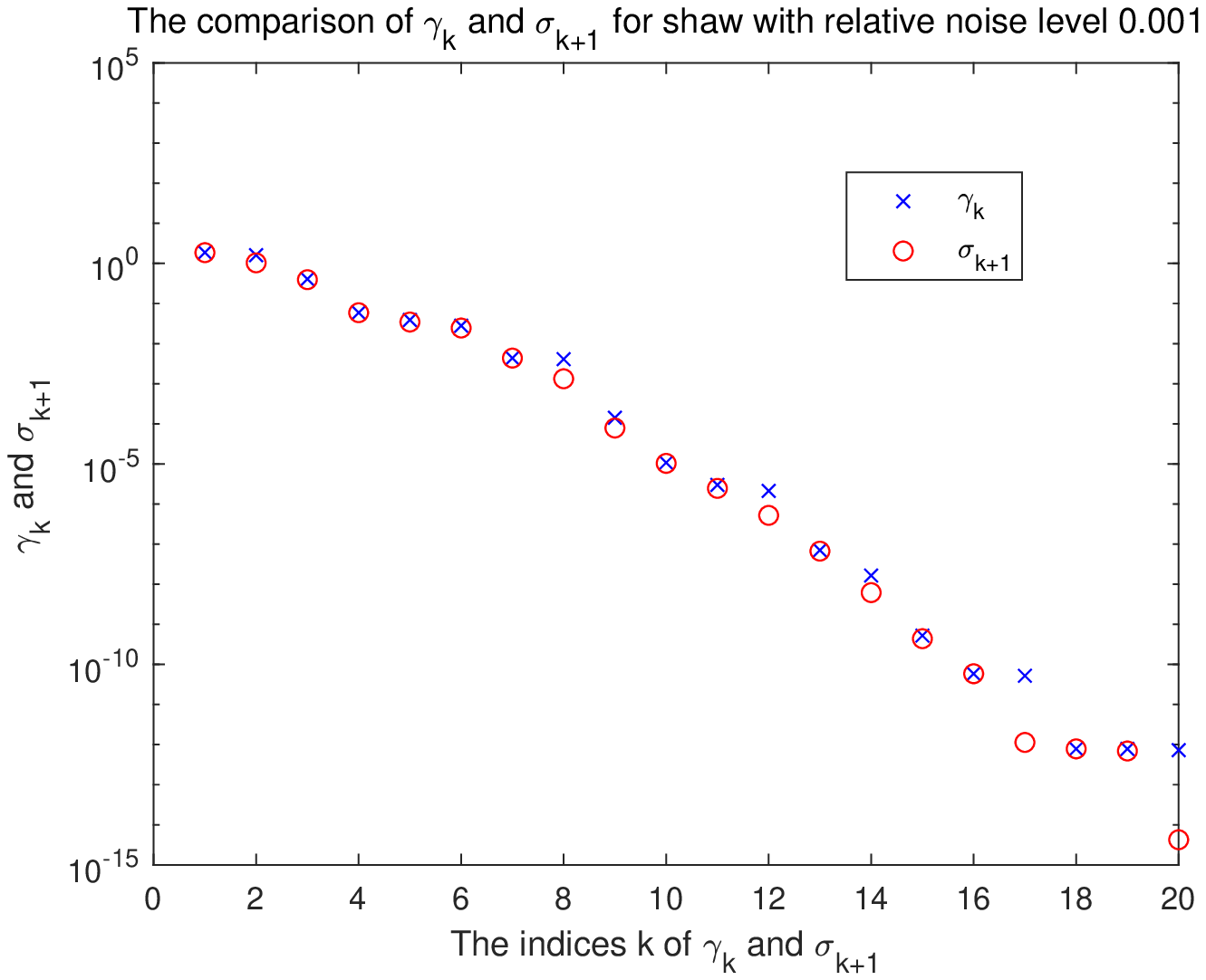}}
  \centerline{(a)}
\end{minipage}
\hfill
\begin{minipage}{0.48\linewidth}
  \centerline{\includegraphics[width=6.0cm,height=4.5cm]{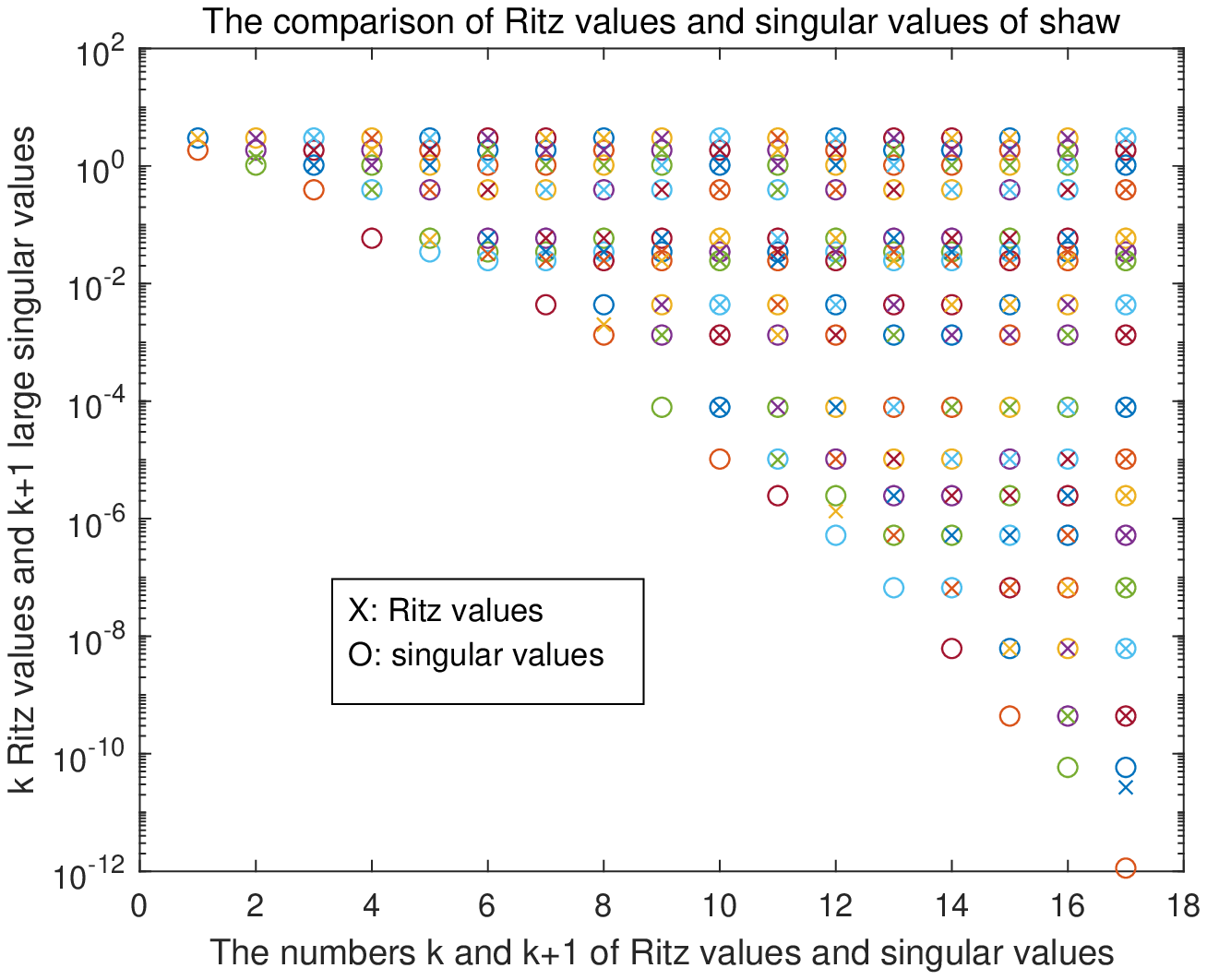}}
  \centerline{(b)}
\end{minipage}
\vfill
\begin{minipage}{0.48\linewidth}
  \centerline{\includegraphics[width=6.0cm,height=4.5cm]{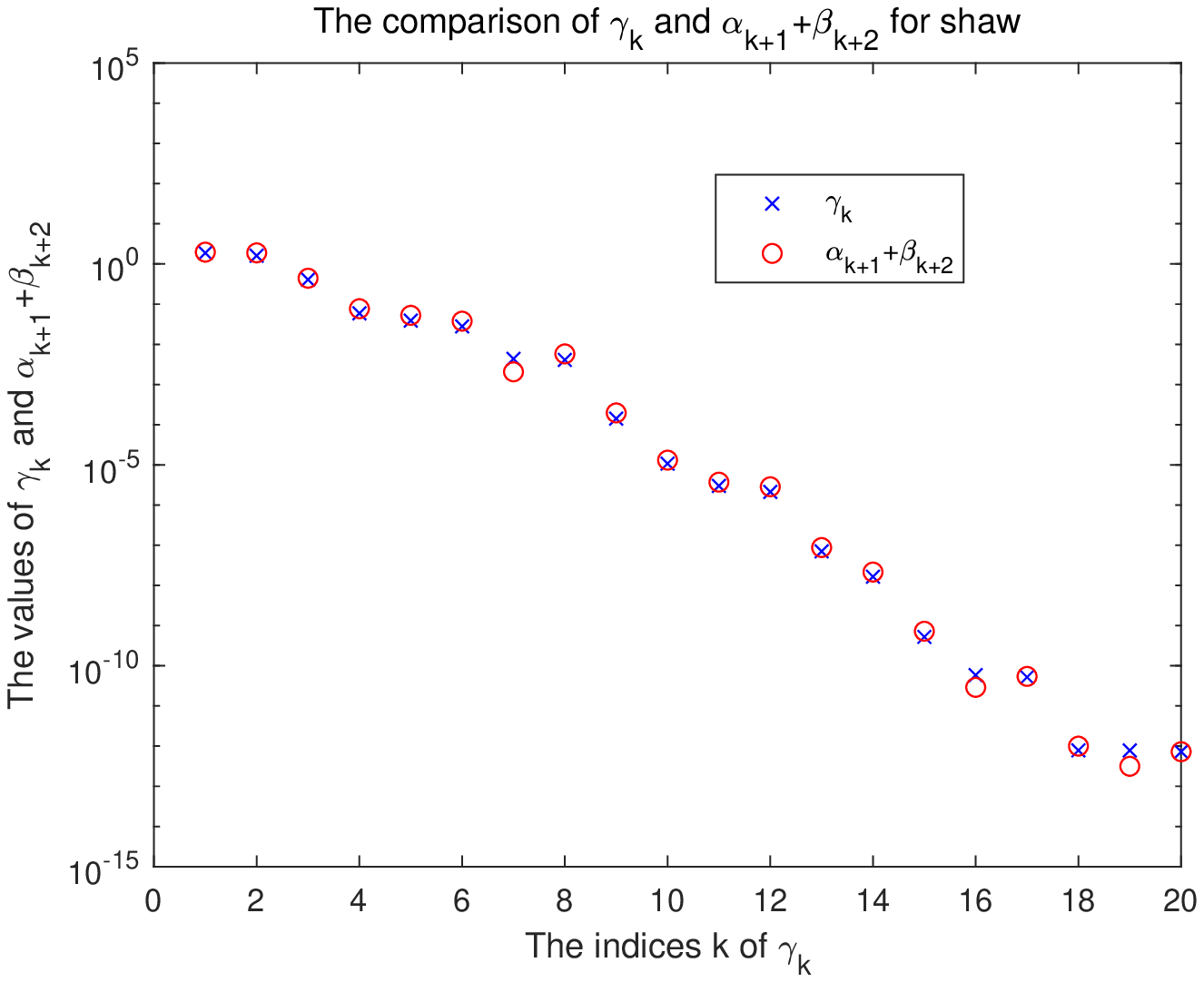}}
  \centerline{(c)}
\end{minipage}
\hfill
\begin{minipage}{0.48\linewidth}
  \centerline{\includegraphics[width=6.0cm,height=4.5cm]{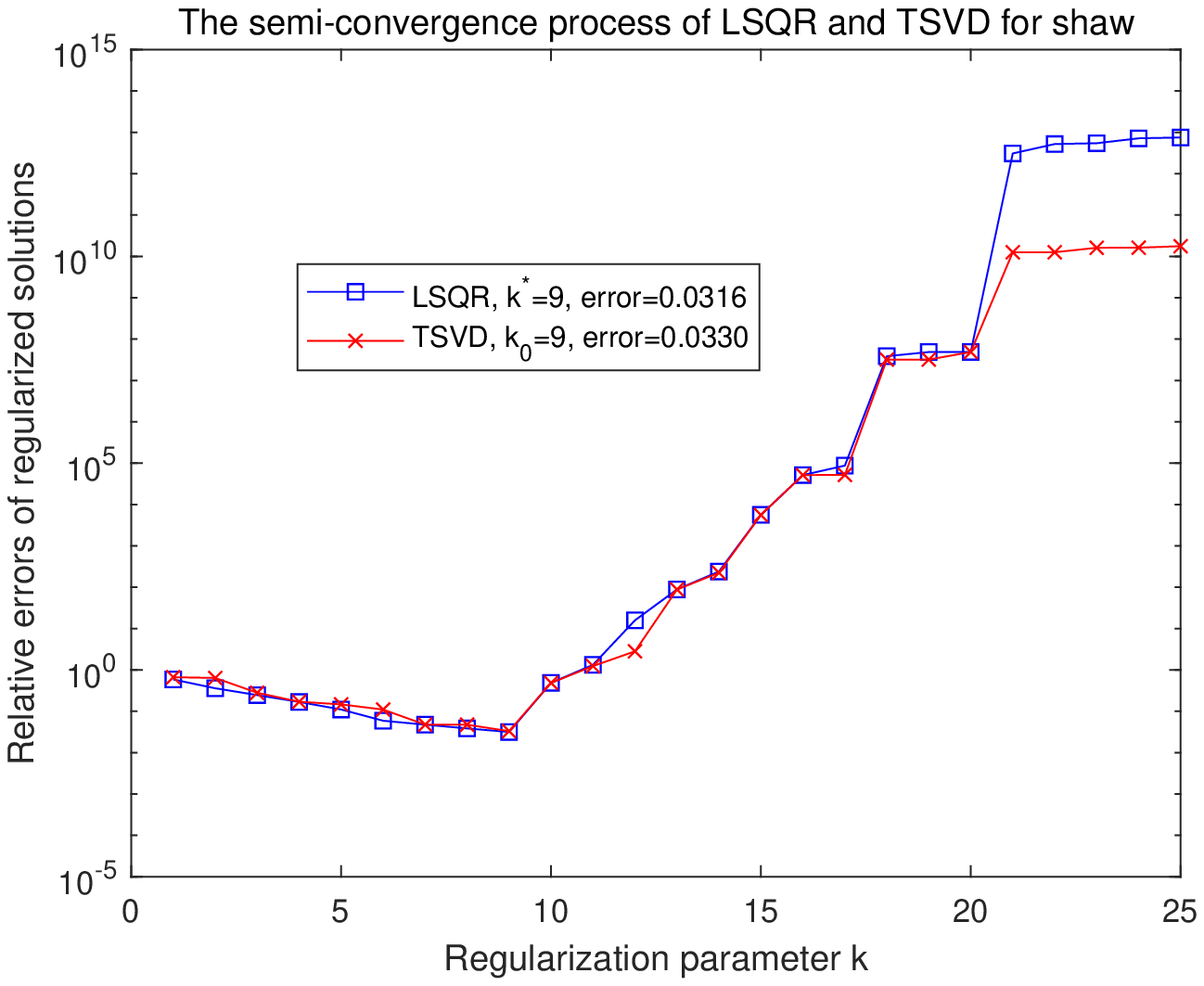}}
  \centerline{(d)}
\end{minipage}
\caption{{\sf shaw} of $n=10240$ with relative noise level $10^{-3}$.} \label{fig1}
\end{figure}

\begin{figure}
\begin{minipage}{0.48\linewidth}
  \centerline{\includegraphics[width=6.0cm,height=4.5cm]{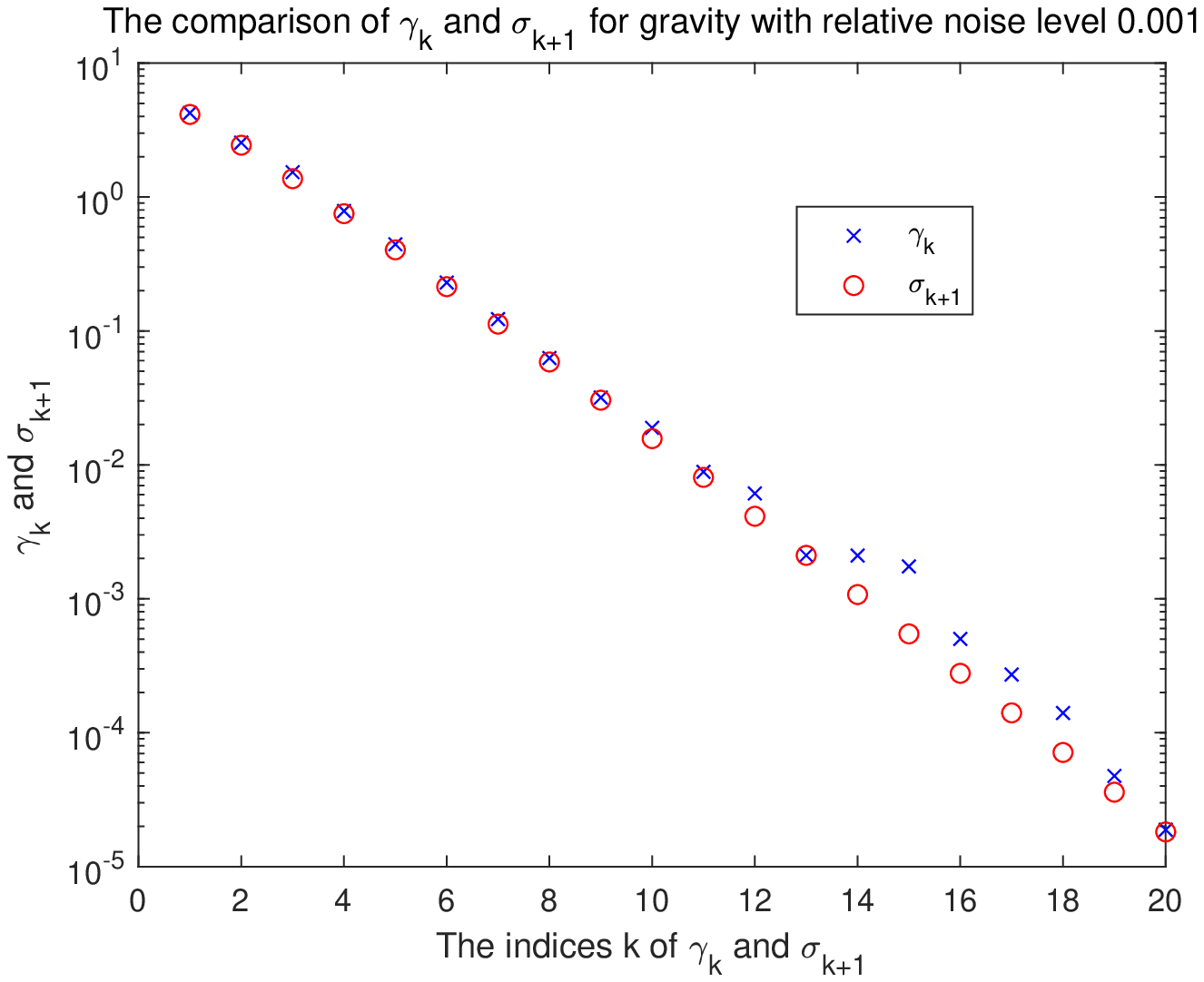}}
  \centerline{(a)}
\end{minipage}
\hfill
\begin{minipage}{0.48\linewidth}
  \centerline{\includegraphics[width=6.0cm,height=4.5cm]{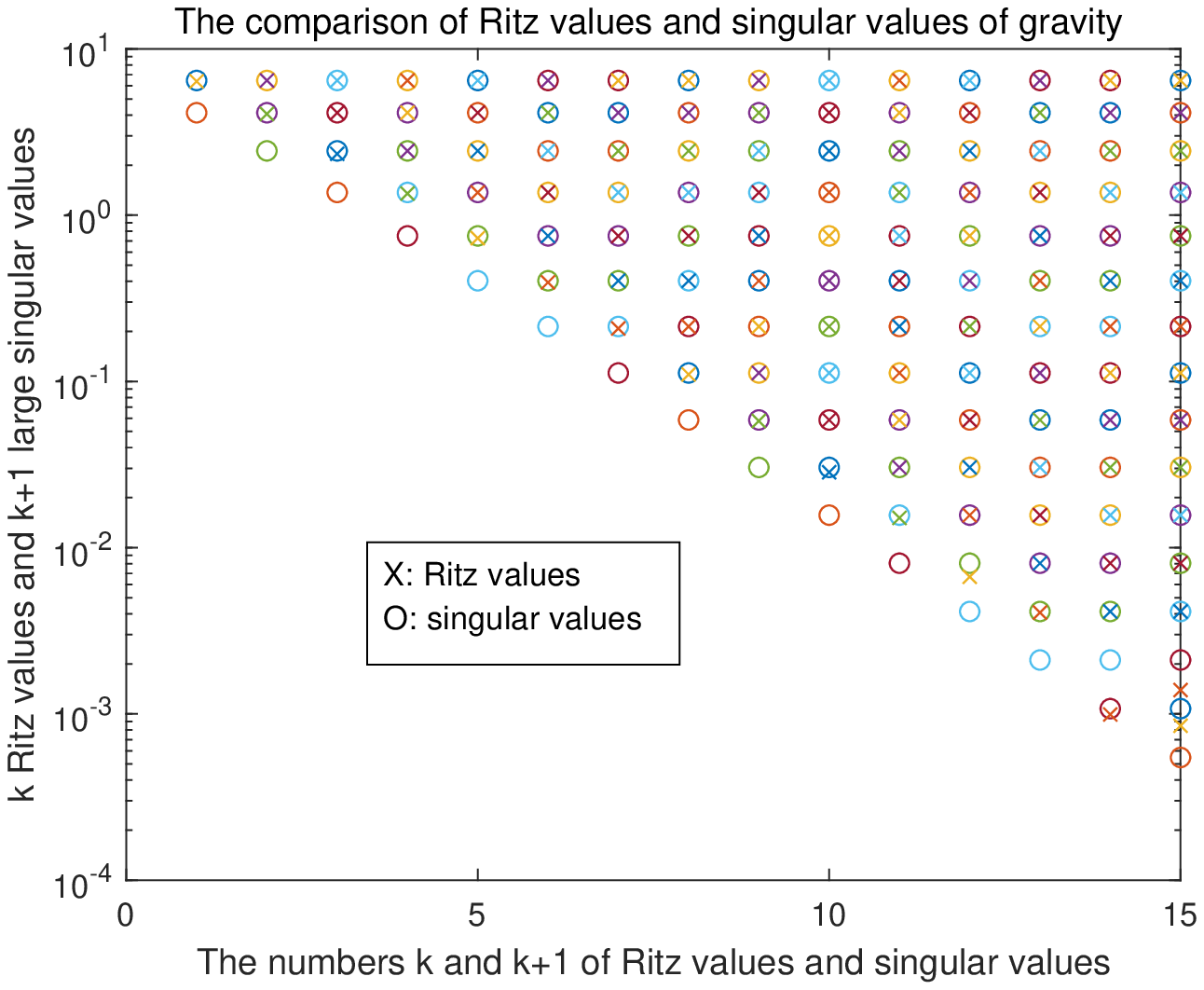}}
  \centerline{(b)}
\end{minipage}
\vfill
\begin{minipage}{0.48\linewidth}
  \centerline{\includegraphics[width=6.0cm,height=4.5cm]{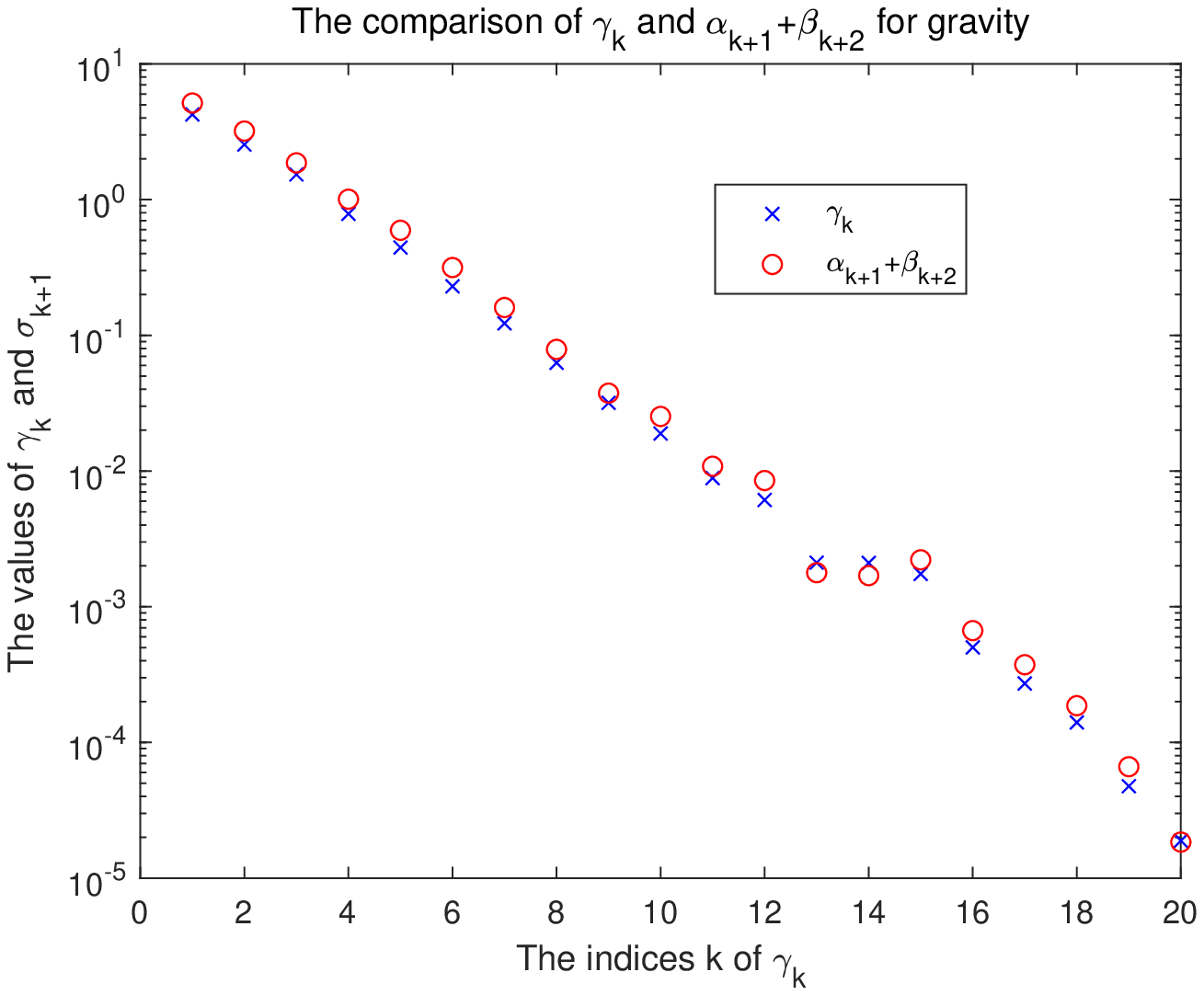}}
  \centerline{(c)}
\end{minipage}
\hfill
\begin{minipage}{0.48\linewidth}
  \centerline{\includegraphics[width=6.0cm,height=4.5cm]{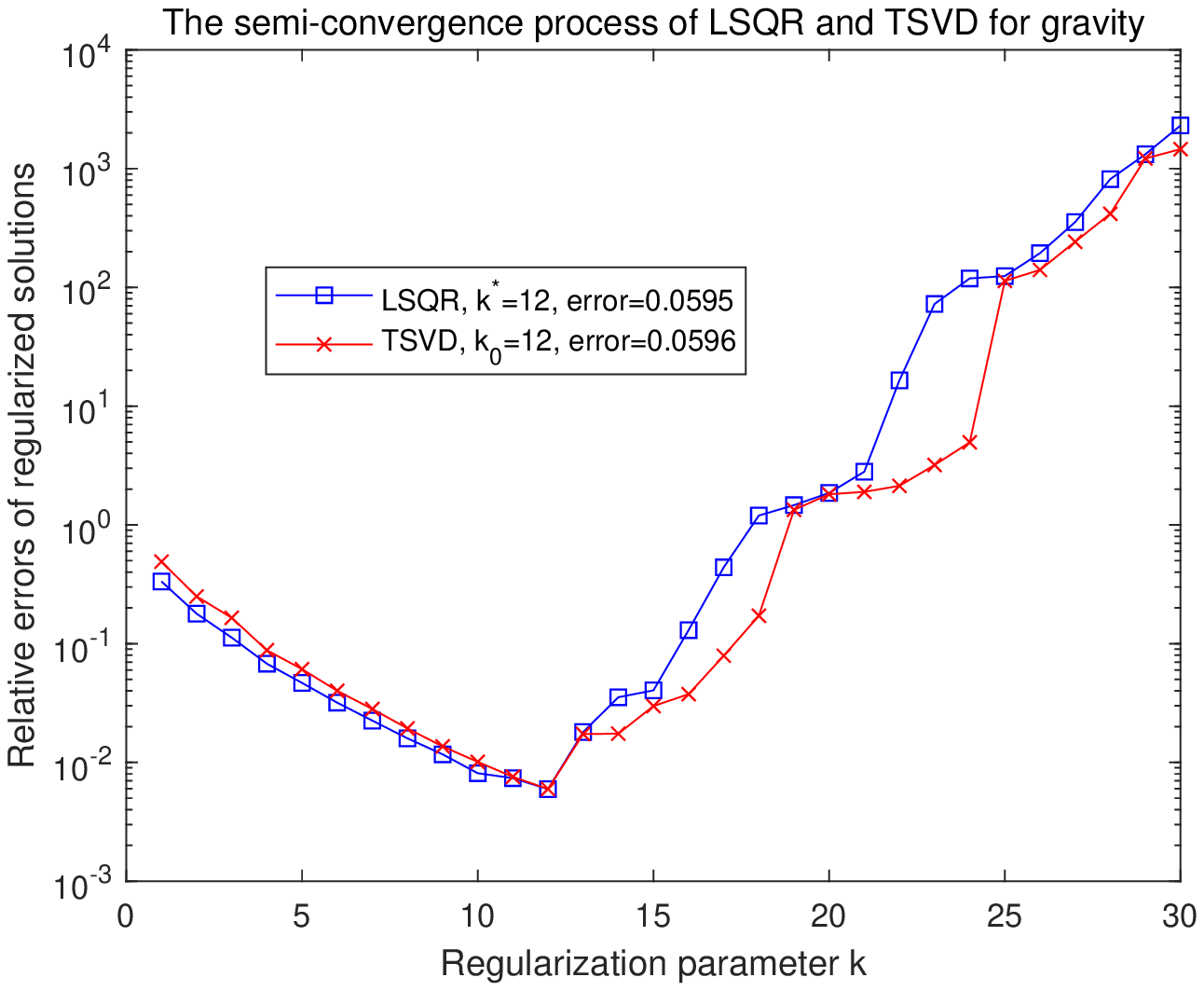}}
  \centerline{(d)}
\end{minipage}
\caption{{\sf gravity} of $n=10240$ with the relative noise level $10^{-3}$.} \label{fig2}
\end{figure}

\begin{figure}
\begin{minipage}{0.48\linewidth}
  \centerline{\includegraphics[width=6.0cm,height=4.5cm]{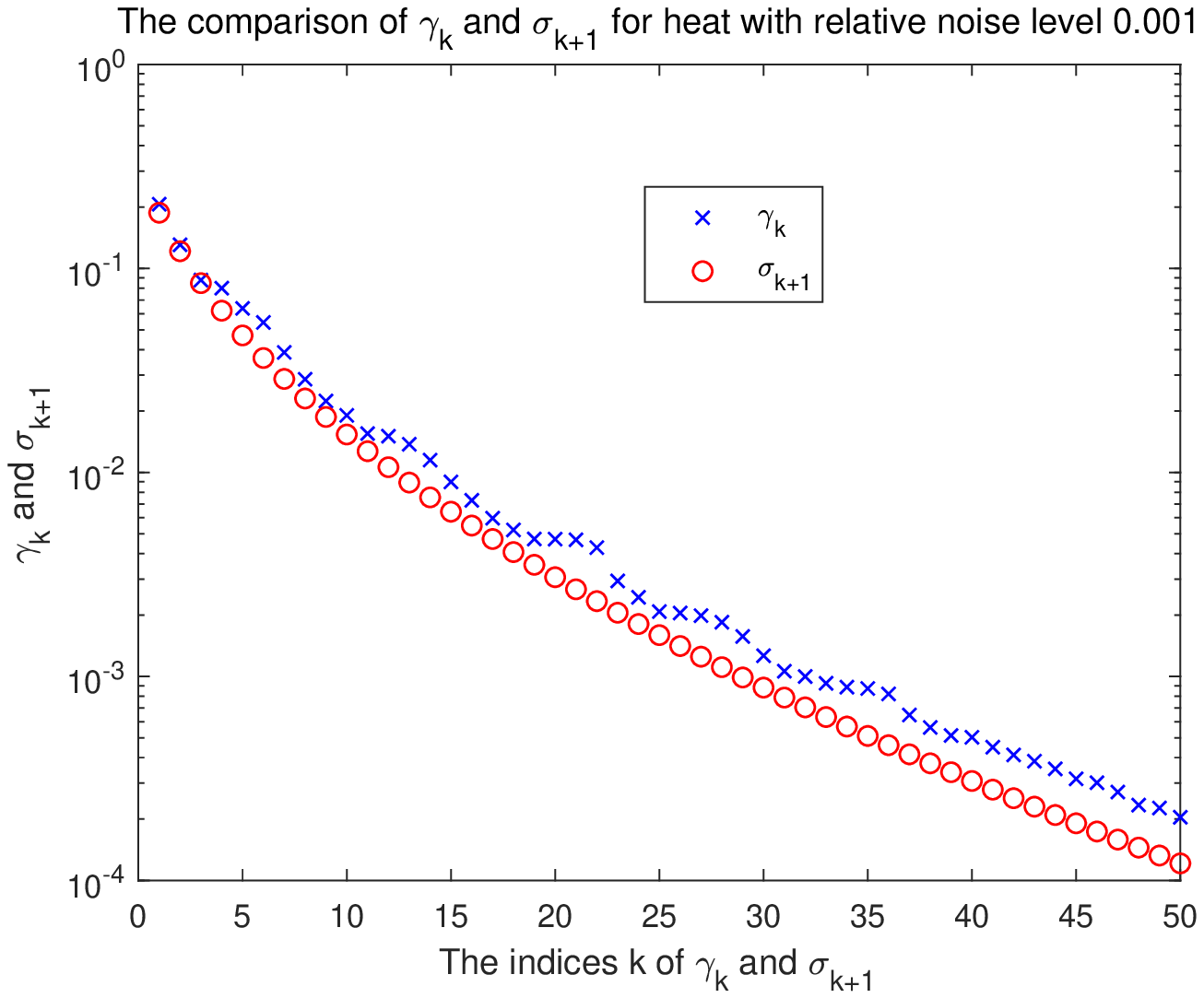}}
  \centerline{(a)}
\end{minipage}
\hfill
\begin{minipage}{0.48\linewidth}
  \centerline{\includegraphics[width=6.0cm,height=4.5cm]{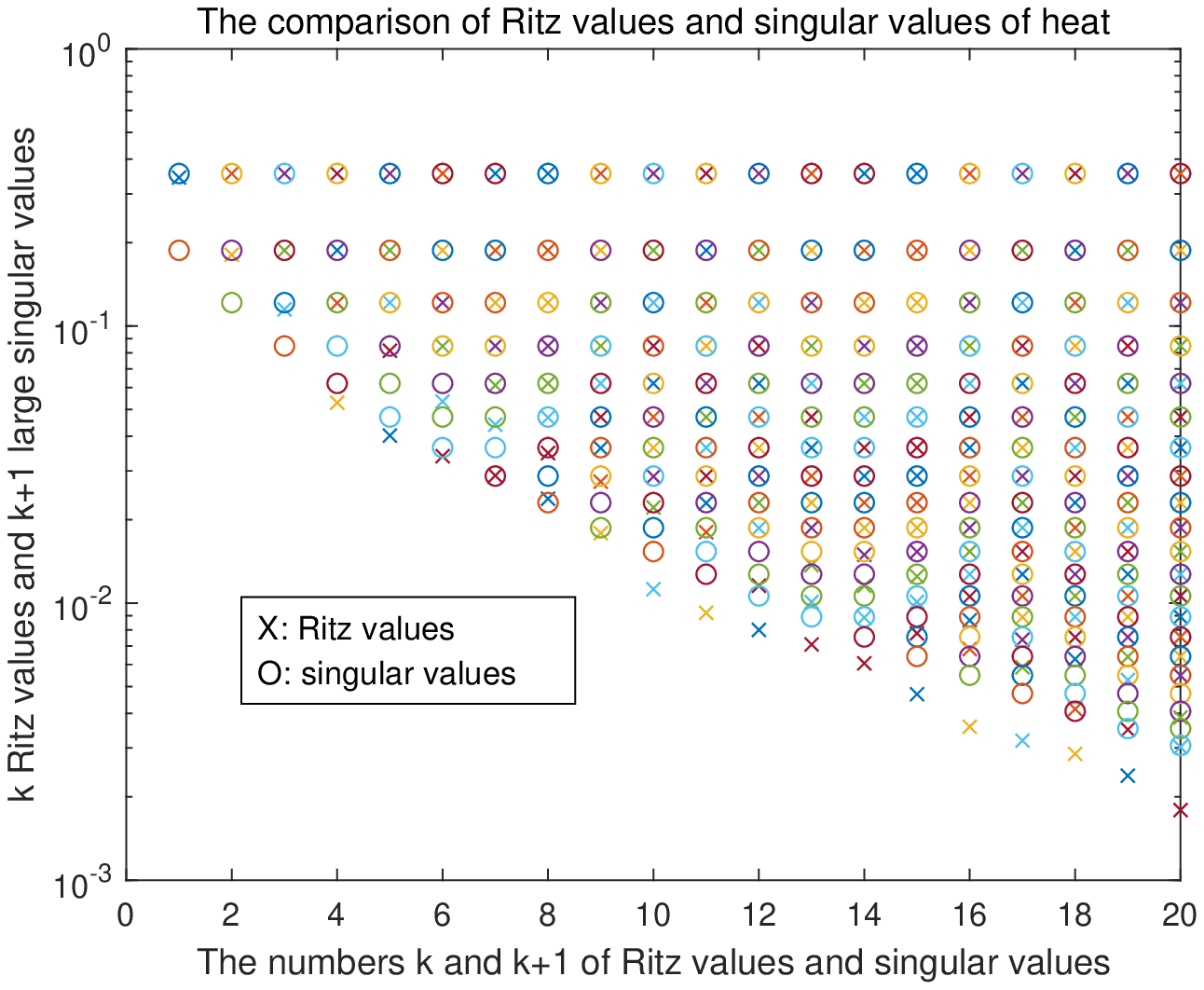}}
  \centerline{(b)}
\end{minipage}
\vfill
\begin{minipage}{0.48\linewidth}
  \centerline{\includegraphics[width=6.0cm,height=4.5cm]{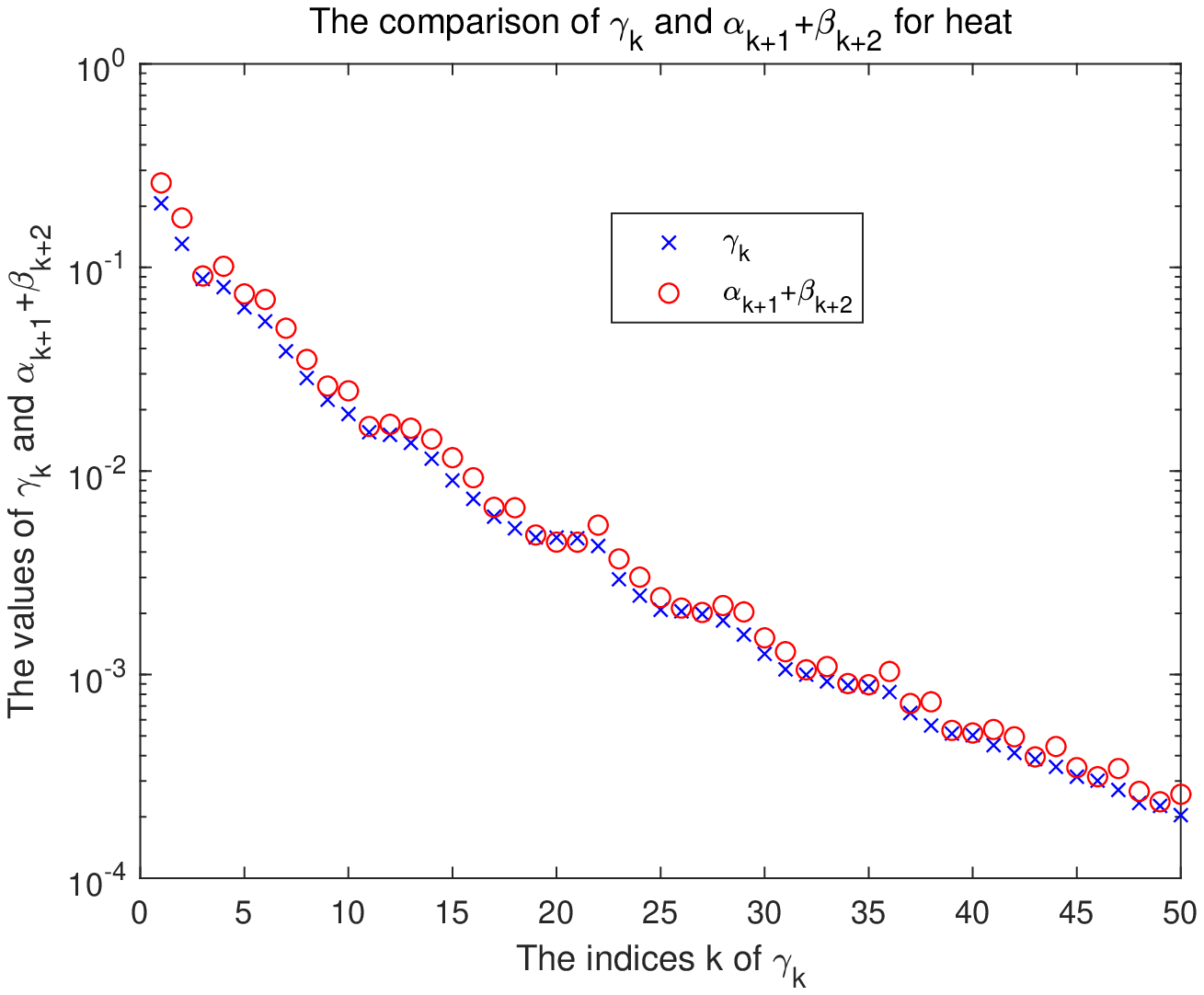}}
  \centerline{(c)}
\end{minipage}
\hfill
\begin{minipage}{0.48\linewidth}
  \centerline{\includegraphics[width=6.0cm,height=4.5cm]{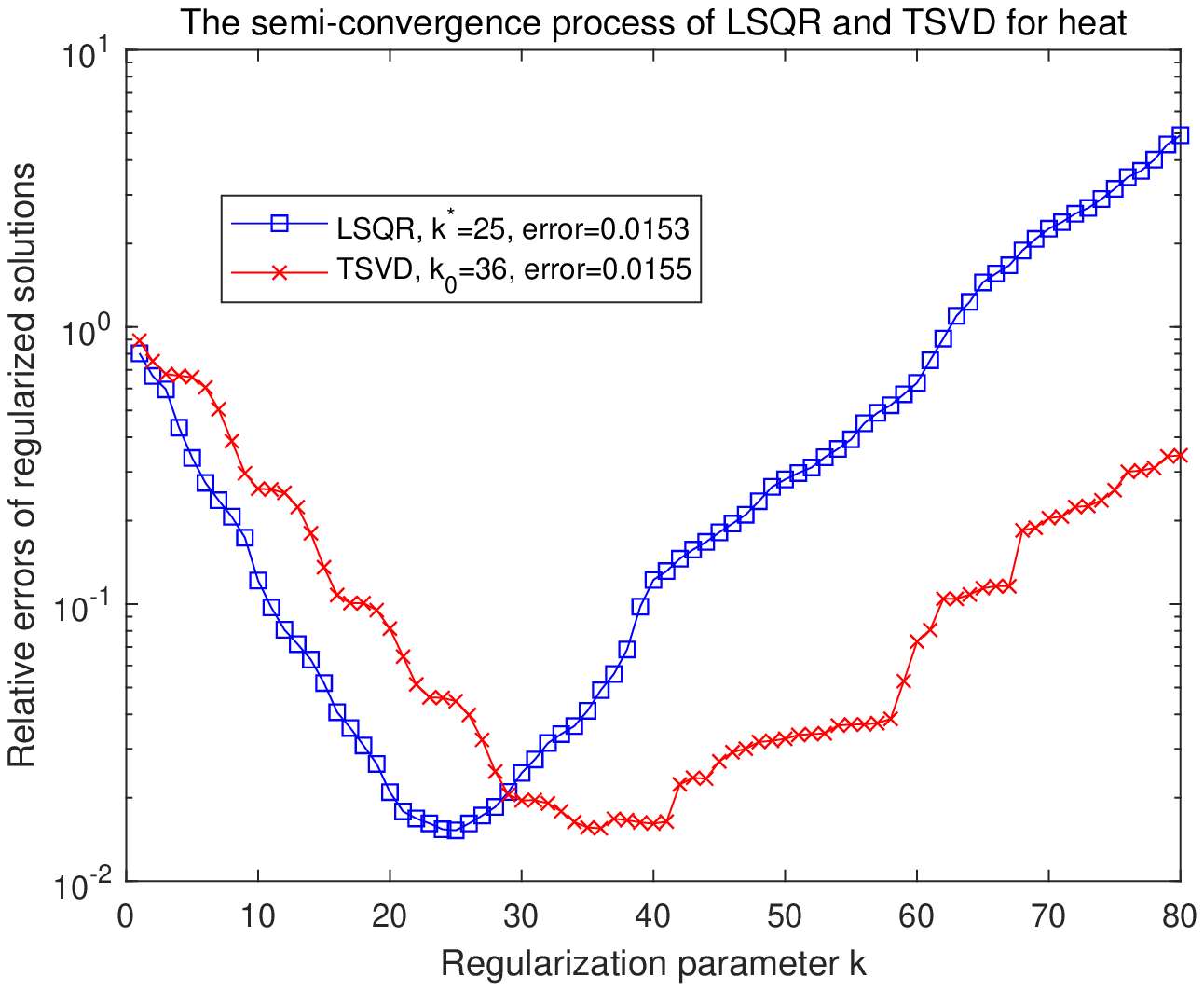}}
  \centerline{(d)}
\end{minipage}
\caption{{\sf heat} of $n=10240$ with the relative noise level $10^{-3}$.} \label{fig3}
\end{figure}

\begin{figure}
\begin{minipage}{0.48\linewidth}
  \centerline{\includegraphics[width=6.0cm,height=4.5cm]{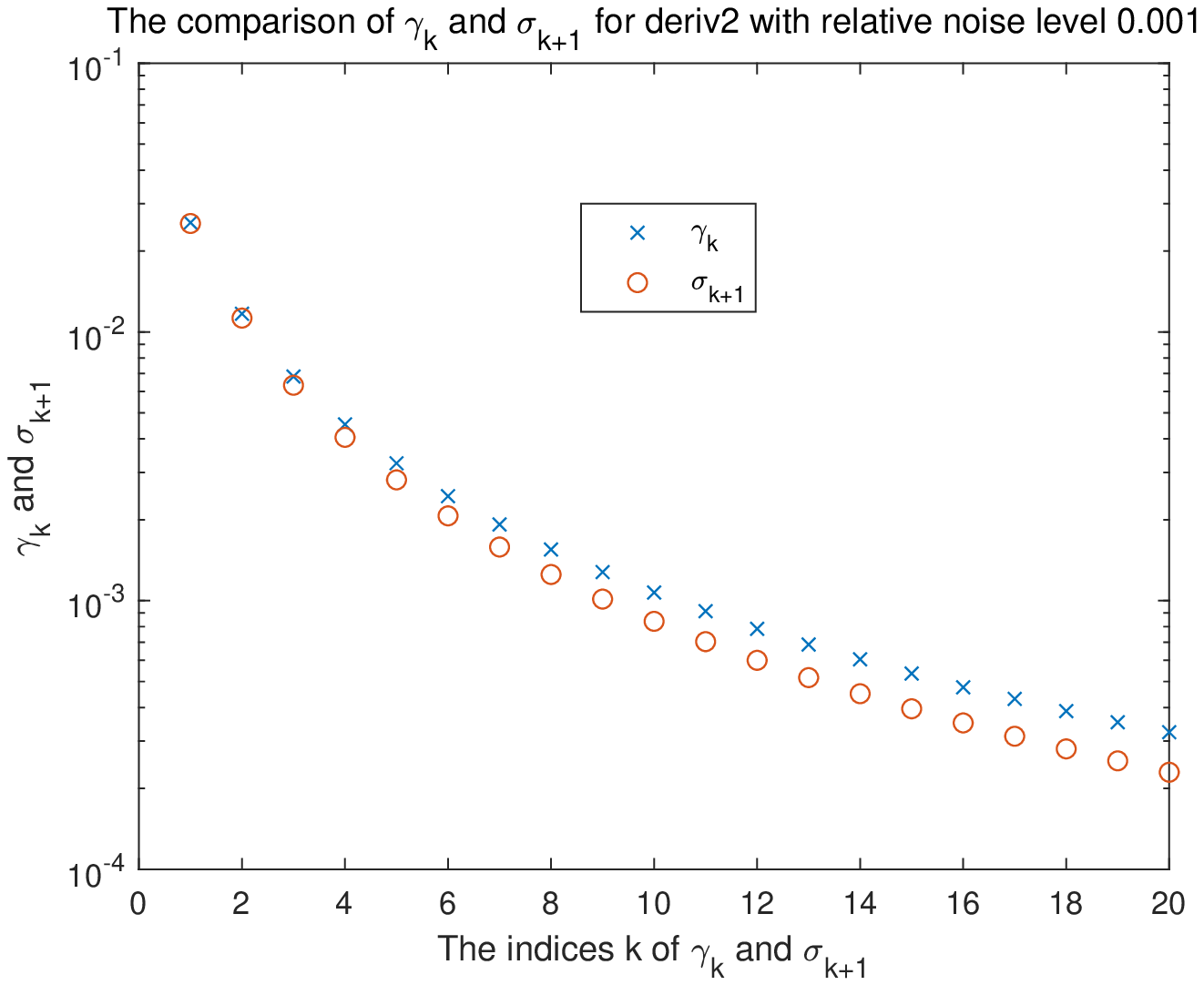}}
  \centerline{(a)}
\end{minipage}
\hfill
\begin{minipage}{0.48\linewidth}
  \centerline{\includegraphics[width=6.0cm,height=4.5cm]{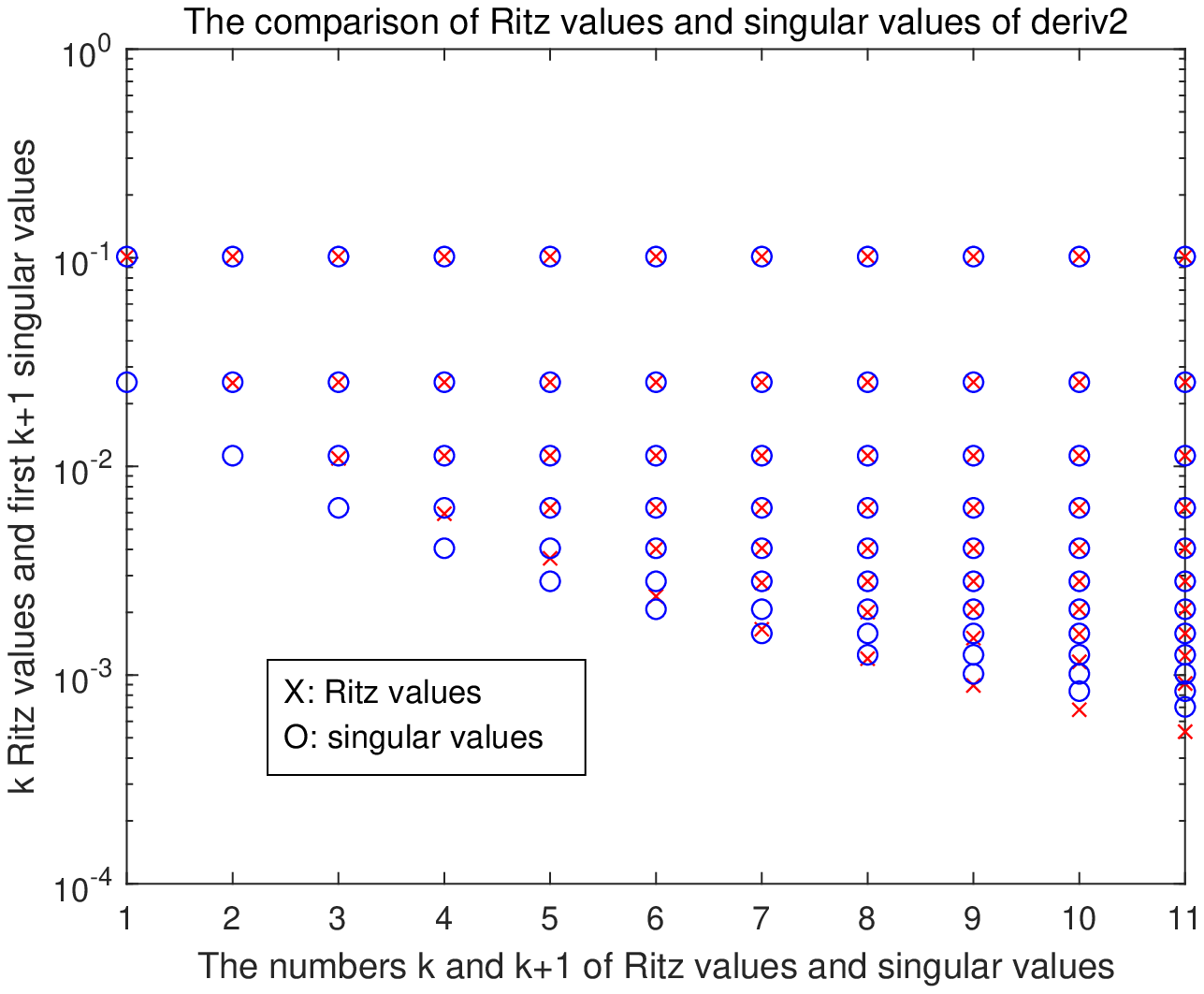}}
  \centerline{(b)}
\end{minipage}
\vfill
\begin{minipage}{0.48\linewidth}
  \centerline{\includegraphics[width=6.0cm,height=4.5cm]{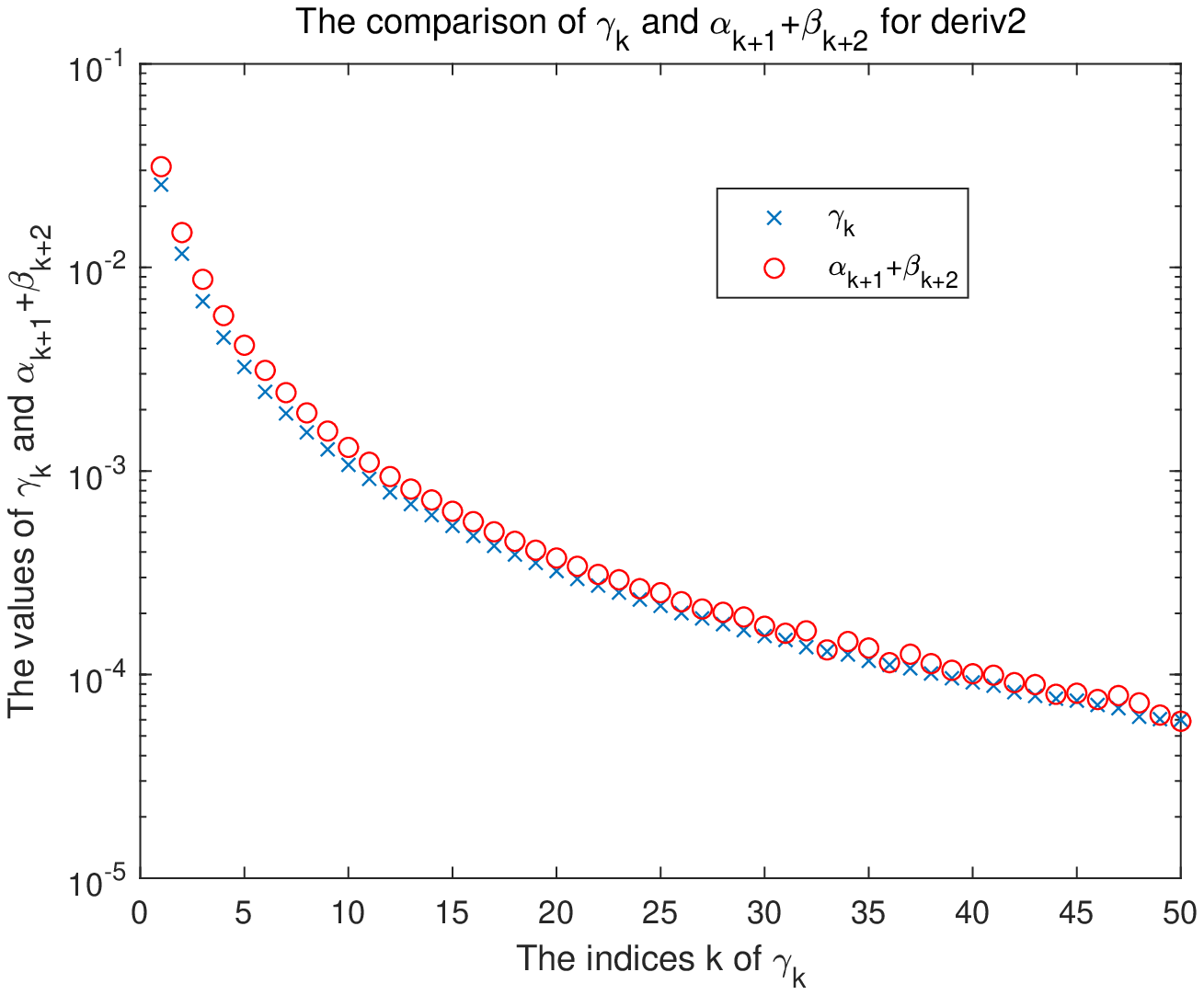}}
  \centerline{(c)}
\end{minipage}
\hfill
\begin{minipage}{0.48\linewidth}
  \centerline{\includegraphics[width=6.0cm,height=4.5cm]{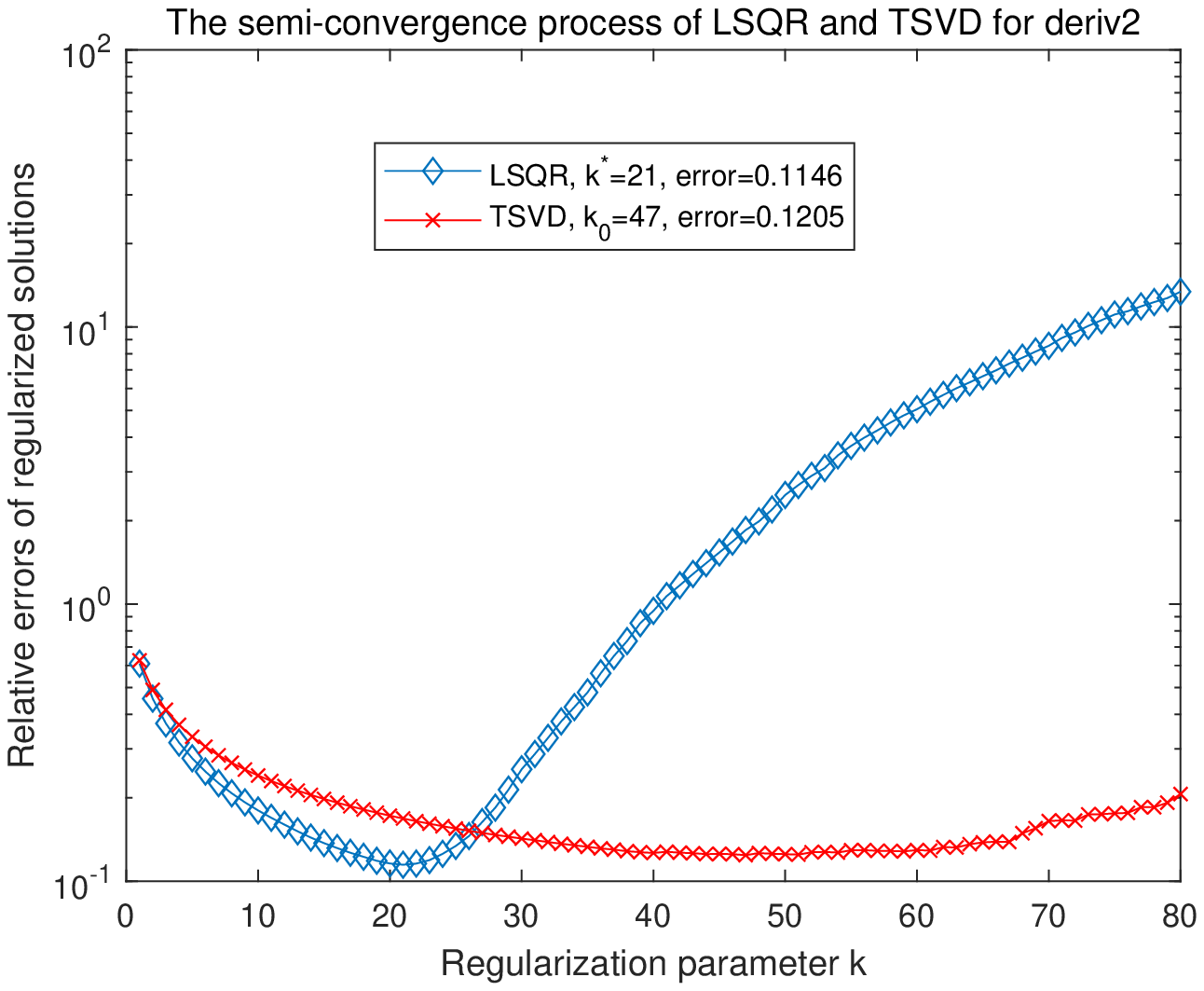}}
  \centerline{(d)}
\end{minipage}
\caption{{\sf deriv2} of $n=10000$ with relative noise level $10^{-3}$.} \label{fig4}
\end{figure}

\begin{figure}
\begin{minipage}{0.48\linewidth}
  \centerline{\includegraphics[width=6.0cm,height=4.5cm]{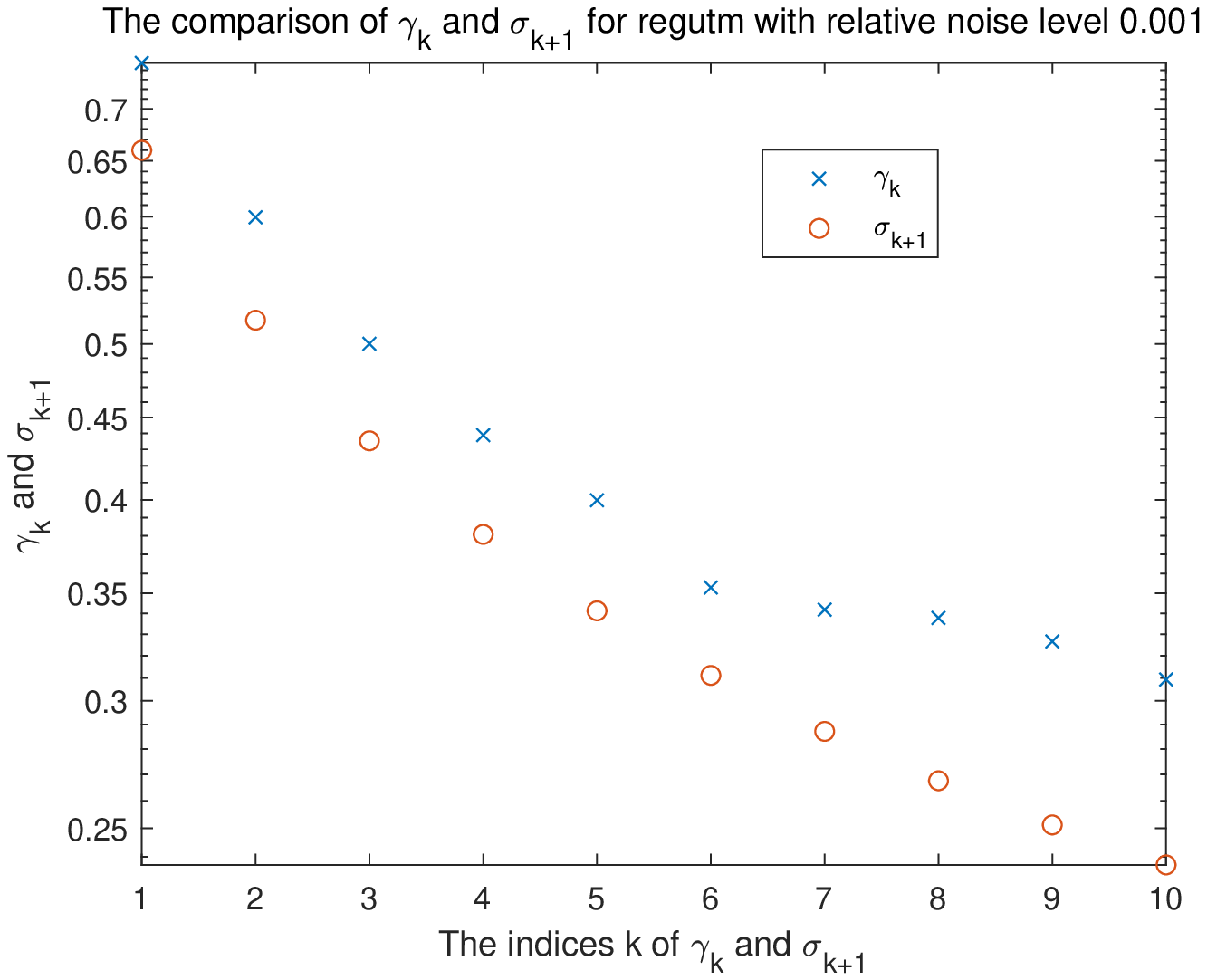}}
  \centerline{(a)}
\end{minipage}
\hfill
\begin{minipage}{0.48\linewidth}
  \centerline{\includegraphics[width=6.0cm,height=4.5cm]{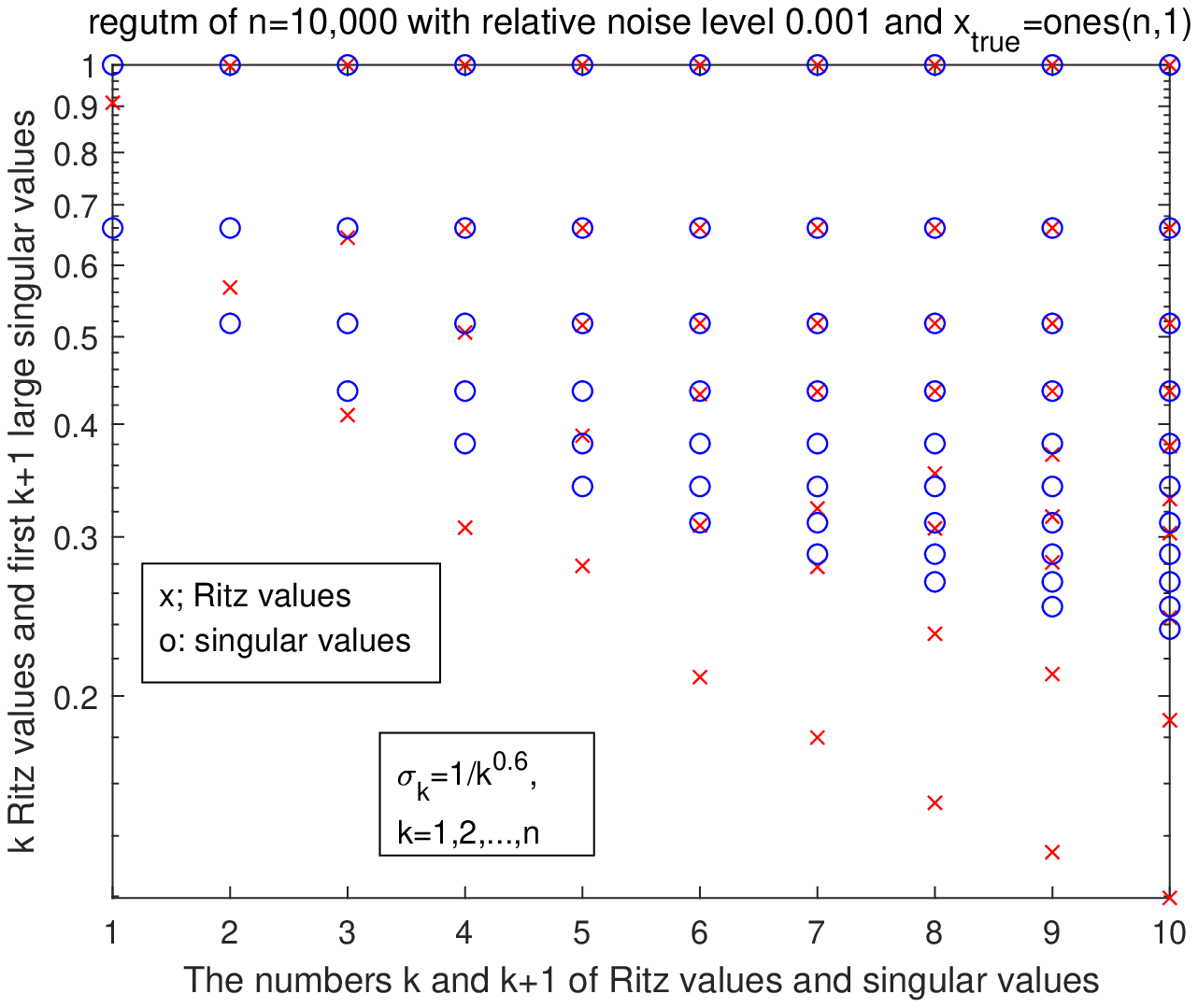}}
  \centerline{(b)}
\end{minipage}
\vfill
\begin{minipage}{0.48\linewidth}
  \centerline{\includegraphics[width=6.0cm,height=4.5cm]{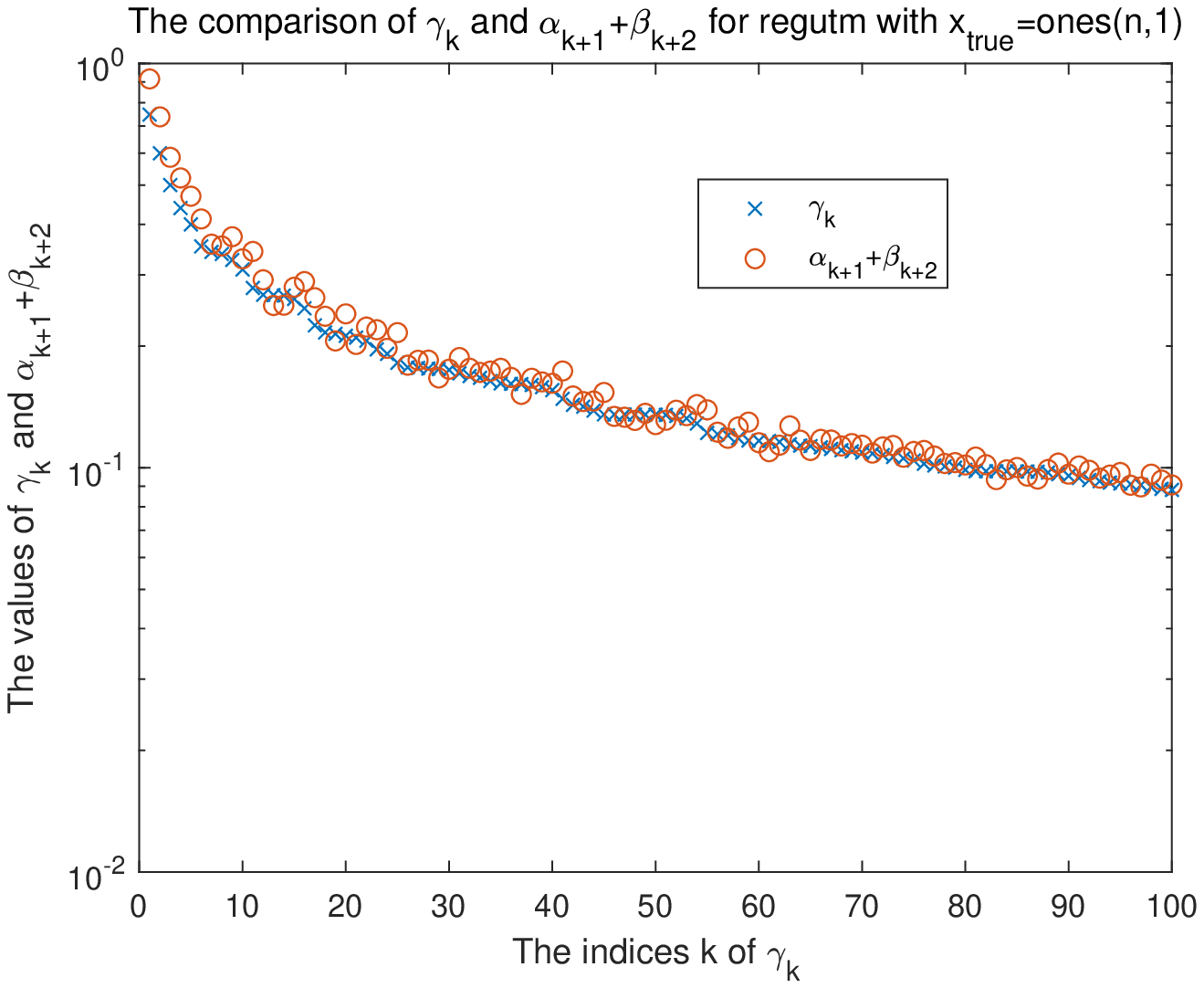}}
  \centerline{(c)}
\end{minipage}
\hfill
\begin{minipage}{0.48\linewidth}
  \centerline{\includegraphics[width=6.0cm,height=4.5cm]{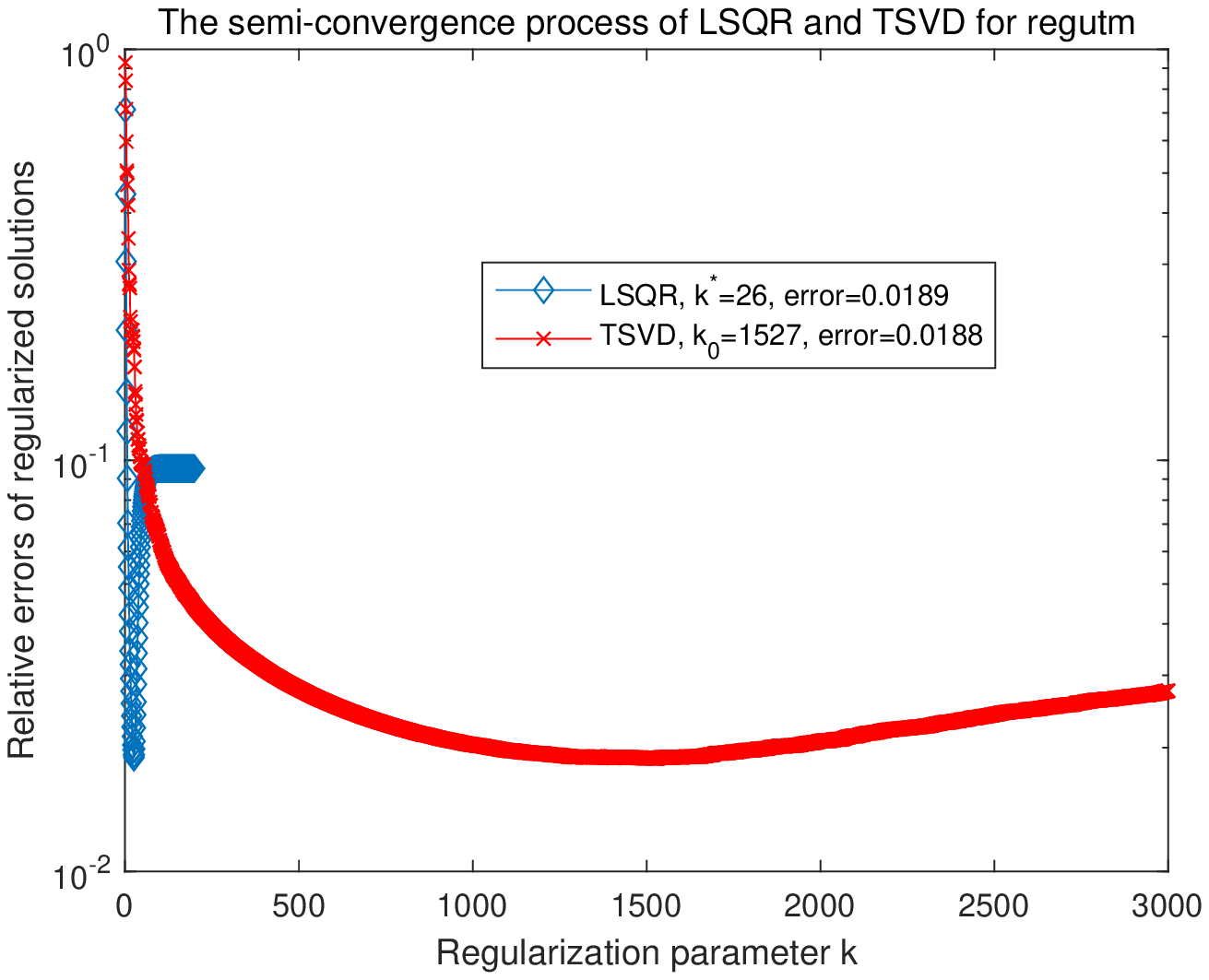}}
  \centerline{(d)}
\end{minipage}
\caption{{\sf regutm} of $m=n=10000$ with the singular values $\sigma_k=\frac{1}{k^{0.6}}$,
$x_{true}=ones(n,1)$ and the relative noise level $10^{-3}$.} \label{fig5}
\end{figure}

\begin{figure}
\begin{minipage}{0.48\linewidth}
  \centerline{\includegraphics[width=6.0cm,height=4.5cm]{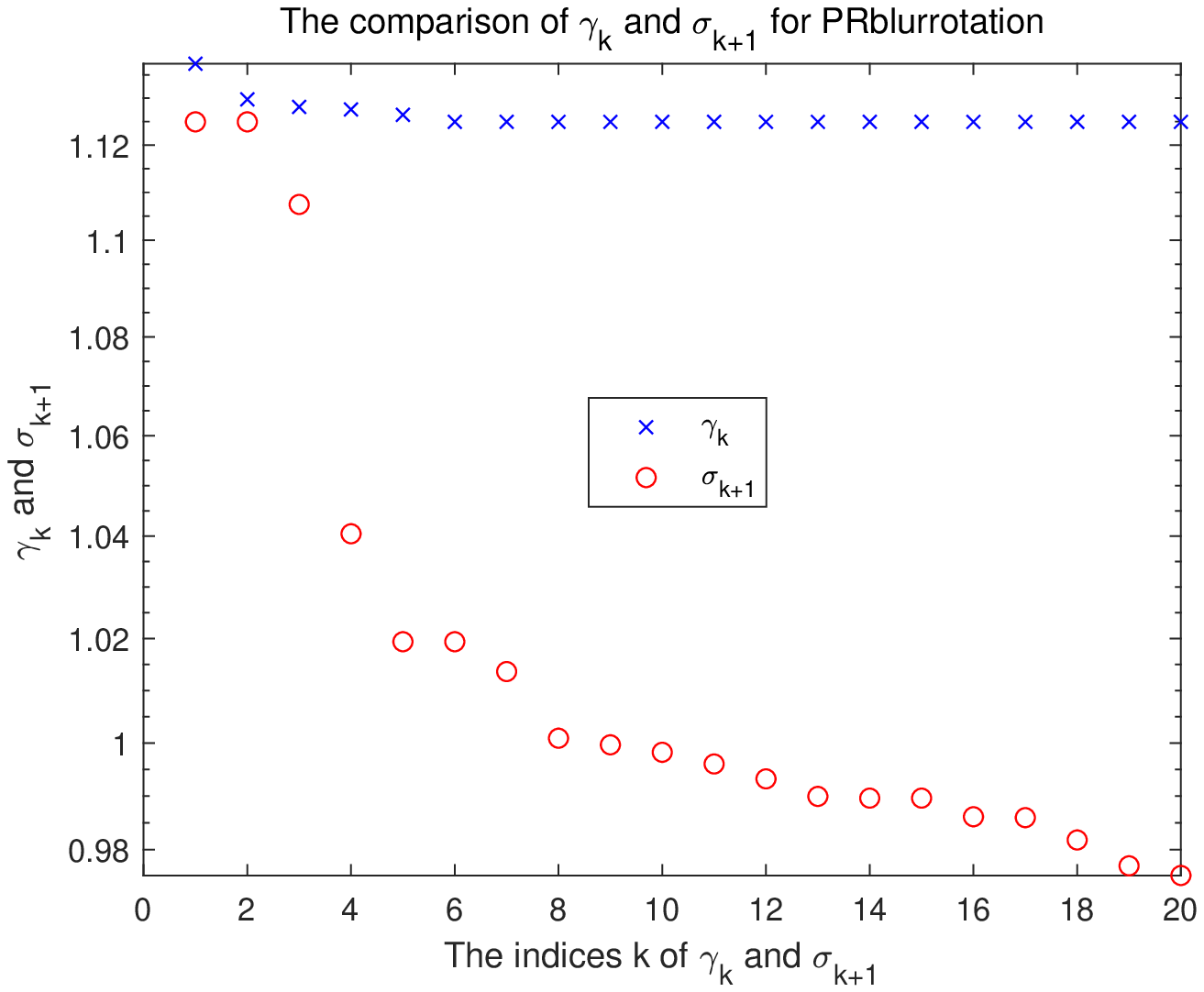}}
  \centerline{(a)}
\end{minipage}
\hfill
\begin{minipage}{0.48\linewidth}
  \centerline{\includegraphics[width=6.0cm,height=4.5cm]{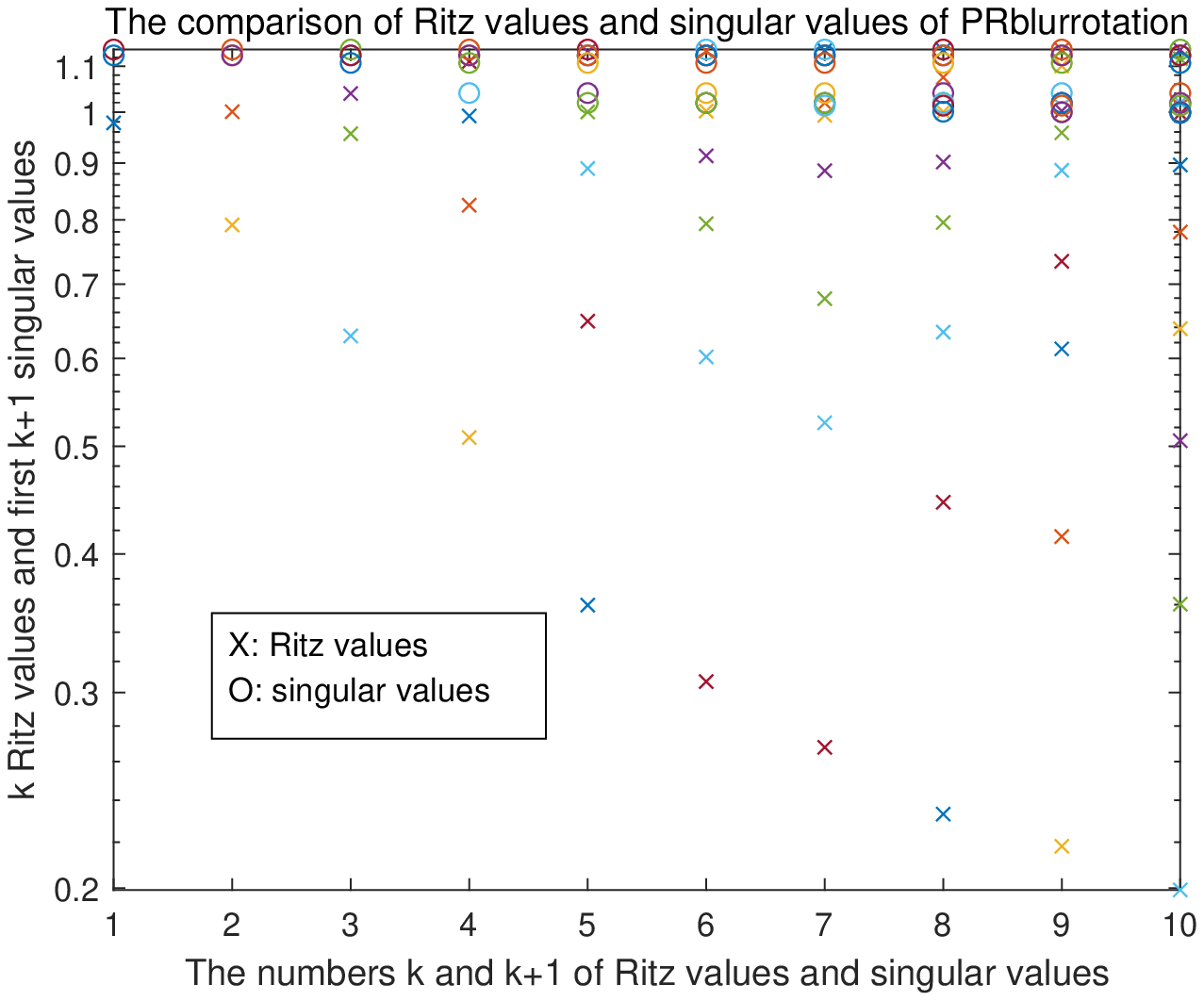}}
  \centerline{(b)}
\end{minipage}
\vfill
\begin{minipage}{0.48\linewidth}
  \centerline{\includegraphics[width=6.0cm,height=4.5cm]{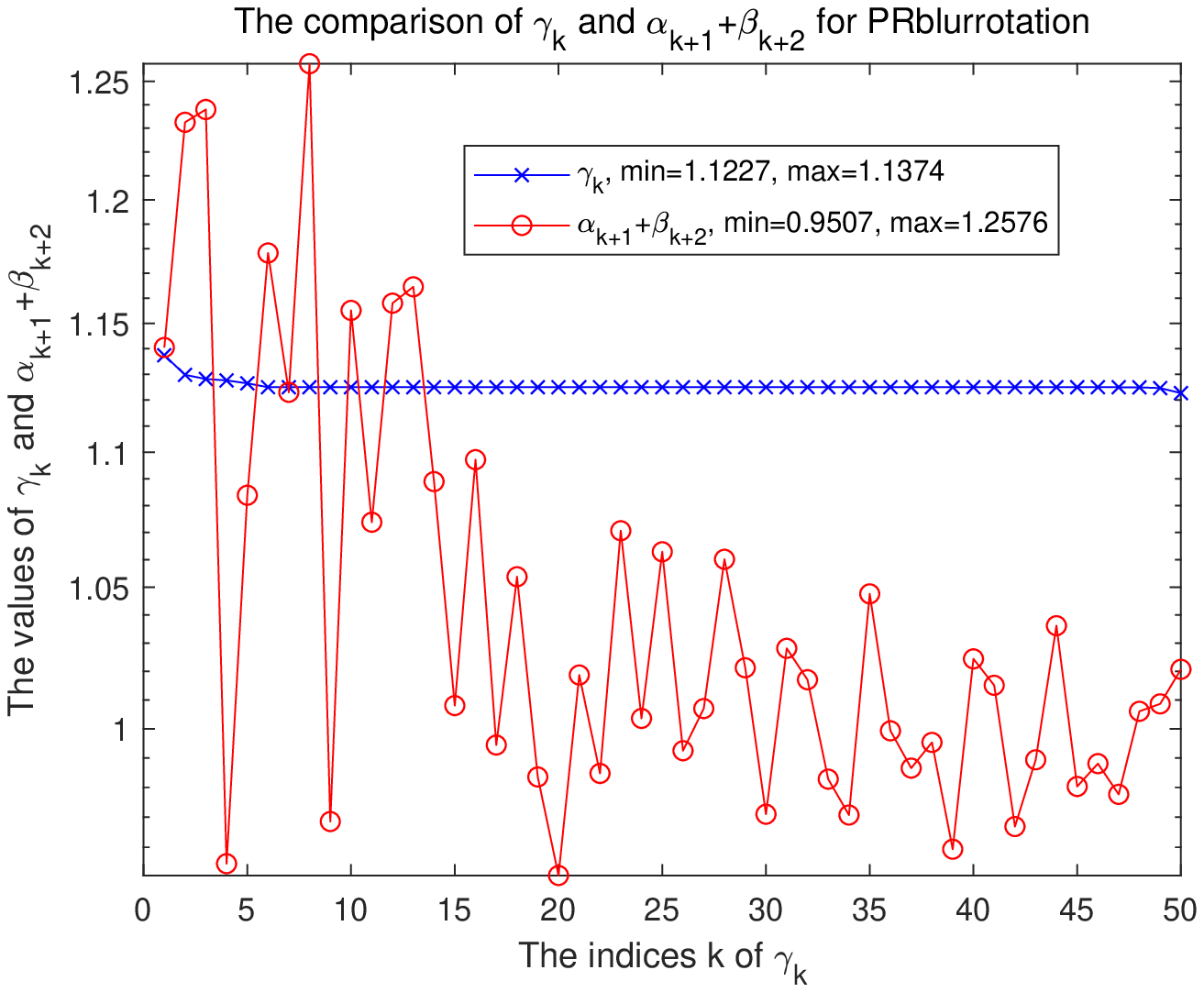}}
  \centerline{(c)}
\end{minipage}
\hfill
\begin{minipage}{0.48\linewidth}
  \centerline{\includegraphics[width=6.0cm,height=4.5cm]{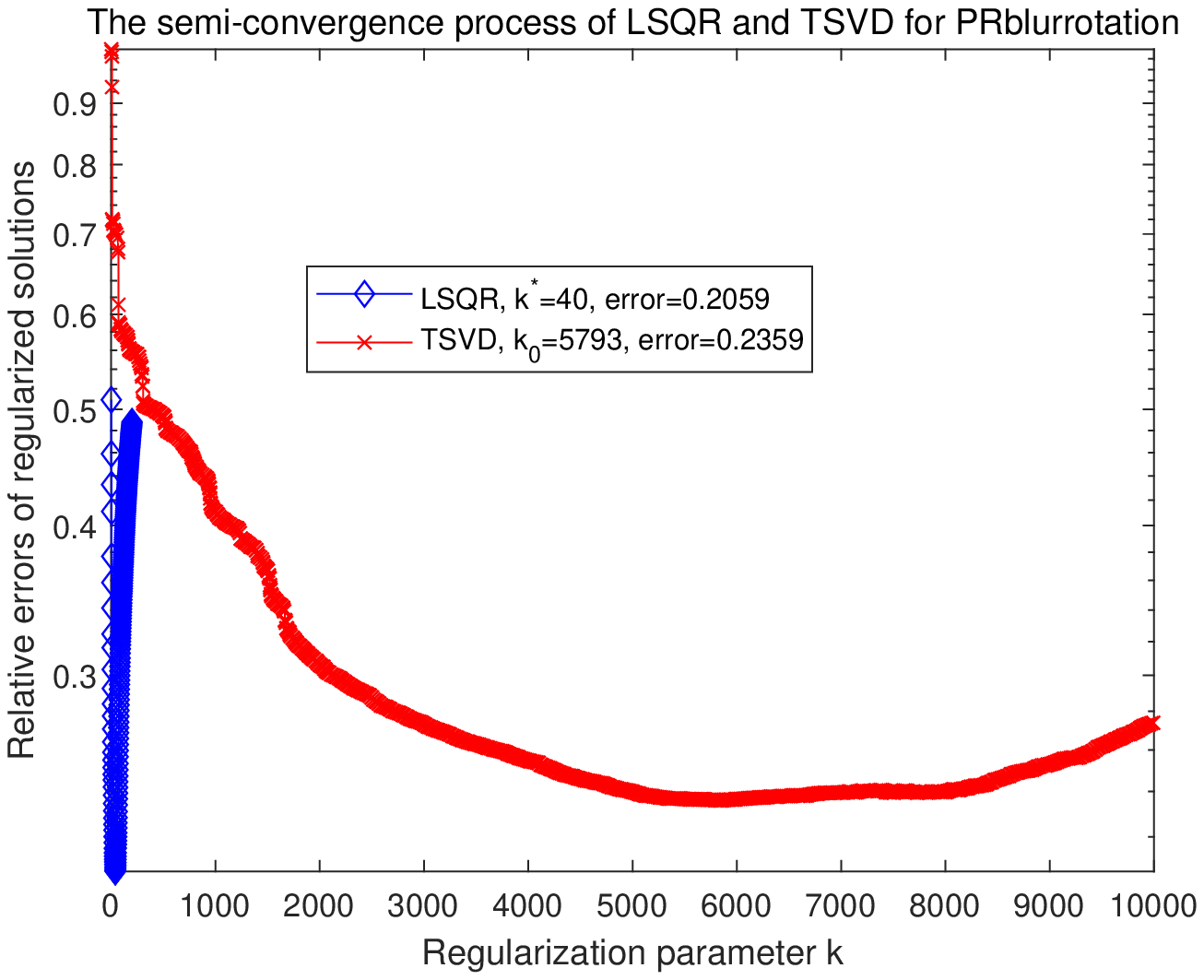}}
  \centerline{(d)}
\end{minipage}
\caption{{\sf PRblurrotation} of $m=n=14400$ with the relative noise level $10^{-2}$.}
\label{figPRblurrotation}
\end{figure}

We first investigate the accuracy $\gamma_k$ of rank $k$ approximations
and the approximation behavior of
$k$ Ritz values $\theta_i^{(k)}$ for the first five test problems
{\sf shaw}, {\sf gravity}, {\sf heat}, {\sf deriv2} and {\sf regutm}.
Given $\varepsilon=10^{-3}$, for each problem and $k=1,2,\ldots$
up to some iteration
bigger than the semi-convergence point $k^*$,
each of Figures~\ref{fig1}--\ref{fig5} depicts the curves of $\gamma_k$ and
$\sigma_{k+1}$, the locations of the $k$ Ritz values $\theta_i^{(k)}$
and the first $k+1$ large singular values $\sigma_i$ of $A$, and the decay
curves of $\gamma_k$ and the sum $\alpha_{k+1}+\beta_{k+2}$. We exhibit them
by (a), (b) and (c) in each figure, respectively.
Figures~\ref{fig1} (d)--\ref{fig5} (d)
draw the semi-convergence processes of LSQR and the TSVD method. We
point out that for different $\varepsilon=10^{-2}$ and $10^{-4}$
we have obtained similar results and observed the same phenomena. So
we omit details on them.
Separately, Figure~\ref{figPRblurrotation} draws the results on
{\sf PRblurrotation} with $\varepsilon=10^{-2}$.


For the severely ill-posed problems {\sf shaw} and {\sf gravity}, their
singular values decay fast, and the $\sigma_1/\sigma_k$
achieve the level of $\frac{1}{\epsilon_{\rm mach}}$ for $k=21\ll n$ and
$k=53\ll n$, respectively. From
Figures~\ref{fig1}--\ref{fig2} (d), we see that the semi-convergence
of TSVD and LSQR occurs at the same steps $k_0=k^*=7$ and the relative
errors
$$
\frac{\|x_{k_0}^{tsvd}-x_{true}\|}{\|x_{true}\|}
\mbox{\ \ and\ \ }
\frac{\|x_{k^*}^{lsqr}-x_{true}\|}{\|x_{true}\|}
$$
of the best regularized solutions $x_{k_0}^{tsvd}$ and
$x_{k^*}^{lsqr}$ are essentially the same for each problem, meaning
LSQR has the full regularization. From Figures~\ref{fig1} (a)--\ref{fig2} (a),
we observe that the $\gamma_k$ are very close to
$\sigma_{k+1}$ and are almost indispensable for $k=1,2,\ldots,k^*$.
This indicates that the $P_{k+1}B_kQ_k^T$ are near
best rank $k$ approximations to $A$ at least until the semi-convergence of
LSQR, confirming Theorems~\ref{main1} and \ref{nearapprox}, which states
that Lanczos bidiagonalization
can generate near best rank $k$ approximations for suitable $\rho>1$
until the semi-convergence of LSQR. Figures~\ref{fig1} (b)--\ref{fig2} (b)
tell us that, for {\sf shaw}, the $k$ Ritz values $\theta_i^{(k)}$ approximate
the first $k$ large singular values of $A$ in natural order, or
interlace the first $k+1$ large ones, at least until the
semi-convergence of LSQR. This confirms Theorem~\ref{ritzvalue}.
Furthermore, we can see from Figure~\ref{fig1} (a)--(b)
that the near best approximations and the approximations of $\theta_i^{(k)}$
in this order are valid at least until $k=16>k^*$  and $k=17>k^*$, respectively,
but for {\sf gravity} both of them are valid only until  $k=13>k^*$, smaller than
those for {\sf shaw}, as is seen from Figure~\ref{fig2} (a)-(b).
This is not surprising because the singular values of {\sf gravity}
do not decay as fast as those of {\sf shaw} though both of them are
severely ill-posed.

For the moderately ill-posed {\sf heat}, we find
$\sigma_1/\sigma_{900}=8.3\times 10^{15}$; for {\sf deriv2}, we
find with $\sigma_1/\sigma_n=1.2\times 10^8$. We can observe from
Figure~\ref{fig3} (d) and Figure~\ref{fig4} (d)
that the semi-convergence of LSQR occurs earlier than that of TSVD,
i.e., $k^*<k_0$.
Theorem~\ref{semicon} states that
in this case the $k$ Ritz values $\theta_i^{(k)}$ must not approximate
the large singular values of $A$ in natural order, or
they do not interlace the first $k+1$ large ones for some $k\leq k^*$.
This is indeed true, and such phenomena occur from $k=4$ and $8$ onwards
for {\sf heat} and {\sf deriv2}, respectively,
as is seen clearly from Figure~\ref{fig3} (b) and Figure~\ref{fig4} (b).
These are in accordance with Theorem~\ref{ritzvalue},
where, for a given moderately ill-posed problem,
the sufficient condition \eqref{condm} for the approximations in this order
fails to meet as $k$ increases up to some point since \eqref{const2} and the
previous results show  that the left-hand side of \eqref{condm}
monotonically increases, while its right-hand side monotonically
decreases with respect to $k$. In the meantime,
Figure~\ref{fig3} (a) and Figure~\ref{fig4} (a) illustrate that Lanczos
bidigonalization
generates near best rank $k$ approximations to $A$ for $k$ no more than 4
and 9 for  {\sf heat} and {\sf deriv2}, respectively,
after which we cannot obtain near best rank $k$
approximations to $A$ any longer. These confirm
Theorems~\ref{main1}--\ref{nearapprox}, which show
that the rank $k$ approximation
becomes poorer as $k$ increases and it is no more a near best one when
$k$ increases up to some point smaller than $k^*$
since the sufficient condition \eqref{condition1} fails to fulfill.

We now check the random mildly ill-posed {\sf regutm}.
Figure~\ref{fig5} (d) shows that the semi-convergence of LSQR
occurs much more early than that of
TSVD, that is, $k^*=26\ll k_0=1527$. Since the singular values $\sigma_k$ decay
to zero for $n$ sufficiently large but substantially more slowly than the
singular values of {\sf shaw}, {\sf gravity}, {\sf heat} and {\sf deriv2},
Theorem~\ref{main1} and the previous analysis
tell us that $\gamma_k$ deviates
from $\sigma_{k+1}$ quickly as $k$ increases. Actually,
the condition \eqref{condition1} may hold only for $k$ very small.
Figure~\ref{fig5} (a)
justifies  our theory, and from it we see that $P_{k+1}B_kQ_k^T$ is
not a near best rank $k$ approximation for $k=3$ onwards by noticing that
$\gamma_3$ is closer to $\sigma_3$ other than $\sigma_4$. Simultaneously,
the Ritz values $\theta_i^{(k)}$ fail to interlace the first $k+1$
large singular values of $A$ for $k=3$ onwards, as indicated by
Figure~\ref{fig5} (b). This demonstrates Theorem~\ref{nearapprox} and
Theorem~\ref{ritzvalue}, where the sufficient condition \eqref{condition1}
for near best rank $k$ approximations
may be satisfied only for $k$ very small and
the sufficient condition \eqref{condm} for the approximations in this order
is not satisfied, which implies that the $k$ Ritz values hardly
approximate the large singular values in natural order even for $k$ very small.

Finally, we look into the results on {\sf PRblurrotation}. It is seen from
Figure~\ref{figPRblurrotation} (a) that $P_{k+1}B_kQ_k^T$ is not a
near best rank $k$ approximation from the first iteration $k=1$ upwards
and $\gamma_k$ is considerably bigger than $\sigma_{k+1}$ for $k=1,2,\ldots,20$.
Figure~\ref{figPRblurrotation} (b) indicates that the $k$ Ritz values
$\theta_i^{(k)}$ never interlace the first $k+1$ large $\sigma_i$ for given
$k=1,2,\ldots,10$, and $\theta_k^{(k)}<\sigma_{k+1}$ substantially.
Figure~\ref{figPRblurrotation} (d) shows that LSQR computes its best regularized
solution $x_{k^*}^{lsqr}$ at $k^*=40$, much more early than the TSVD
method, which obtains its best solution
$x_{k_0}^{tsvd}$ at $k_0=5793$.

\begin{figure}
\begin{minipage}{0.48\linewidth}
  \centerline{\includegraphics[width=6.0cm,height=4.5cm]{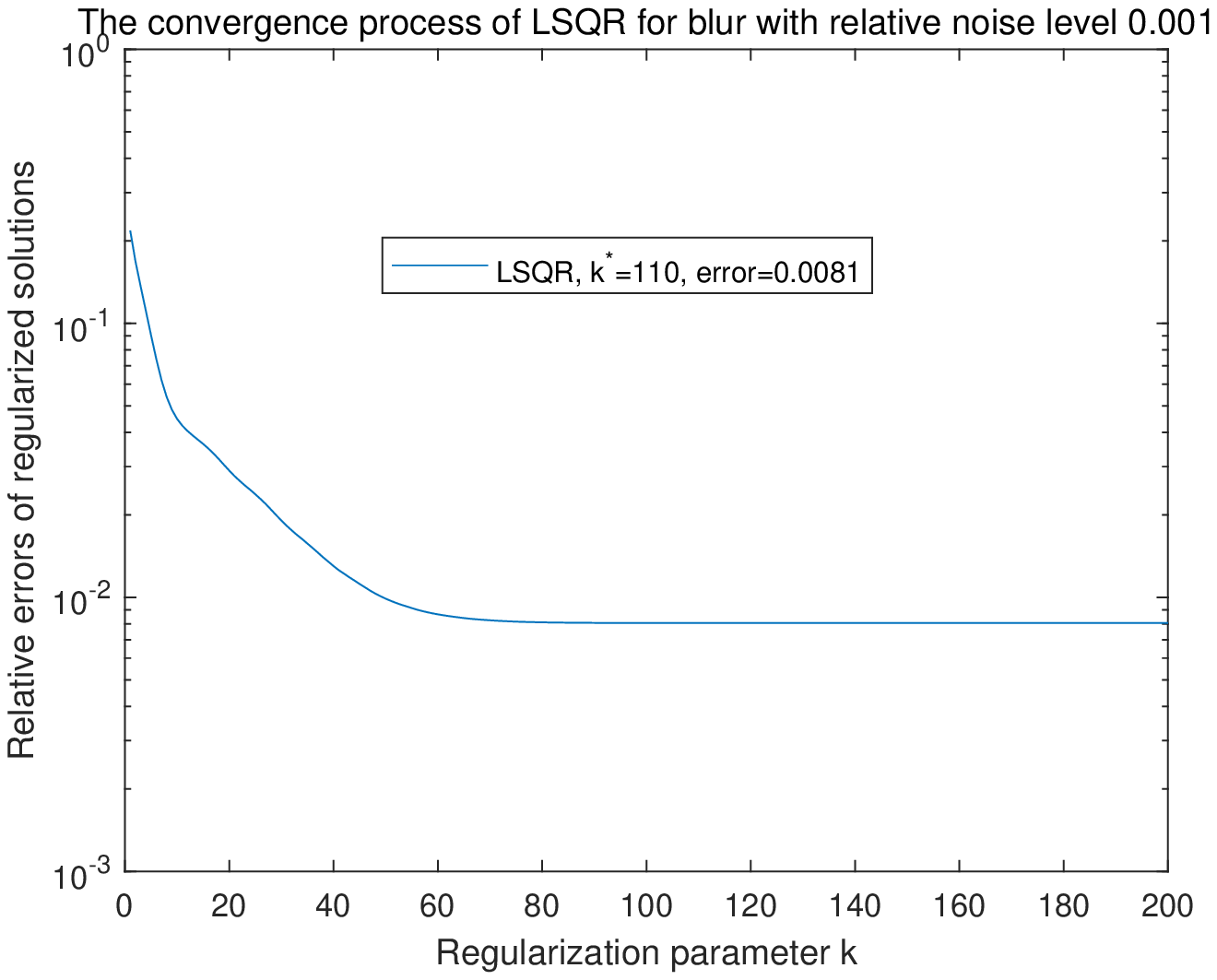}}
  \centerline{(a)}
\end{minipage}
\hfill
\begin{minipage}{0.48\linewidth}
  \centerline{\includegraphics[width=6.0cm,height=4.5cm]{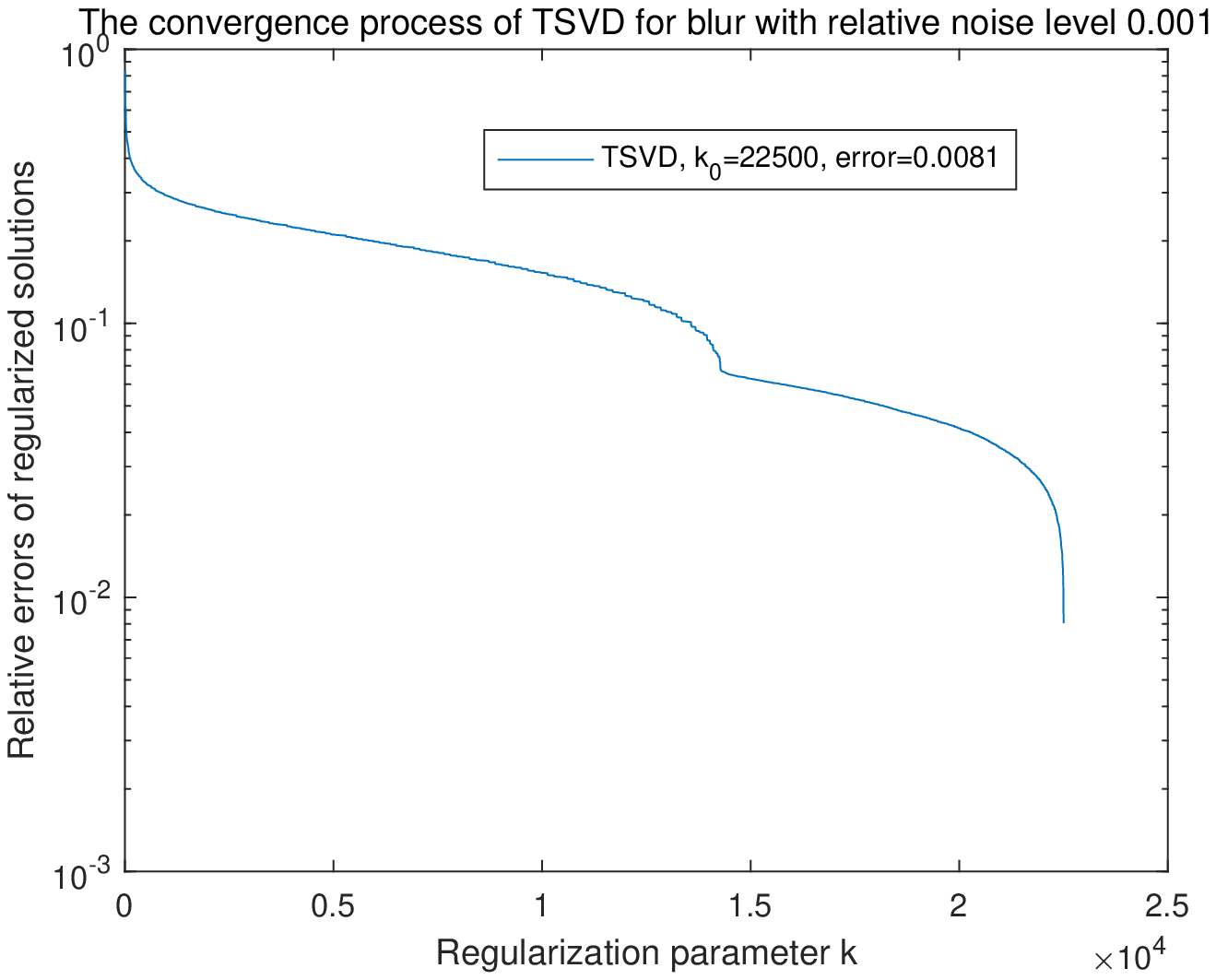}}
  \centerline{(b)}
\end{minipage}
\vfill
\begin{minipage}{0.48\linewidth}
  \centerline{\includegraphics[width=6.0cm,height=4.5cm]{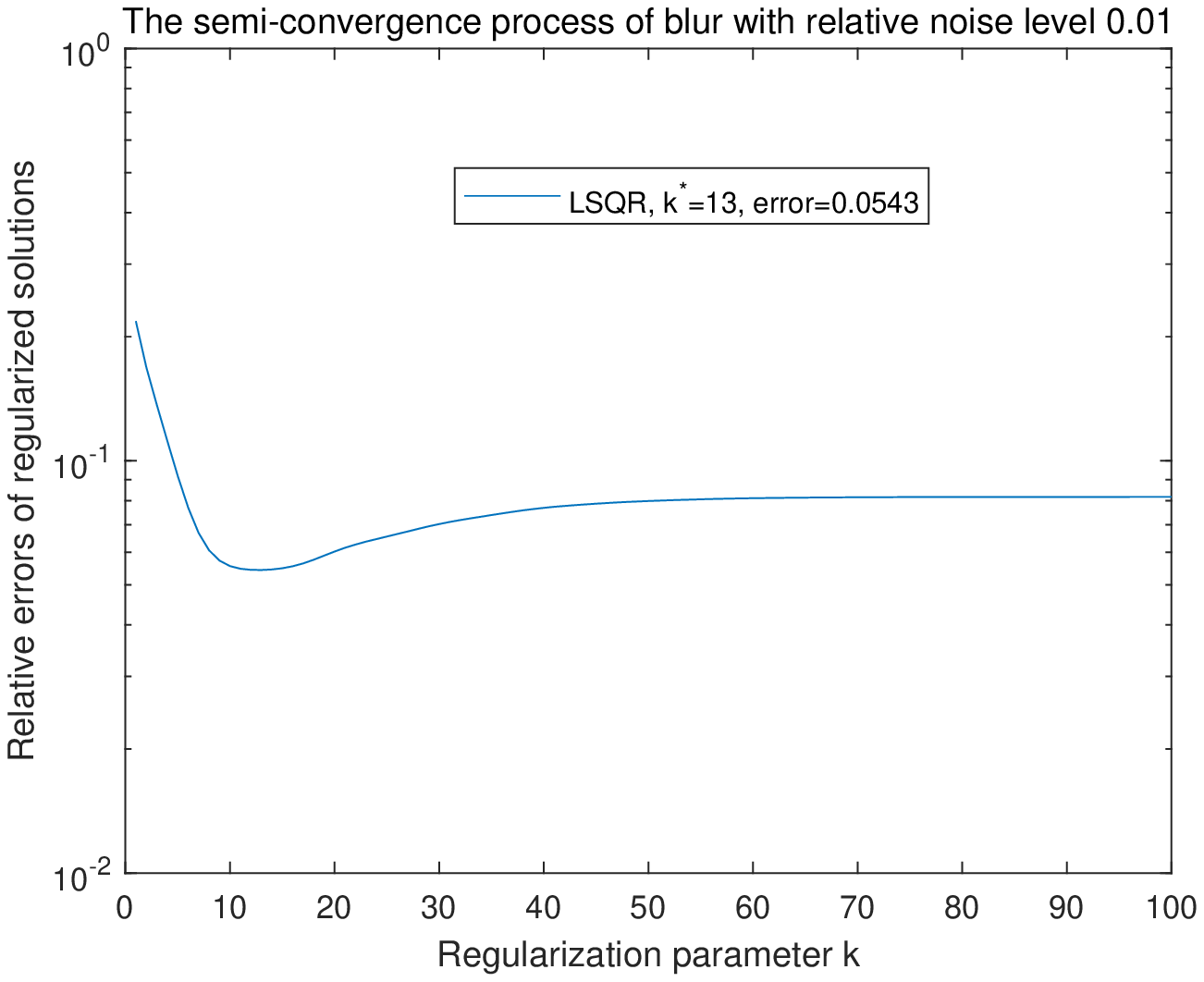}}
  \centerline{(c)}
\end{minipage}
\hfill
\begin{minipage}{0.48\linewidth}
  \centerline{\includegraphics[width=6.0cm,height=4.5cm]{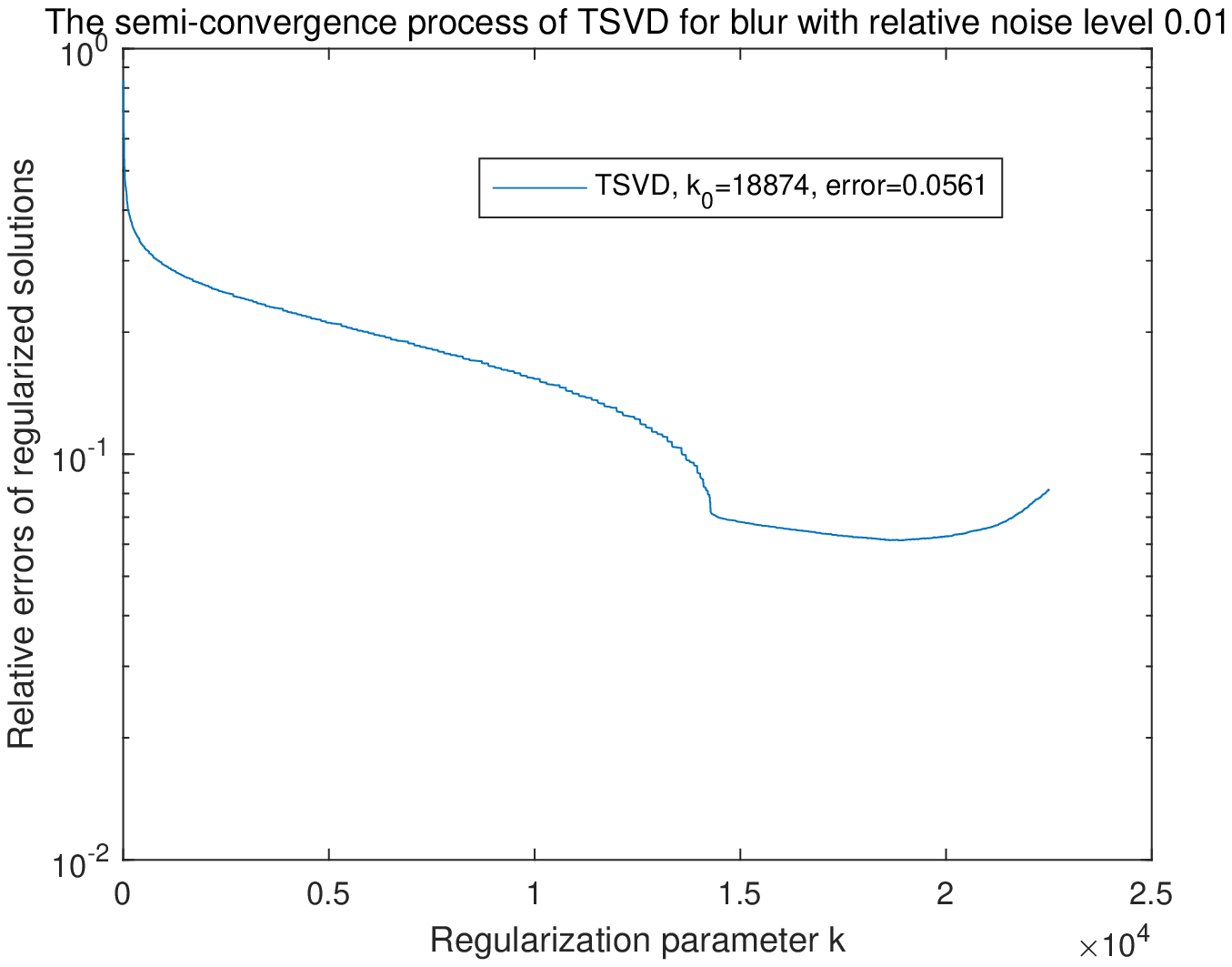}}
  \centerline{(d)}
\end{minipage}
\caption{{\sf blur} of $n=22500$ with $\frac{\sigma_1}{\sigma_n}=31.5$ and
$\varepsilon=10^{-3}$ and $10^{-2}$.} \label{fig6}
\end{figure}

\begin{figure}
\begin{minipage}{0.48\linewidth}
  \centerline{\includegraphics[width=6.0cm,height=4.5cm]{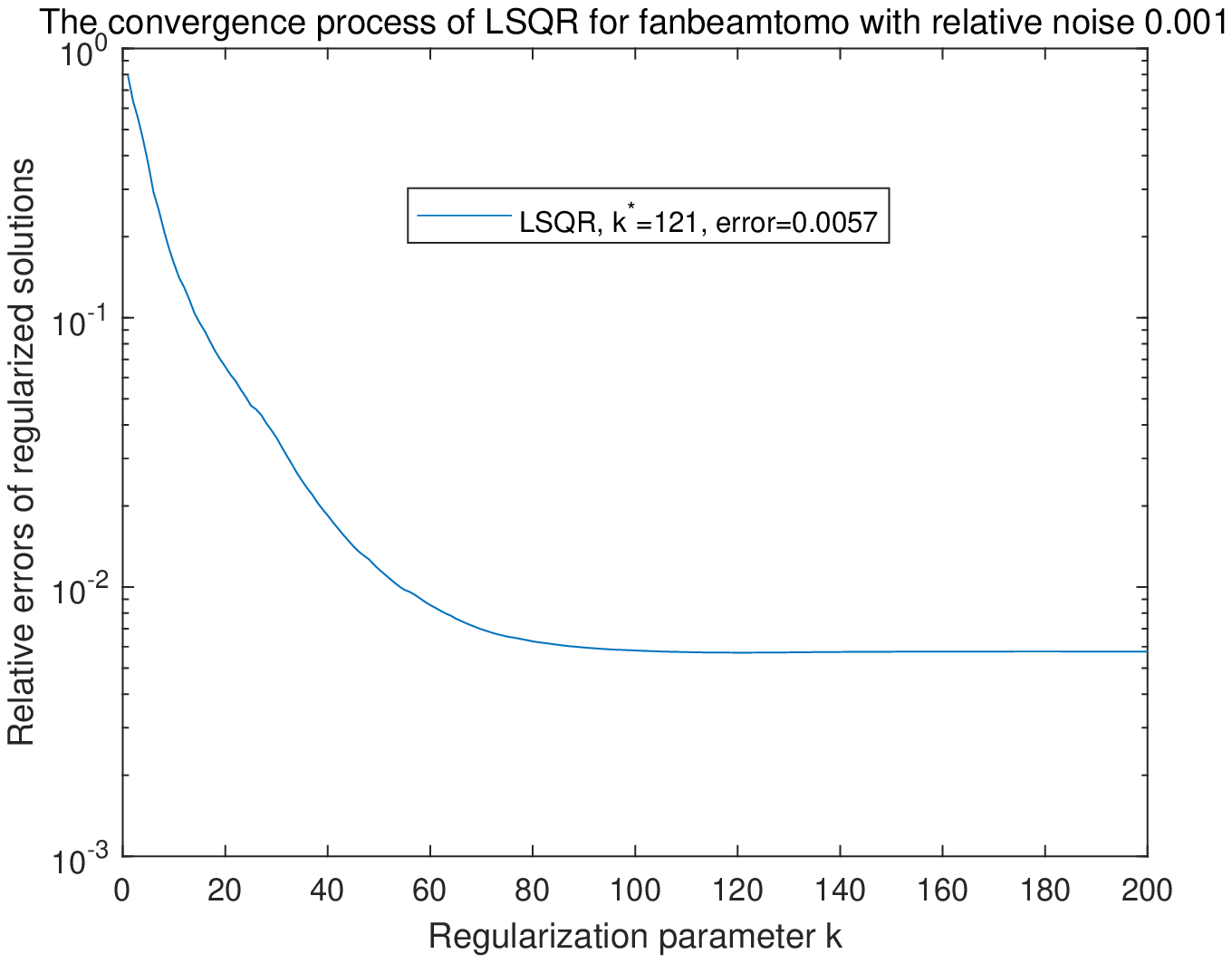}}
  \centerline{(a)}
\end{minipage}
\hfill
\begin{minipage}{0.48\linewidth}
  \centerline{\includegraphics[width=6.0cm,height=4.5cm]{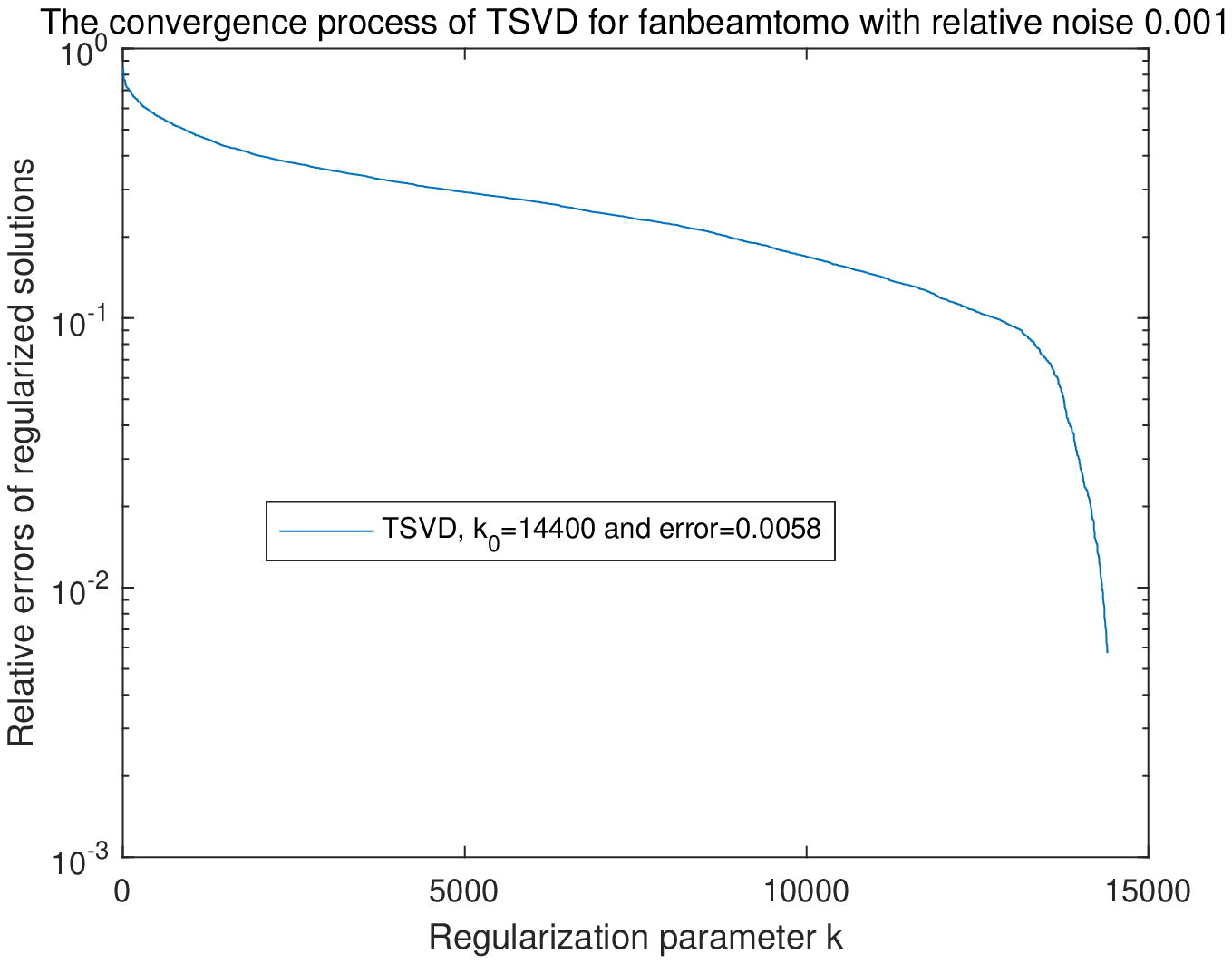}}
  \centerline{(b)}
\end{minipage}
\vfill
\begin{minipage}{0.48\linewidth}
  \centerline{\includegraphics[width=6.0cm,height=4.5cm]{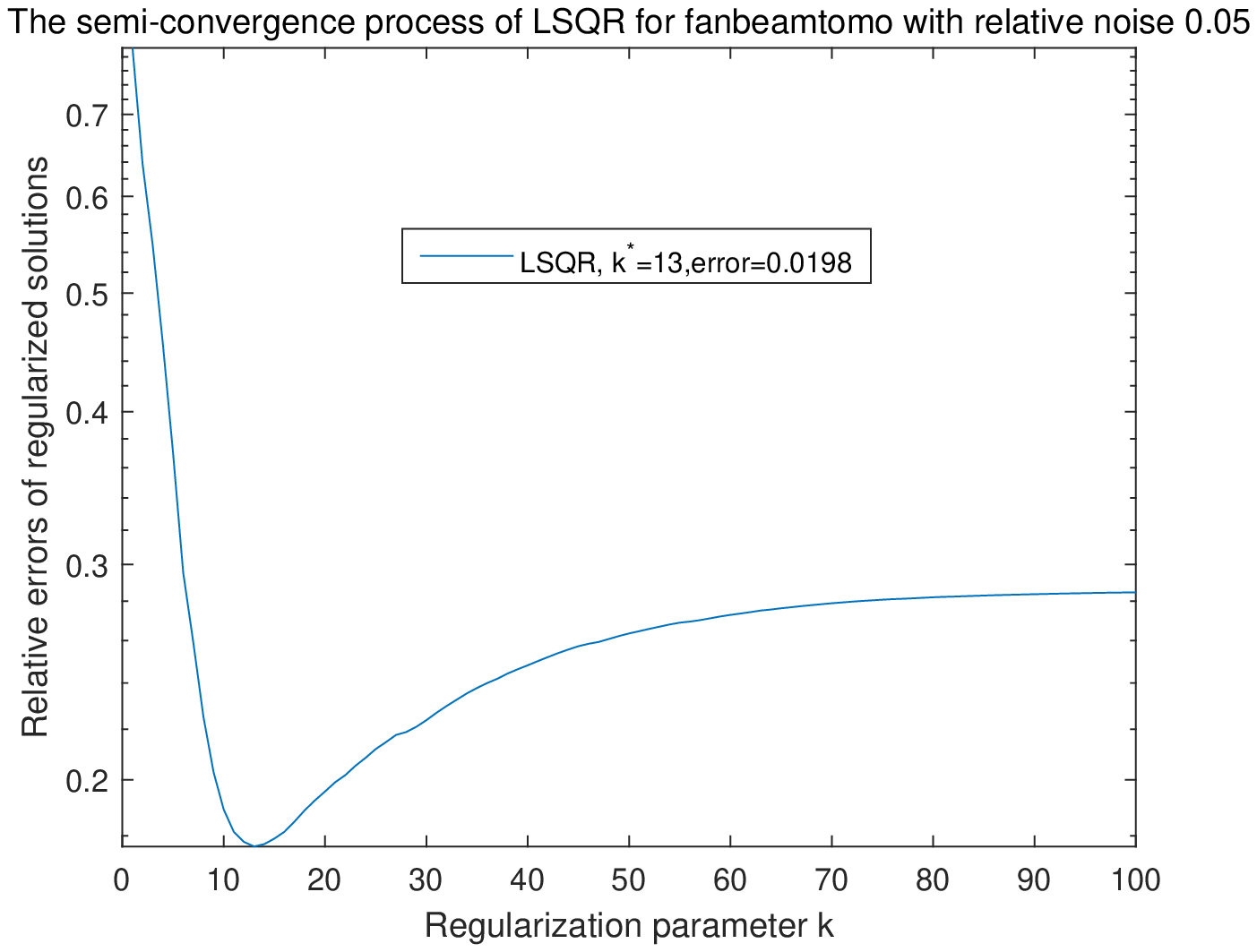}}
  \centerline{(c)}
\end{minipage}
\hfill
\begin{minipage}{0.48\linewidth}
  \centerline{\includegraphics[width=6.0cm,height=4.5cm]{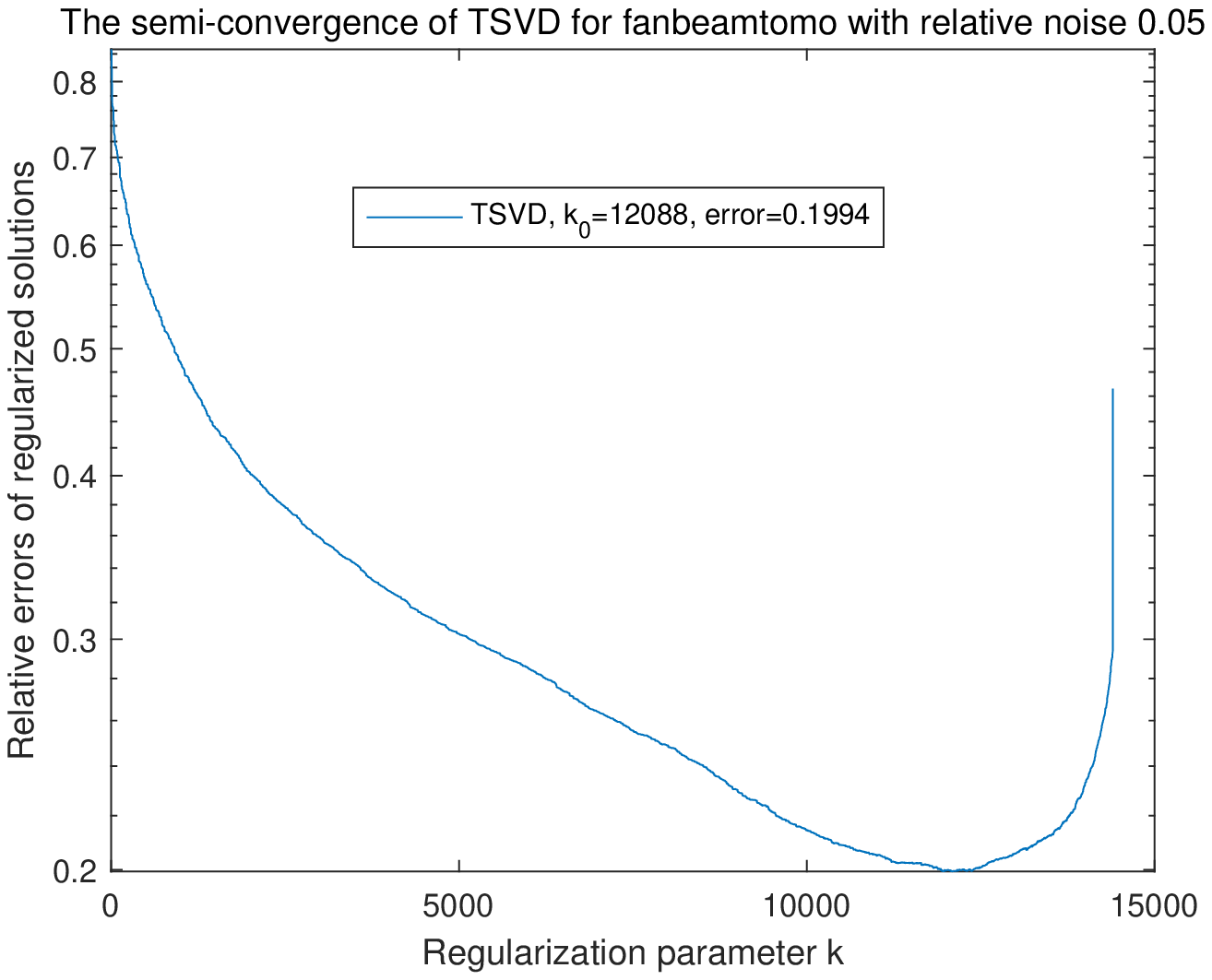}}
  \centerline{(d)}
\end{minipage}
\caption{{\sf fanbeamtomo} of $m=61200,\ n=14400$ with $\frac{\sigma_1}{\sigma_n}=2472$
and $\varepsilon=10^{-3}$ and 0.05.}
\label{fig7}
\end{figure}

\begin{figure}
\begin{minipage}{0.48\linewidth}
 \centerline{\includegraphics[width=6.0cm,height=4.5cm]{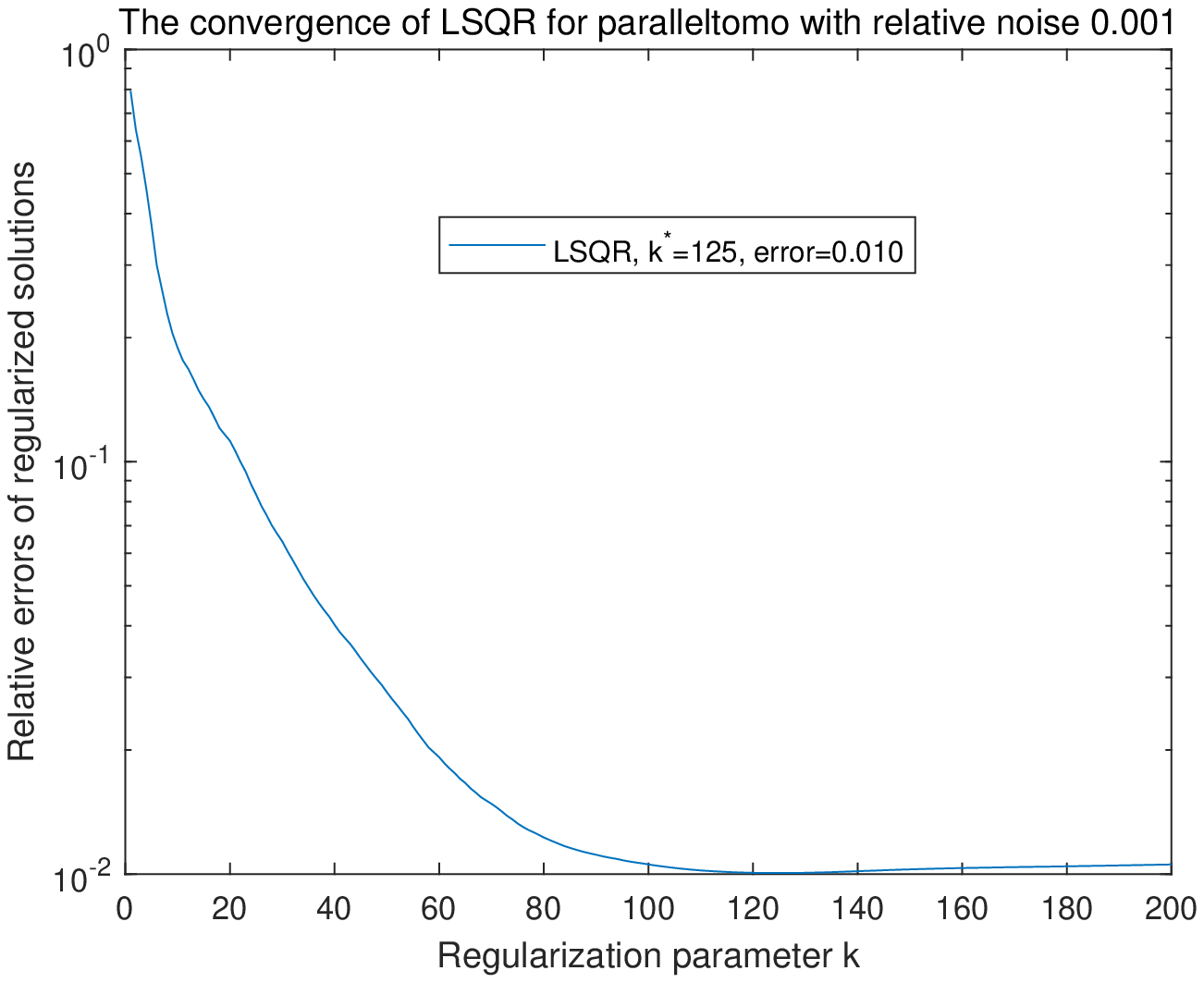}}
  \centerline{(a)}
\end{minipage}
\hfill
\begin{minipage}{0.48\linewidth}
  \centerline{\includegraphics[width=6.0cm,height=4.5cm]{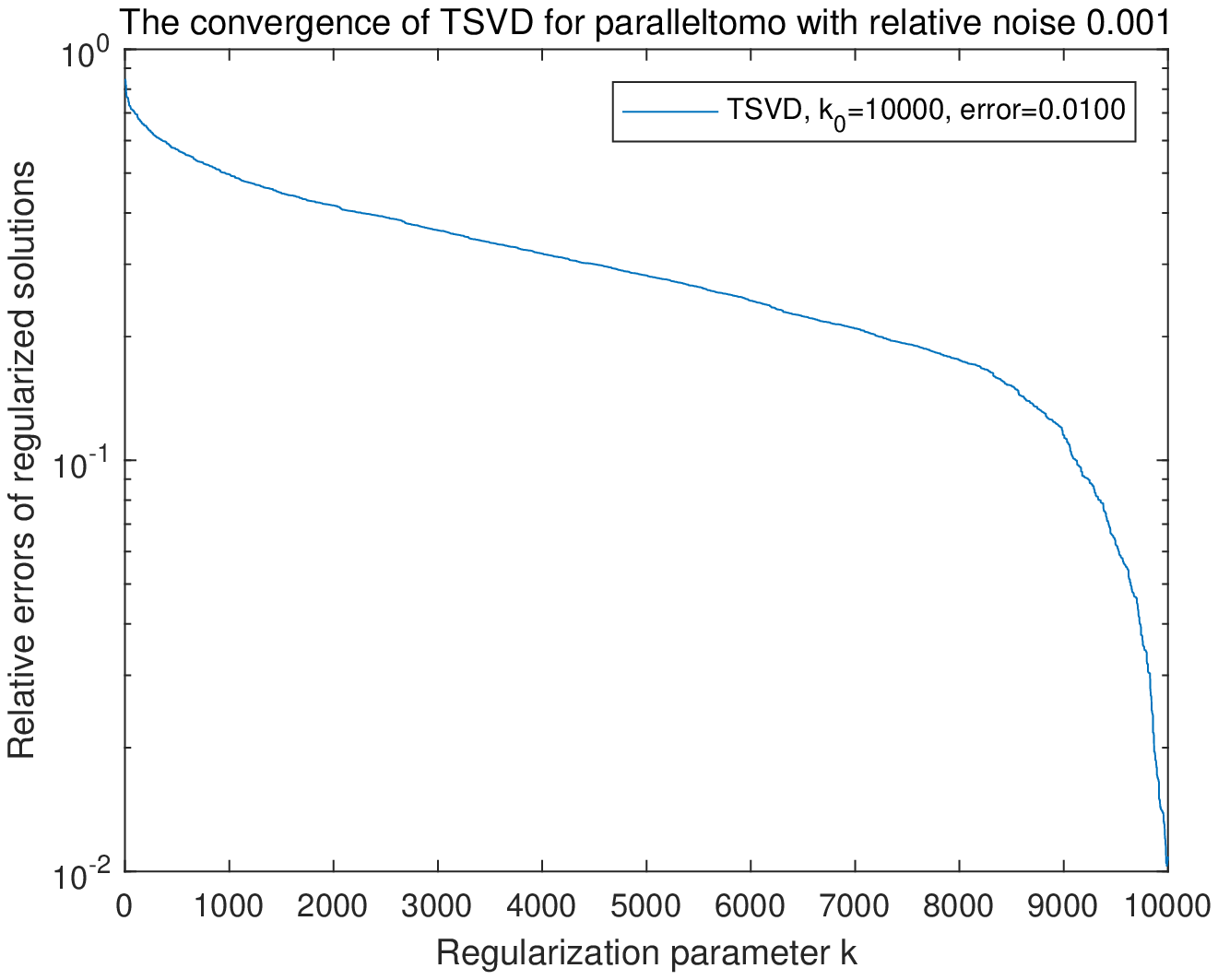}}
  \centerline{(b)}
\end{minipage}
\vfill
\begin{minipage}{0.48\linewidth}
 \centerline{\includegraphics[width=6.0cm,height=4.5cm]{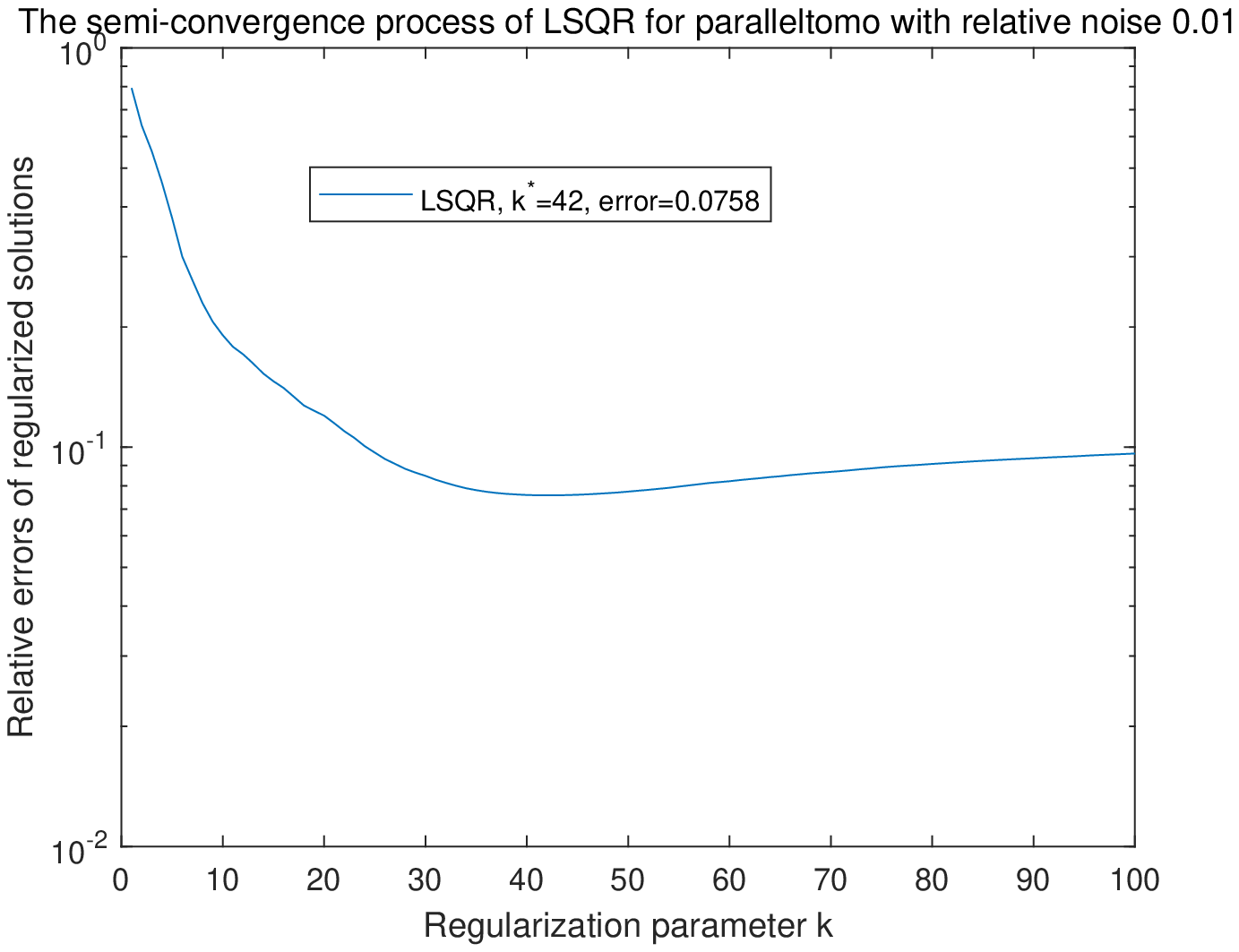}}
  \centerline{(c)}
\end{minipage}
\hfill
\begin{minipage}{0.48\linewidth}
  \centerline{\includegraphics[width=6.0cm,height=4.5cm]{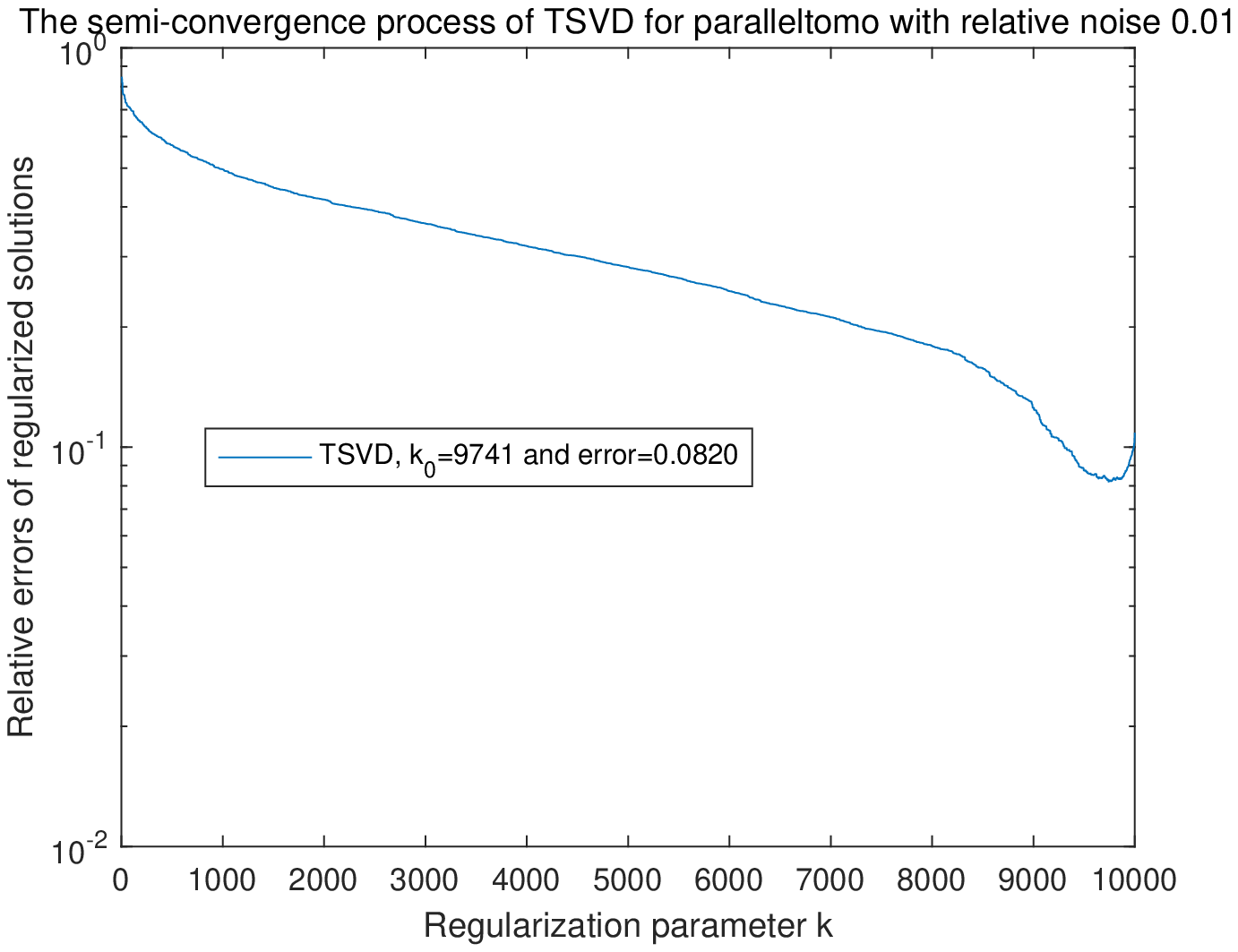}}
  \centerline{(d)}
\end{minipage}
\caption{{\sf paralleltomo} of $m=25380,\ n=10000$ with $\frac{\sigma_1}{\sigma_n}=408.8$
and $\varepsilon=10^{-3}$ and $10^{-2}$.} \label{fig8}
\end{figure}

\subsection{Decay properties of $\gamma_k$ and $\alpha_{k+1}+\beta_{k+2}$}

In Figures~\ref{fig1} (c)--\ref{fig5} (c), we depict
the decay curves of $\gamma_k$ and  $\alpha_{k+1}+\beta_{k+2}$.
From them, we see that $\alpha_{k+1}+\beta_{k+2}$ behaves very
similar to $\gamma_k$ and matches $\gamma_k$ very well, and their
decay curves highly resemble,
independent of the degree of ill-posedness. Therefore, we can use
$\alpha_{k+1}+\beta_{k+2}$ to reliably determine the decay rate
of $\gamma_k$. We also see that $\gamma_k$ monotonically decreases
with respect to $k$. These results justify Theorem~\ref{main2} and the first
two remarks followed.

Figure~\ref{figPRblurrotation} (c) shows that $\gamma_k$ monotonically
decays very slowly and $\gamma_1=1.1374$ and $\gamma_{50}=1.1227$, respectively.
It is due to the scale of vertical ordinate that, at first glance, gives one
an illusion that $\alpha_{k+1}+\beta_{k+2}$ deviates from $\gamma_k$ very much.
As a matter of fact, the $\alpha_{k+1}+\beta_{k+2}$
estimate the $\gamma_k$ quite accurately,
and the minimum and maximum of $\alpha_{k+1}+\beta_{k+2}$ are 0.9507 and 1.2576,
respectively. These confirm our theoretical results.

\subsection{Observations on the regularization ability of LSQR}

We now report the results on
the 2D problems {\sf blur}, {\sf fanbeamtomo} and {\sf paralleltomo}
in Table~\ref{tab1}. Particularly, we investigate the regularization
behavior and ability of LSQR on them and all the previous problems.

We give some details in Figures~\ref{fig6}--\ref{fig8}. Although
the orders $m$ and $n$ are already tens of thousands, the ratios
$\sigma_1/\sigma_n$ on these three problems are only 31.5, 2472 and 408.8,
and their singular values are far from small and are not yet clustered at zero,
which, intuitively, do not satisfy the definition of a discrete
ill-posed problem since the ratio $\sigma_1/\sigma_n$ is modest.
Therefore, the existing
regularization theory does not suit well for such practical problems.

Indeed, it is not prompt to
regard such problems as discrete ill-posed ones because,
in the context of solving least squares problems or linear systems,
such problems are
quite well conditioned and at least not ill conditioned at first glance.
Nevertheless, the situation is subtle, and we will have more findings.
With different relative noise levels $\varepsilon$, the TSVD method and
LSQR may exhibit very different behavior. For $\varepsilon=10^{-3}$,
we have observed that the best TSVD regularized solutions are
simply $x_{k_0}^{tsvd}=x_{n}^{tsvd}=A^{\dagger}b=x_{naive}$, as indicated
by Figure~\ref{fig6} (b)--Figure~\ref{fig8} (b). This means that $e$ does not
exert influence on regularization and we have solved the problems as if
they are ordinary ones.
LSQR treats them as ordinary ones and solves them in its regular way too; it
is seen from Figure~\ref{fig6} (a)--Figure~\ref{fig8} (a) that the solution
errors decrease until they stabilize and no semi-convergence occurs.

For {\sf blur}  and {\sf paralleltomo} with a larger $\varepsilon=10^{-2}$
and {\sf fanbeamtomo} with $\varepsilon=0.05$,
the situation changes. The noises $e$ critically affect
the solution processes, and both TSVD and LSQR exhibit semi-convergence
phenomena, as we
see from (c)-(d) in Figures~\ref{fig6}--\ref{fig8}. We have
observed that TSVD takes many $k_0$ SVD dominant components to form
the best regularized solutions but LSQR uses much fewer
$k^*$ iterations to obtain the best regularized solutions with the same
accuracy as TSVD does.

We also see from Figure~\ref{fig1} (d)--Figure~\ref{fig5} (d) that LSQR
takes $k^*\leq k_0$
iterations to compute the best regularized solutions $x_{k^*}^{lsqr}$ as
accurately as $x_{k_0}^{tsvd}$ for the five test problems
{\sf shaw, gravity, heat, deriv2} and {\sf regutm} that have different degrees
of ill-posedness; the weaker is the degree of ill-posedness, the smaller $k^*$ is
relative to $k_0$. For $k^*=k_0$, we have established the rigorous regularization
theory in this paper and \cite{jia18a}
and shown that LSQR has the full regularization; for $k^*<k_0$, the full or partial
regularization is not yet revealed theoretically.
However, beyond one's common expectation,
as we have seen from Figures~\ref{fig3} (d)--\ref{fig8} (d),
the experiments on all the other ill-posed problems show that
the best regularized solutions
$x_{k^*}^{lsqr}$ by LSQR are as accurate as the best solutions $x_{k_0}^{tsvd}$
by the TSVD method.
It is worthwhile to notice that both $\varepsilon=10^{-3}$ and $0.05$ are
practical. For these problems, the experiments have demonstrated that
LSQR has the full regularization. Numerical experiments on
a number of 2D mildly ill-posed image deblurring
problems and tomography problems from \cite{gazzola18}
have also demonstrated that LSQR has the full regularization \cite{jia18a}.

\section{Conclusions}\label{concl}

For the large-scale \eqref{eq1}, iterative solvers
are generally the only viable approaches. Of them, the mathematically
equivalent LSQR and CGLS are most
popularly used  Krylov iterative solvers for general purposes.
They have general regularizing effects and exhibit semi-convergence.
It has long been known that if the Ritz values converge to the large singular
values of $A$ in natural order until the occurrence of semi-convergence
of LSQR then the regularized solution $x_{k^*}^{lsqr}$ is as accurate as the TSVD
solution $x_{k_0}^{tsvd}$, that is, LSQR has the full
regularization.

For severely and moderately ill-posed problems, we have proved that,
with suitable $\rho>1$ or $\alpha>1$,
a $k$-step Lanczos bidiagonalization produces a near best
rank $k$ approximation of $A$ and the $k$ Ritz values approximate
the first $k$ large singular values of $A$ in natural order until
the semi-convergence of LSQR, so that
LSQR has the full regularization.
But for moderately ill-posed problems with $\alpha>1$ not enough
and mildly ill-posed problems, we have proved
that the above results generally do not hold for some $k\leq k^*$.
These results have given accurate and definitive solutions of
the highly concerned and challenging problems on
the convergence behavior of Ritz values for the three kinds
of ill-posed problems.

We have proved that the accuracy $\gamma_k$ of rank $k$ approximation generated
by Lanczos bidiagonalization monotonically increases with $k$.
We have also derived bounds for the diagonals and subdiagonals of bidiagonal
matrices generated by Lanczos bidiagonalization. Particularly,
we have proved that they decay as fast as the singular values of $A$
for severely or moderately ill-posed problems
with suitable $\rho>1$ or $\alpha>1$.
These bounds are of theoretical and practical importance since we
have shown that $\alpha_{k+1}+\beta_{k+2}$ can be used to reliably
judge the decay rate of $\gamma_k$ of the rank $k$ approximations
during computation without extra cost.

We have made illuminating numerical experiments and confirmed our
theory. In addition, we have investigated the regularizing effects of LSQR
and the TSVD method on some practical discrete problems arising from 2D image
continuous deblurring problems. The 2D test problems, though large scale,
are quite well conditioned. We have found that LSQR and the TSVD methods work as
if they solved ordinary linear systems, in which a noise $e$ with practical
level $\varepsilon$ may not play a role in regularization; if $e$ is larger, the
two methods have semi-convergence phenomena. In any case, the best regularized
solutions obtained by LSQR and the TSVD methods essentially have the same
accuracy, meaning that LSQR has the full regularization.

As the numerical experiments in this paper and \cite{jia18a}
have demonstrated, LSQR has the full regularization for all the test problems
in \cite{gazzola18,hansen07,hansen12}, independent of the degree of
ill-posedness. These draw us to the conjecture that LSQR
has the full regularization for any kind of ill-posed problem
with the discrete Picard condition satisfied.




\begin{thebibliography}{10}

\bibitem{aster}
R.~C. Aster, B.~Borchers, and C.~H. Thurber, \emph{Parameter {E}stimation and
  {I}nverse {P}roblems}, second ed., Elsevier, New York, 2013.

\bibitem{berisha}
S.~Berisha and J.~G. Nagy, \emph{Restore tools: Iterative methods for image
  restoration}, 2012, available from
  http://www.mathcs.emory.edu/$^\sim$nagy/RestoreTools.

\bibitem{bjorck96}
{\AA}.~Bj{\"{o}}rck, \emph{Numerical {M}ethods for {L}east {S}quares
  {P}roblems}, SIAM, Philadelphia, PA, 1996.

\bibitem{bjorck15}
\sameauthor, \emph{Numerical {M}ethods in {M}atrix {C}omputations}, Texts in
  Applied Mathematics, Springer, Cham, 2015.

\bibitem{bjorck79}
{\AA}.~Bj{\"{o}}rck and L.~Eld{\'{e}}n, \emph{Methods in numerical algebra for
  ill-posed problems}, Report LiTH-R-33-1979, Dept. of Mathematics,
  Link\"{o}ping Univeristy, Sweden, 1979.

\bibitem{eicke}
B.~Eicke, A.~K. Lious, and R.~Plato, \emph{The instability of some gradient
  methods for ill-posed problems}, Numer. Math., 58 (1990),
  pp.~129--134.

\bibitem{engl93}
H.~W. Engl, \emph{Regularization methods for the stable solution of inverse
  problems}, Surveys Math. Indust., 3 (1993), pp.~71--143.

\bibitem{engl00}
H.~W. Engl, M.~Hanke, and A.~Neubauer, \emph{{R}egularization of {I}nverse
  {P}roblems}, Kluwer Academic Publishers, 2000.

\bibitem{gazzola18}
S.~Gazzola, P.~C. Hansen, and J.~G. Nagy, \emph{{IR} tools: A {MATLAB} package
  of iterative regularization methods and large-scale test problems}, Numer.
  Algor., doi.org/10.1007/s11075-018-0570-7.

\bibitem{gazzola15}
S.~Gazzola and P.~Novati, \emph{Inheritance of the discrete {P}icard condition
  in {K}rylov subspace methods}, BIT Numer. Math., 56 (2016),
  pp.~893--918.

\bibitem{gilyazov}
S.~F. Gilyazov and N.~L. Gol'dman, \emph{Regularization of {I}ll-{P}osed
  {P}roblems by {I}teration {M}ethods}, Mathematics and its Applications,
  Kluwer Academic Publishers, Dordrecht, 2000.

\bibitem{hanke95}
M.~Hanke, \emph{Conjugate {gr}adient {T}ype {M}ethods for {I}ll-{P}osed
  {P}roblems}, Pitman Research Notes in Mathematics Series, Longman,
  Essex, 1995.

\bibitem{hanke01}
\sameauthor, \emph{On {L}anczos based methods for the regularization of discrete
  ill-posed problems}, BIT Numer. Math., 41 (2001), Suppl., pp.~1008--1018.

\bibitem{hanke93}
M.~Hanke and P.~C. Hansen, \emph{Regularization methods for large-scale
  problems}, Surveys Math. Indust., 3 (1993), pp.~253--315.

\bibitem{hansen90}
P.~C. Hansen, \emph{The discrete {P}icard condition for discrete ill-posed
  problems}, BIT, 30 (1990), pp.~658--672.

\bibitem{hansen90b}
\sameauthor, \emph{Truncated singular value decomposition solutions to discrete
  ill-posed problems with ill-determined numerical rank}, SIAM J. Sci. Statist.
  Comput., 11 (1990), pp.~503--518.

\bibitem{hansen95}
\sameauthor, \emph{Test matrices for regularization methods}, SIAM J. Sci. Comput.,
16 (1995), pp.~506--512.

\bibitem{hansen98}
\sameauthor, \emph{Rank-{D}eficient and {D}iscrete {I}ll-{P}osed {P}roblems:
  {N}umerical {A}spects of {L}inear {I}nversion}, SIAM Monographs on
  Mathematical Modeling and Computation, SIAM, Philadelphia, PA, 1998.

\bibitem{hansen07}
\sameauthor, \emph{Regularization {T}ools version 4.0 for {M}atlab 7.3}, Numer.
  Algor., 46 (2007), pp.~189--194.

\bibitem{hansen10}
\sameauthor, \emph{Discrete {I}nverse {P}roblems: {I}nsight and {A}lgorithms},
  Fundamentals of Algorithms, SIAM, Philadelphia, PA, 2010.

\bibitem{hansen12}
P.~C. Hansen and M.~Saxild-Hansen, \emph{{AIR} tools--a {MATLAB} package of
  algebraic iterative reconstruction methods}, J. Comput. Appl. Math.,
  236 (2012), pp.~2167--2178.

\bibitem{hps16}
M.~R. Hn\v{e}tynkov\'{a}, Marie Kub\'{i}nov\'{a}, and M.~Ple\v{s}inger,
  \emph{Noise representation in residuals of {LSQR}, {LSMR}, and {C}raig
  regularization}, Linear Algebra Appl., 533 (2017), pp.~357--379.

\bibitem{hps09}
M.~R. Hn\v{e}tynkov\'{a}, M.~Ple\v{s}inger, and Z.~Strako\v{s}, \emph{The
  regularizing effect of the {G}olub-{K}ahan iterative bidiagonalization and
  revealing the noise level in the data}, BIT Numer. Math., 49 (2009),
  pp.~669--696.

\bibitem{hofmann86}
B.~Hofmann, \emph{{R}egularization for {A}pplied {I}nverse and {I}ll-{P}osed
  {P}roblems}, Teubner, Stuttgart, Germany, 1986.

\bibitem{huangjia}
Y.~Huang and Z.~Jia, \emph{Some results on the regularization of {LSQR} for
  large-scale ill-posed problems}, Science China Math., 60 (2017), pp.~701--718.

\bibitem{ito15}
K.~Ito and B.~Jin, \emph{Inverse {P}roblems: Tikhonov {T}heory and
  {A}lgorithms}, Series on Applied Mathematics, World Scientific
  Publishing Co. Pte. Ltd., Hackensack, NJ, 2015.

\bibitem{jia18a}
Z.~Jia, \emph{{A}pproximation accuracy of the {K}rylov subspaces for linear
  discrete ill-posed problems},  (2018), arXiv:math.NA/1805.10132.

\bibitem{kaipio}
J.~Kaipio and E.~Somersalo, \emph{Statistical and {C}omputational {I}nverse
  {P}roblems}, Applied Mathematical Sciences, Springer-Verlag, New
  York, 2005.

\bibitem{kern}
M.~Kern, \emph{Numerical {M}ethods for {I}nverse {P}roblems}, John Wiley \&
  Sons, Inc., 2016.

\bibitem{kirsch}
A.~Kirsch, \emph{An {I}ntroduction to the {M}athematical {T}heory of {I}nverse
  {P}roblems}, second ed., Applied Mathematical Sciences, Springer,
  New York, 2011.

\bibitem{natterer}
F.~Natterer, \emph{The {M}athematics of {C}omputerized {T}omography}, Classics
  in Applied Mathematics, SIAM, Philadelphia, PA, 2001, Reprint of the
  1986 original edition.

\bibitem{paige82}
C.~C. Paige and M.~A. Saunders, \emph{{LSQR}: an algorithm for sparse linear
  equations and sparse least squares}, ACM Trans. Math. Software, 8
  (1982), pp.~43--71.

\bibitem{stewart01}
G.~W. Stewart, \emph{Matrix {A}lgorithms {II}: {E}igensystems}, SIAM,
  Philadelphia, PA, 2001.

\bibitem{stewartsun}
G.~W. Stewart and J.-G Sun, \emph{Matrix {P}erturbation {T}heory}, Computer
  Science and Scientific Computing, Academic Press, Inc., Boston, MA, 1990.

\bibitem{vogel02}
C.~R. Vogel, \emph{Computational {M}ethods for {I}nverse {P}roblems}, Frontiers
  in Applied Mathematics, SIAM, Philadelphia, PA, 2002.

\end{thebibliography}

\end{document}